\magnification 1200

\catcode`\@=11
\def\displaylinesno #1{\displ@y\halign{
\hbox to\displaywidth{$\@lign\hfil\displaystyle##\hfil$}&
\llap{$##$}\crcr#1\crcr}}

\def\ldisplaylinesno #1{\displ@y\halign{
\hbox to\displaywidth{$\@lign\hfil\displaystyle##\hfil$}&
\kern-\displaywidth\rlap{$##$}
\tabskip\displaywidth\crcr#1\crcr}}
\catcode`\@=12

\font \para=cmr10 at 14.4 pt
\let \ol=\overline
\let \ul=\underline
\let \ti=\widetilde

\font \nr=eufb7 at 10pt

\font \srm=cmr10 at 7.5pt

\font\ens=msbm10
\font\sens=msbm10 at 7.5pt
\font\main=cmsy10 at 10pt
\font\smain=cmsy10 at 7.5pt
\font\ssmain=cmsy10 at 5.625pt

\font \fin=lasy8 at 15.4 pt
\def \X{\mathop{\hbox{\nr X}^{\hbox{\srm nr}}}\nolimits}
\def \CC{\mathop{\hbox{\ens C}}\nolimits}
\def \sCC{\mathop{\hbox{\sens C}}\nolimits}
\def \RR{\mathop{\hbox{\ens R}}\nolimits}

\def \Xim{\mathop{\hbox{\nr X}^{\hbox{\srm nr}}_0}\nolimits}

\def \o{\mathop{\hbox{\main O}}\nolimits}
\def \so{\mathop{\hbox{\smain O}}\nolimits}
\def \sso{\mathop{\hbox{\ssmain O}}\nolimits}
\def \mes{\mathop{\hbox{mes}}\nolimits}
\def \C{\mathop{\hbox{\main C}}\nolimits}
\def \W{\mathop{\hbox{\main W}}\nolimits}

\def \A{\mathop{\hbox{\main A}}\nolimits}
\def \D{\mathop{\hbox{\main D}}\nolimits}
\def \H{\mathop{\hbox{\main H}}\nolimits}
\def \L{\mathop{\hbox{\main L}}\nolimits}

\def \P{\mathop{\hbox{\main P}}\nolimits}
\def \I{\mathop{\hbox{\main I}}\nolimits}

\def \R{\mathop{\hbox{\main R}}\nolimits}
\def \S{\mathop{\hbox{\main S}}\nolimits}
\def \sA{\mathop{\hbox{\smain A}}\nolimits}
\def \sD{\mathop{\hbox{\smain D}}\nolimits}

\def \E{\mathop{\hbox{\main E}}\nolimits}
\def \sE{\mathop{\hbox{\smain E}}\nolimits}
\def \sH{\mathop{\hbox{\smain H}}\nolimits}
\def \sI{\mathop{\hbox{\smain I}}\nolimits}

\def \sL{\mathop{\hbox{\smain L}}\nolimits}
\def \sP{\mathop{\hbox{\smain P}}\nolimits}

\def \sS{\mathop{\hbox{\smain S}}\nolimits}
\def \ssS{\mathop{\hbox{\ssmain S}}\nolimits}
\def \id{\mathop{\hbox{\rm id}}\nolimits}

\def \tr{\mathop{\hbox{\rm tr}}\nolimits}
\def \GL{\mathop{\hbox{\rm GL}}\nolimits}
\def \End{\mathop{\hbox{\rm End}}\nolimits}

\def \Hom{\mathop{\hbox{\rm Hom}}\nolimits}
\def \Rat{\mathop{\hbox{\rm Rat}}\nolimits}
\def \sRat{\mathop{\hbox{\srm Rat}}\nolimits}
\def \Res{\mathop{\hbox{\rm Res}}\nolimits}
\def \rg{\mathop{\hbox{\rm rg}}\nolimits}

\def \supp{\mathop{\hbox{\rm supp}}\nolimits}

\def \St{\mathop{\hbox{\rm St}}\nolimits}
\def \GL{\mathop{\hbox{\rm GL}}\nolimits}

\def \deg{\mathop{\hbox{\rm deg}}\nolimits}
\def \sreg{\mathop{\hbox{\srm reg}}\nolimits}
\def \Stab{\mathop{\hbox{\rm Stab}}\nolimits}
\def \sStab{\mathop{\hbox{\srm Stab}}\nolimits}
\centerline {\para D\'ecomposition spectrale et repr\'esentations sp\'eciales}
\centerline {\para d'un groupe r\'eductif p-adique}

\null
\null
\centerline {par Volker Heiermann}

\vfill
Institut f\"ur Mathematik

Humboldt-Universit\"at zu Berlin

10099 Berlin

Germany

heierman@mathematik.hu-berlin.de
\eject
{\bf R\'esum\'e:} \it Soit $G$ un groupe r\'eductif p-adique connexe. Nous effectuons une
d\'ecompo-sition spectrale sur $G$ \`a partir de la formule d'inversion de Fourier utilis\'ee dans \rm
[\it Une formule de Plancherel pour l'alg\`ebre de Hecke d'un groupe r\'eductif p-adique - \rm V.
Heiermann, Comm. Math. Helv. {\bf 76}, 388-415, 2001]. \it Nous en d\'eduisons essentiellement qu'une
repr\'esentation cuspidale d'un sous-groupe de Levi $M$ appartient au support cuspidal d'une
repr\'esentation de carr\'e int\'egrable de $G$ si et seulement si c'est un p\^ole de la fonction $\mu
$ de Harish-Chandra d'ordre \'egal au rang parabolique de $M$. Ces p\^oles sont d'ordre maximal. Plus
pr\'ecis\'ement, nous montrons que cette condition est n\'ecessaire et que sa suffisance \'equivaut \`a
une propri\'et\'e combinatoire de la fonction $\mu $ de Harish-Chandra qui s'av\`ere \^etre une 
cons\'equence d'un r\'esultat de E. Opdam. En outre, nous obtenons des identit\'es entre des combinaisons 
lin\'eaires de coefficients matriciels. Ces identit\'es contiennent des informations sur le degr\'e formel 
des repr\'esentations de carr\'e int\'egrable ainsi que sur leur position dans la repr\'esentation induite.
\rm

\null {\bf Summary:} \it Let $G$ be a reductiv connected $p$-adic group. With help of the Fourier
inversion formula used in \rm [\it Une formule de Plancherel pour l'alg\`ebre de Hecke d'un groupe
r\'eductif p-adique - \rm V. Heiermann, Comm. Math. Helv. {\bf 76}, 388-415, 2001] \it we give a
spectral decomposition on $G$. In particular we deduce from it essentially that a cuspidal
representation of a Levi subgroup $M$ is in the cuspidal support of a square integrable representation
of $G$, if and only if it is a pole of Harish-Chandra's $\mu $-function of order equal to the parabolic
rank of $M$. These poles are of maximal order. In more explicit terms, we show that this condition is
necessary and that its sufficiency is equivalent to a combinatorical property of Harish-Chandra's $\mu
$-function which appears to be a consequence of a result of E. Opdam. We get also identities between some 
linear combinations of matrix coefficients. These identities contain informations on the formel degree of
square integrable representations and on their position in the induced representation. \rm

\vfill
\eject 
{\para Table des mati\`eres}

\vskip 3cm 
{\bf Introduction}\hfill { 4}

\null {\bf 1. Notations et pr\'eliminaires} \hfill { 8}

\null {\bf 2. Formule de Plancherel et th\'eor\`eme de Paley-Wiener}\hfill { 13}

\null {\bf 3. Th\'eor\`eme des r\'esidus}\hfill { 15}

\null {\bf 4. La fonction $\mu $ de Harish-Chandra, son syst\`eme de}\hfill { 25}

{\bf \ \ \ \ racines singuli\`eres et son groupe de Weyl}

\null {\bf 5. Les d\'eriv\'ees d'un coefficient matriciel}\hfill { 31}

\null {\bf 6. Prolongement rationnel d'une identit\'e polynomiale}\hfill { 34}

\null {\bf 7. L'\'el\'ement $\varphi _{P,\so }^H$ de l'espace de Paley-Wiener}\hfill { 40}

\null {\bf 8. La d\'ecomposition spectrale}\hfill { 42}

\null {\bf Annexe A: Le cas d'un groupe semi-simple de type $G_2$}\hfill { 54}

\null {\bf Annexe B: Une preuve du th\'eor\`eme de Opdam}\hfill 62

\null
{\bf Index de notations}\hfill 69

\null
{\bf Index terminologique}\hfill 71

\null 
{\bf R\'ef\'erences bibliographiques}\hfill {72}

\vfill
\eject
\null
\null
\null
{\bf Introduction:} 

\null
\null  Soit $G$ le groupe des points rationnels d'un groupe alg\'ebrique, r\'eductif et connexe
d\'efini sur un corps local non archim\'edien $F$. Dans la terminologie de [H1] et [W] (qui est
rappel\'e dans les paragraphes 1 et 2) on se fixe un couple $(P, \o)$
form\'e d'un sous-groupe parabolique standard $P=MU$ et de l'orbite inertielle $\o $ d'une
repr\'esentation irr\'eductible cuspidale de $M$. Soit $\varphi _{P,\so }$ la composante en $(P,\o )$
d'un \'el\'ement de l'espace de Paley-Wiener matriciel de $G$. On a montr\'e dans [H1] (cf. proposition
0.2) qu'il existe une application polynomiale $\xi _{\varphi _{P, \sso}}:\o \rightarrow i_{\overline
{P}\cap K}^KE_{\so }\otimes i_{P\cap K}^K E_{\so }^{\vee }$ \`a image dans un espace de dimension
finie, telle que
$$\varphi _{P,\so }(\sigma )=\sum _{w\in W; w\so =\so } (J_{P\vert
\overline {wP}}(\sigma )\lambda (w)\otimes J_{P\vert wP}(\sigma  ^{\vee })\lambda (w))\ \xi _{\varphi
_{P,\sso }}(w^{-1}\sigma ),\eqno {\hfill (\#)}$$  en tout point $\sigma $ de $\o $ en lequel les
op\'erateurs d'entrelacement
$J_{P\vert\overline {wP}}(\sigma )\lambda (w)$ et $J_{P\vert wP}(\sigma ^{\vee })$ $\lambda (w)$ sont
d\'efinis. En outre, la fonction $f_{\varphi _{P,\sso }}$ d\'efinie sur $G$ par
$$f_{\varphi _{P,\sso }}(g)=\gamma(G/M)\int _{\Re (\sigma )=r\gg _P0}\ \deg(\sigma )E_P^G((J_{P\vert\ol{P}}(\sigma )\otimes 1)
\xi _{\varphi _{P, \sso }}(\sigma ))(g^{-1}) \mu (\sigma )d\Im (\sigma ),\eqno {\hfill (\#\#)}$$  
est lisse \`a support compact et sa transform\'ee de Fourier en $(P,\o )$ est donn\'ee par $\varphi _{P,\so }$ (la fonction $\mu $ \'etant celle de Harish-Chandra).

Dans cet article nous allons calculer l'int\'egrale (\#\#) \`a l'aide de la formule int\'egrale de
Cauchy et l'\'ecrire comme somme d'int\'egrales, chaque int\'egrale portant sur l'orbite unitaire d'une
repr\'esentation irr\'eductible de carr\'e int\'egrable dont le support cuspidal est donn\'e par un
\'el\'ement de la classe de conjugaison de $\o $. Les termes qui apparaissent seront identifi\'es avec
ceux dans la formule de Plancherel.

On montrera \`a cette occasion que seuls les p\^oles venant de la fonction $\mu $ de Harish-Chandra
subsisteront apr\`es avoir effectu\'e le changement de contours. On en d\'eduira en particulier que,
pour qu'une repr\'esentation irr\'eductible cuspidale de $M$ appartienne au support cuspidal d'une
repr\'esentation de carr\'e int\'egrable de $G$, il faut qu'elle soit un p\^ole   de la fonction $\mu $
de Harish-Chandra d'ordre \'egal au rang parabolique de $M$ (cf. corollaire {\bf 8.6}). Ces p\^oles
sont d'ordre maximal (ceci r\'esulte d'un r\'esultat dans [O2] cf. {\bf 8.1}). Par ailleurs, nous
montrons que la r\'eciproque est \'equivalente \`a une propri\'et\'e combinatoire de la fonction
$\mu $ de Harish-Chandra qui s'av\`ere \^etre une cons\'equence d'un r\'esultat  de E. Opdam (cf. {\bf
8.7} pour un \'enonc\'e plus pr\'ecis). Le tout r\'epond \`a une conjecture qui nous a \'et\'e
communiqu\'e par A. Silberger et qui a \'et\'e formul\'e par lui dans un projet de recherche en 1978.

Dans le cas de la s\'erie principale non ramifi\'ee d'un groupe d\'eploy\'e, o\`u la localisation des
p\^oles et z\'eros de la fonction $\mu $ de Harish-Chandra est connue, cette conjecture sur le support
cuspidal d'une s\'erie discr\`ete non cuspidale est \'equivalente \`a celle de Deligne-Langlands (cf.
[H2]). Dans ce cas le r\'esultat est d\'ej\`a connu gr\^ace aux travaux de Kazhdan et Lusztig (au moins
si $G$ est semi-simple \`a centre connexe) (cf. [KL]).

Signalons par ailleurs les travaux de Bernstein-Zelevinsky [Z] qui d\'ecrivent les repr\'e-sentations
irr\'eductibles de carr\'e int\'egrable non cuspidales de $GL_N (F)$ et ceux de Moeglin et
Moeglin-Tadic (cf. [M2] et [MT]) qui traitent - modulo une conjecture sur les p\^oles de la fonction
$\mu $ de Harish-Chandra (\'enonc\'e et justifi\'e dans [M1]) - les autres groupes classiques
d\'eploy\'es connexes (et un petit peu plus).

Nous obtenons en outre des identit\'es entre des combinaisons lin\'eaires de d\'eriv\'ees de
coefficients matriciels de repr\'esentations induites - d\'eduits d'op\'erateurs d'entrelacement - et
des combinaisons lin\'eaires de coefficients matriciels de repr\'esentations irr\'eductibles de carr\'e
int\'egrable de m\^eme support cuspidal. Ces identit\'es contiennent des informations sur le degr\'e
formel et sur la position des repr\'esentations de carr\'e int\'egrable dans la repr\'esentation
induite. On explique ceci dans {\bf 8.6} pour le cas "r\'egulier". Pour extraire ces informations dans
le cas g\'en\'eral, il faudrait analyser de plus pr\`es les combinaisons lin\'eaires de coefficients
matriciels qui apparaissent, ce qui semble \^etre un probl\`eme difficile. A titre d'exemple, nous
traitons le cas d'un groupe semi-simple d\'eploy\'e de type $G_2$ dans l'annexe. C'est le groupe de
rang minimal, o\`u des probl\`emes non n\'egligeables apparaissent. Notre \'etude utilise des arguments
de [MW1] Annexe III. (Elle est toutefois plus proche du cas des corps des fonctions.) Nous ne poussons
cependant pas l'\'etude de notre identit\'e dans le cas difficile aussi loin que cela a \'et\'e fait
dans [MW1], puisque d'un c\^ot\'e les r\'esultats correspondants sont d\'ej\`a connus dans ce cas (il
s'agit de repr\'esentations ayant un vecteur invariant  par un sous-groupe d'Iwahori) et d'autre part
que les calculs peuvent se faire dans le cas local \`a l'aide d'un logiciel. Il nous a donc sembl\'e
que cela ne repr\'esente pas beaucoup d'int\'er\^et d'exposer un calcul qui ne comporte aucune id\'ee
susceptible de se g\'en\'eraliser.

Le principe de ce travail est en quelque sorte analogue \`a la d\'ecomposition spectrale de l'espace
des formes automorphes sur un corps global qui a son origine dans les travaux de R.P. Langlands (cf.
[L], [MW1]), l'id\'ee de d\'ecrire des s\'eries discr\`etes par des r\'esidus ayant apparamment apparu
la premi\`ere fois dans des travaux de A. Selberg.  Signalons \'egalement le travail de E. Opdam [O2]
qui effectue une d\'ecomposition spectrale dans le contexte des alg\`ebres de Hecke affines \`a partir
de sa formule g\'en\'eratrice pour la trace [O1]. La partie des r\'esultats de [O2] qui se traduit en 
termes de la th\'eorie des repr\'esentations de groupes $p$-adiques s'applique donc en particulier \`a 
la s\'erie principale non ramifi\'ee d'un groupe $p$-adique semi-simple, en consid\'erant les 
\'el\'ements de l'alg\`ebre d'Iwahori-Hecke. L'auteur de [O2] en d\'eduit des informations sur le degr\'e 
formel de ces repr\'esentations (voir \'egalement [HO2]). On abordera cette question dans le cadre donn\'e 
par l'article pr\'esent \'eventuellement dans un travail futur.

Remarquons que l'analogue de nos r\'esultats dans le cadre de la d\'ecomposition spectrale de Langlands
est toujours un probl\`eme ouvert, des r\'esultats complets n'existant que pour le groupe $G=GL_N(F)$
(cf. [MW2]).

Survolons alors le contenu de cet article: au paragraphe 1 sont fix\'ees les notations. On y
trouve \'egalement les d\'efinitions des diff\'erentes notions attach\'ees aux hyperplans d'une orbite
inertielle $\o $ qui seront centrales dans le reste de ce travail. Le paragraphe 2 sert \`a des rappels
des r\'esultats de [H1] et sur la formule de Plancherel de Harish-Chandra, adapt\'es \`a nos besoins. 

Nous pr\'esentons au paragraphe 3  le th\'eor\`eme des r\'esidus que l'on appliquera dans la suite. Sa
pr\'esentation et sa preuve suivent la m\'ethode employ\'ee dans [BS]. Cet article traitant
d'int\'egrales sur des espaces vectoriels et non pas sur un produit de tores de $\CC ^*$, notre concept
de r\'esidu ne peut toutefois \^etre le m\^eme. Par ailleurs, les p\^oles que nous rencontrons ne sont
pas toujours r\'eels. Nous donnons donc une preuve assez d\'etaill\'ee. Remarquons que l'article [BS]
utilise une nouvelle approche inspir\'ee de [HO1] et cite [A], [MW1] et [L], la preuve de l'unicit\'e 
des donn\'ees de r\'esidu \'etant plus indirecte dans le travail de Langlands (cf. [MW1] V.2.2 Remarques 
(2) et V.3.13 Remarques (ii)). Mentionnons \'egalement le travail [Mo].

Dans le paragraphe 4, nous analysons de plus pr\`es le cadre dans lequel le th\'eor\`eme des r\'esidus
sera appliqu\'e ensuite. Nous donnons en particulier plus de d\'etails sur les hyperplans r\'esiduels
qui apparaissent. Gr\^ace \`a une observation que nous avons trouv\'ee dans un article de A. Silberger
[S3], on peut associer \`a ces hyperplans un syst\`eme de racines. En faisant usage des propri\'et\'es
des syst\`emes de racines, nous r\'e\'ecrivons le th\'eor\`eme de r\'esidus sous une forme qui
s'appliquera directement \`a notre situation. Le proc\'ed\'e et l'\'enonc\'e sont analogues au
r\'esultat principal de [BS]. En passant, on \'enoncera quelques lemmes combinatoires sur le groupe de
Weyl.

Au paragraphe 5, nous prouvons que, pour $D=D_{\lambda }$ un op\'erateur diff\'erentiel holomorphe sur
$a_{M,\sCC}^*$ \`a coefficients constants et $\tau $ un \'el\'ement de $i_P^GE_{\so}\otimes i_P^GE_{\so
}^{\vee }$, la fonction
$D_{\lambda }E_{P,\sigma _{\lambda }}^G(\tau )$, d\'efinie sur $G$, correspond au coefficient
matriciel d'une repr\'esentation admissible de longueur finie dont les sous-quotients irr\'eductibles
sont \`a support cuspidal dans la classe de $W$-conjugaison de $\o $. Ce sont de tels termes qui
interviennent lors de la d\'esint\'egration de la formule (\#\#).

Dans le paragraphe 6, nous g\'en\'eralisons une propri\'et\'e (\`a cet endroit  conjecturale) de la
transformation (\#) \`a certaines applications rationnelles $\ti {\xi }$ \`a la place de $\xi $. Ceci
nous permettra d'identifier par r\'ecurrence la partie continue qui appara\^\i t apr\`es
d\'esint\'egration de (\#\#).

Pour montrer que certains p\^oles que nous obtenons apr\`es le calcul de (\#\#) s'annulent, nous
utilisons des \'el\'ements sp\'eciaux de l'espace de Paley-Wiener de $G$. Une id\'ee semblable a
\'et\'e utilis\'ee dans [O2] pour se d\'ebarasser des r\'esidus inutiles. Les \'el\'ements que nous
utilisons sont d\'ej\`a apparus dans [H1] lors de la preuve de la proposition 0.2. Un rappel de la
d\'efinition de ces \'el\'ements ainsi que la preuve de leurs propri\'et\'es particuli\`eres est le
sujet du paragraphe 7. 

L'\'enonc\'e et la preuve de notre th\'eor\`eme principal font l'objet du paragraphe 8. Nous y
pr\'ecisons \'egalement la notion de l'ordre d'un p\^ole de la fonction $\mu $. Les r\'esultats
concernant les supports cuspidaux des repr\'esentations de carr\'e int\'egrable d\'ecouleront alors de
la d\'ecomposition spectrale et de sa preuve.

On trouvera dans l'annexe {\bf A} le traitement d'un groupe semi-simple d\'eploy\'e de type $G_2$ "\`a la
main". Ceci est l'analogue local du cas \'etudi\'e dans [MW1] A.III pour les corps de nombres comme
cela a \'et\'e remarqu\'e ci-dessus. 

A la fin, dans l'annexe {\bf B}, on donne pour la commodit\'e du lecteur une preuve du th\'eor\`eme 
{\bf 8.7} qui est, comme d\'ej\`a signal\'e, d\^u \`a E. Opdam (cf. [O2] theorem 3.29). (On l'a 
ins\'er\'e ici avec l'autorisation de E. Opdam, parce que la preuve dans l'article [O2] est assez 
concise et r\'edig\'ee dans un langage diff\'erent du n\^otre.)

Nous recommandons au lecteur de continuer la lecture de cet article directement avec le paragraphe 8
(apr\`es avoir \'eventuellement jet\'e un coup d'oeil sur les notations pr\'esent\'ees au paragraphe 1
et sur une revue des r\'esultats de [H1] et de la formule de Plancherel au paragraphe 2), en
consultant les paragraphes ant\'erieurs au fur et \`a mesure. 
 
Lors des diff\'erentes \'etapes de ce travail l'auteur a b\'en\'efici\'e d'un s\'ejour \`a l'Institute
for Advanced Study in Princeton en correspondance avec R.P. Langlands (financ\'e par une bourse
Feodor-Lynen de la fondation Alex. von Humboldt et par la bourse DMS 97-29992 de la NSF) et d'un
s\'ejour \`a l'Universit\'e Paris VII en correspondance avec M.-F. Vign\'eras (financ\'e par le
Nachkontaktprogramm de la fondation Alex. von Humboldt). Mon tuteur aupr\`es de la fondation Alex. von
Humboldt \'etait E.-W. Zink.

A ces occasions et lors d'autres s\'ejours \`a Paris, l'auteur a b\'en\'efici\'e de quelques
discussions avec R.P. Langlands, C. Moeglin et J.-L. Waldspurger, ce qui l'a aid\'e \`a guider ses
d\'emarches. Il les remercie vivement de lui avoir consacr\'e du temps.

Je remercie \'egalement P. G\'erardin et les autres membres de l'\'equipe de Th\'eorie de Groupes pour
le bon accueil en leur sein lors de mes diff\'erents s\'ejours \`a Paris, ainsi que L. Clozel pour des
encouragements.

Mes remerciements vont de m\^eme \`a W. Hoffmann pour une remarque sur les d\'eriv\'ees de coefficients
matriciels et s\^urtout \`a E. Opdam qui m'a consacr\'e du temps en m'expliquant son travail et qui a
autoris\'e l'annexe {\bf B}. 

Ma reconnaissance va finalement \`a A. Silberger pour m'avoir communiqu\'e la conjecture sur les
supports cuspidaux des repr\'esentations sp\'eciales \`a un moment o\`u j'\'etais ignorant et qui
m'avait motiv\'e \`a aller au bout de cette d\'ecomposition spectrale. 

Pour l'examen d\'etaill\'e de mon travail, je remercie tout particuli\`erement J.-L. Waldspurger.

\null
\null
\null
\vfill
\eject
{\bf 1. Notations et pr\'eliminaires}

\null
\null {\bf 1.1} Fixons un tore d\'eploy\'e maximal $T_0$ dans $G$ et un sous-groupe ouvert compact
maximal $K$ de $G$ qui est en bonne position relative \`a $T_0$. (De tels groupes s'obtiennent comme
fixateur d'un point sp\'ecial de l'appartement associ\'e \`a $T_0$ dans l'immeuble de
$G$.) On notera $M_0$ le centralisateur de $T_0$ dans $G$ et $W=W^G=\hbox{\rm Norm}_G(T_0)/M_0$ le
groupe de Weyl d\'efini relatif \`a $T_0$.

Un sous-groupe parabolique $P$ de $G$ sera dit semi-standard, s'il contient $T_0$. Il existe alors un
unique sous-groupe de Levi $M$ de $P$ qui contient $T_0$. On appellera un tel sous-groupe $M$ de $G$
dans la suite simplement un sous-groupe de Levi semi-standard de $G$. Remarquons que
$M_0$ est le plus petit sous-groupe de Levi semi-standard. On \'ecrira souvent $P=MU$, si $P$ est un
sous-groupe parabolique semi-standard et si
$M$ est le sous-groupe de Levi semi-standard de $P$. On d\'esignera par $T_M$ le tore d\'eploy\'e
maximal dans le centre de $M$. Le sous-groupe parabolique de $G$ qui est oppos\'e \`a $P$ sera not\'e
$\ol {P}$. (On a donc
$P\cap\ol{P}=M$.)

Si $M$ est un sous-groupe de Levi semi-standard de $G$, on pose $\W (M)=\{w\in W\vert wMw^{-1}=M\}$, et
on \'ecrira
$W(M)$ pour le groupe quotient $\W(M)/W^M$. On identifiera ses \'el\'ements avec certains \'el\'ements
de $\W (M)$ (par exemple de longueur minimale dans leur classe \`a gauche modulo $W^M$). Les
identit\'es qui suivent seront toutefois essentiellement ind\'ependantes du choix des repr\'esentants.

On notera $\Sigma ^G(T_M)$, (resp. $\Sigma ^G(P)$ (resp. $\Sigma _{red}^G(P)$)) l'ensemble des racines
(resp. racines r\'eduites) de $T_M$ dans l'alg\`ebre de Lie de $G$ (resp. $U$). On omettra souvent
l'indice $G$ si cela ne pr\^ete pas \`a confusion.

\null {\bf 1.2} L'ensemble des caract\`eres $F$-rationnels d'un sous-groupe de Levi semi-standard de
$G$ sera not\'e
$\Rat (M)$. On pose
$a_M^*=\Rat(M)\otimes _{\hbox{\sens Z}}\hbox {\ens R}$ et on notera $a_M$ l'espace dual. Le choix d'un
sous-groupe parabolique semi-standard $P$ de Levi $M$ de $G$ est \'equivalent \`a celui d'un certain
ordre sur $a_M^*$. On notera
$>_P$ l'ordre qui correspond \`a $P$.

On d\'esignera par $\langle\ ,\ \rangle$ le produit de dualit\'e $a_M\times a_M^* \rightarrow
\hbox{\ens R}$. Avec
$q$ \'egal au cardinal du corps r\'esiduel de $F$, on d\'efinit un homomorphisme $H_M: M\rightarrow
a_M^*$ par la relation $\langle \chi\otimes 1, H_M(m) \rangle =-\log_q(\chi(m))$ pour tout
$\chi\in\Rat(M)$ et tout $m\in M$. 

Lorsque $M=M_0$, on \'ecrit $a_0^*$ \`a la place de $a_M^*$. On munit l'espace $a_0^*$ d'un produit
scalaire invariant pour l'action par $W$. Pour tout sous-groupe de Levi semi-standard $M$, l'espace
$a_M^*$ sera muni de la structure euclidienne induite par celle de $a_0^*$. Si $M'$ est un sous-groupe
de Levi qui contient $M$, on a une d\'ecomposition orthogonale $a_M^*=a_{M'}^*\oplus a_M^{M'*}$. On
notera parfois $\lambda =\lambda _{M'}+\lambda ^{M'}$ la d\'ecomposition d'un \'el\'ement $\lambda $ de
$a_M^*$ selon cette d\'ecomposition. On a une d\'ecomposition duale
$a_M=a_{M'}\oplus a_M^{M'}$. Le compos\'e de $H_M$ avec la projection de $a_M$ sur $a_M^{M'}$ sera
not\'e $H_M^{M'}$.
 
Pour $H$ un groupe, posons $X(H)=\Hom(H,\hbox{\ens C}^{\times })$. Notons $\vert . \vert _F$ la valeur
absolue normalis\'ee de $F$ dont le groupe des valeurs est $q^{\hbox{\sens Z}}$. Pour $M$ un
sous-groupe de Levi semi-standard, posons $M^1=\bigcap _{\chi\in\sRat (M)}$ $\ker (\vert\chi
\vert _F)$. Le groupe quotient $M/M^1$ est un $\hbox{\ens Z}$-module libre de rang \'egal \`a celui de
$T_M$. On note
$\X (M)=X(M/M^1)$. C'est une vari\'et\'e alg\'ebrique complexe affine dont les \'el\'ements sont
appel\'es les caract\`eres non ramifi\'es de $M$. On a une application surjective du complexifi\'e
$a_{M,\sCC}^*$ de $a_M^*$ dans
$\X(M)$ qui associe \`a $\lambda $ le caract\`ere non ramifi\'e $\chi _{\lambda }$ de
$M$, d\'efini par $\chi _{\lambda }(m)=q^{-\langle H_M(m),\lambda \rangle}$. Sa restriction \`a $a_M^*$
est injective, ce qui permet de d\'efinir la \it partie r\'eelle \rm d'un caract\`ere non ramifi\'e, en
posant $\Re (\chi _{\lambda }):=\Re(\lambda )$ pour tout $\lambda \in a_{M,\sCC}^*$.

Le sous-groupe de $\X (M)$ form\'e des caract\`eres unitaires sera not\'e $\Xim (M)$.

Un ensemble ordonn\'e $\underline {m}=(m_1,\cdots ,m_d)$ d'\'el\'ements de $M$ dont les images modulo
$M^1$ forment une base de ce $\hbox {\ens Z}$-module sera appel\'e une base de $M$ modulo $M^1$. Le
choix d'une telle base d\'etermine un isomorphisme $\X (M)\rightarrow \CC ^{\times }$, $\chi \mapsto
(\chi (m_1),..., \chi (m_d))$, not\'e $\eta _{\underline {m}}$ ci-dessous. Cet isomorphisme est
compatible avec la structure de vari\'et\'e affine complexe sur $\X (M)$. On pourra alors pr\'eciser la
notion de fonction polynomiale et rationnelle sur $\X (M)$: une fonction $f:\X (M)\rightarrow\CC$ sera
dite polynomiale, si et seulement s'il existe un polyn\^ome $p\in\CC [x_1,x_1^{-1},...,x_d,x_d^{-1}]$,
tel que $f(\chi )=p(\eta _{\underline {m}}(\chi ))$ pour tout $\chi \in \X(M)$. Le polyn\^ome $p$
d\'epend \'evidemment du choix de $\underline {m}$. Le cas d'une fonction rationnelle se traite de
mani\`ere analogue.

Les sous-groupes de Levi minimaux de $G$ contenant $M$ sont en bijection avec les \'el\'ements de
$\Sigma _{red}(P)$. On notera $M_{\alpha }$ le sous-groupe de Levi semi-standard de $G$ qui correspond
\`a un \'el\'ement $\alpha $ de
$\Sigma _{red}(P)$. On a une suite exacte courte
$$1\rightarrow M_{\alpha }^1\cap M/M^1\rightarrow M/M^1\rightarrow M_{\alpha }/M_{\alpha
}^1\rightarrow 1$$ qui prouve que $M_{\alpha }^1\cap M/M^1$ est un $\hbox{\ens Z}$-module libre de rang
$1$. Notons $h_{\alpha }$ un
\'el\'ement de $M_{\alpha }^1\cap M$ qui v\'erifie $\langle\alpha ,H_M(h_{\alpha })\rangle>0$ et dont
l'image dans
$M_{\alpha }^1\cap M/M^1$ d\'efinit une base de ce module. Notons par ailleurs $\alpha ^{\vee }$
l'unique \'el\'ement de $a_M$ qui est un multiple de $H_M(h_{\alpha })$ par un r\'eel $>0$ et qui
v\'erifie $\langle\alpha ,\alpha ^{\vee }\rangle =2$, ainsi que $\alpha ^*$ l'unique multiple de
$\alpha $ par un r\'eel $>0$ qui v\'erifie $\langle H_M(h_{\alpha }), \alpha ^*\rangle =1$.

On a une suite exacte duale de la pr\'ec\'edente donn\'ee par
$$1\rightarrow \X (M_{\alpha })\rightarrow\X(M)\rightarrow X(M_{\alpha }^1\cap M/M^1)\rightarrow 1,$$
la premi\`ere application \'etant la restriction \`a $M$ et la deuxi\`eme l'\'evaluation en $h_{\alpha
}$. Cette suite exacte courte est scind\'ee: il suffit de prolonger $h_{\alpha }$ en une base
$\underline {m}=(h_{\alpha },m_2,...,m_d)$ de $M$ modulo $M^1$. Une section $s_{\underline {m}}:
X(M_{\alpha }^1\cap M/M^1)\rightarrow\X(M)$ s'obtient alors, en associant \`a $\chi $ l'\'el\'ement de
$\X (M)$ qui vaut $\chi (h_{\alpha })$ en $h_{\alpha }$ et $1$ en tout $m_j$,
$j=2,\cdots, d$. Il n'est pas possible de d\'efinir un scindage "canonique".

\null {\bf 1.3} Fixons une repr\'esentation irr\'eductible cuspidale $(\sigma ,E)$ d'un sous-groupe de
Levi semi-standard $M$ de $G$. L'orbite inertielle de $\sigma $ est l'ensemble $\o=\o _{\sigma }$
form\'e des classes d'\'equivalence des repr\'esentations $\sigma\otimes\chi $ avec $\chi $ parcourant
$\X (M)$. On \'ecrira $\o _0=\o _{\sigma ,0}$ pour le sous-ensemble de $\o $ form\'e des classes
d'\'equivalence des repr\'esentations unitaires. On a sur $\o $ (resp. $\o _0$) une structure de
vari\'et\'e alg\'ebrique complexe (resp. analytiques r\'eelle) d\'eduite de sa structure d'espace
principal homog\`ene sous $\X (M)$ (resp. $\Xim (M)$). Pour tout $\sigma\in\o $, il existe un unique
caract\`ere non ramifi\'e $\chi _{\lambda }$ de $\X (M)$ avec $\lambda\in a_M^*$, tel que
$\sigma\otimes\chi _{\lambda }^{-1}$ soit une repr\'esentation unitaire. On \'ecrira $\Re(\sigma
):=\lambda $, $\Im (\sigma )=\sigma\otimes\chi_{-\lambda }$. On notera $\Re _{M'}(\sigma )$ la
projection de $\Re (\sigma )$ sur $a_{M'}^*$ et $\Im _{M'}(\sigma )=\sigma\otimes\chi _{-\Re_{M'}
(\sigma )}$.

Le symbole $E_{\so }$ d\'esignera parfois un espace vectoriel complexe dans lequel se r\'ealisent les
repr\'esentations dans $\o $. On notera $W(M,\o )$ le sous-groupe de $W(M)$ form\'e des \'el\'ements
$w$ qui v\'erifient $w\o =\o$. Si $M'$ est un sous-groupe de Levi qui contient $M$, $\Stab _{M'}(\sigma
)$ ou $\Stab _{M'}(\o )$ d\'esignera le sous-groupe de $\X (M')$ form\'e des caract\`eres non
ramifi\'es $\chi $ v\'erifiant $\sigma\otimes\chi\simeq\sigma $ (pour un \'el\'ement $\sigma $ de $\o
$). Si $M'=M$, on omettra parfois l'indice $M$.

Une fonction complexe $f$ sur $\o $ sera dite polynomiale (resp. rationnelle), si la fonction $\X
(M)\rightarrow\CC$, $\chi\mapsto f(\sigma\otimes\chi )$ est polynomiale (resp. rationnelle). De m\^eme,
on d\'efinira la notion de fonction lisse sur $\o _0$.

On notera $i_P^G$ le foncteur de l'induction parabolique normalis\'ee qui envoie les repr\'esen-tations
unitaires sur des repr\'esentations unitaires. D\'esignons par $E_{\chi }$ l'espace $E$ muni de la
repr\'esentation
$\sigma\otimes\chi $. Notons $i_{P\cap K}^KE$ l'espace des fonctions $f:K\rightarrow E$ invariantes \`a
droite par un sous-groupe ouvert de $K$ et v\'erifiant $f(muk)=\sigma (m)f(k)$ pour tout $m\in M\cap K$,
$u\in U\cap K$ et $k\in K$. La restriction \`a $K$ d\'efinit pour tout caract\`ere non ramifi\'e $\chi
$ de $M$ un isomorphisme entre
$i_P^GE_{\chi }$ et $i_{P\cap K}^KE$. Ainsi toutes les repr\'esentations $i_P^G(\sigma\otimes\chi )$ se
r\'ealisent dans le m\^eme espace $i_{P\cap K}^KE$. 

Ceci permet de d\'efinir la notion un peu subtile d'une application polynomiale ou rationnelle sur $\o
$ \`a valeurs dans $i_{P\cap K}^KE\otimes i_{P'\cap K}^KE^{\vee }$ ($E^{\vee }$ d\'esignant la
contragr\'ediente de $E$) etc.  (cf. [W] VI.1). (Ce sont des applications v\'erifiant une certaine
propri\'et\'e d'invariance par rapport \`a $\Stab (\o )$.) Si $V$ est un tel espace, on notera $\P
(\o,V)$ l'espace des applications polynomiales de $\o $ dans $V$ dont l'image est incluse dans un
sous-espace de dimension finie.

\null
\null {\bf 1.4} Lorsque $M$ est un sous-groupe de Levi semi-standard de $G$ et que $\sigma $ est une
repr\'esentation irr\'eductible lisse d'un sous-groupe de Levi semi-standard de $M$, on pose $\o
_{\sigma ,M}=\{\sigma\otimes\chi
\vert \chi \in\X (M)\}$. On appelle cet espace \it l'orbite inertielle de $\sigma $ sous $M$ ou
relative \`a $M$. \rm Une application $\psi :\o _{\sigma ,M}\rightarrow\CC $ sera dite polynomiale
(resp. rationnelle), si l'application $\X (M) \rightarrow \CC$, $\chi \mapsto \psi(\sigma\otimes\chi )$
est polynomiale (resp. rationnelle).

Fixons maintenant $M$ et $\sigma $. Un \it sous-espace affine (radiciel) \rm de $\o _{\sigma ,M}$ sera
pour nous un espace de la forme $\o _{\sigma ',M'}$ avec $\sigma '\in\o_{\sigma, M}$ et $M'$ un
sous-groupe de Levi  de $G$ qui contient $M$. Lorsque $M'=M_{\alpha }$ avec $\alpha\in\Sigma(T_M)$ et
que la r\'ef\'erence \`a $M$ para\^\i t \'evidente, on \'ecrira simplement $\o _{\sigma', \alpha }$.
(C'est par exemple le cas, si $\sigma '$ est une repr\'esentation de $M$.)  Un tel espace sera appel\'e
un \it hyperplan affine (radiciel) \rm de $\o _{\sigma ,M}$. On munira $\o _{\sigma ,\alpha }$
implicitement de l'orientation dans $\o $ d\'etermin\'ee par $\alpha $.

On pose $\Re (\o _{\sigma ',M'})=\Re (\sigma ')+a_{M'}^*$. (C'est l'ensemble des parties r\'eelles des
\'el\'ements dans $\o _{\sigma ,M'}$.) On dira que $M'$ est le \it sous-groupe de Levi associ\'e \`a
$A=\o _{\sigma ',M'}$ ou \`a $L=\Re (\o _{\sigma ',M'})$. \rm On le notera $M_A$ ou $M_L$.

Remarquons que, par d\'efinition de $\o $, les \'el\'ements de $\Re(\o )$ ont tous la m\^eme projection
sur $a_0^{M*}$. On appellera cet \'el\'ement l'origine de $\Re(\o )$ et on le notera $r(\o )$. (Si $\o
$ est form\'e de repr\'esentations de $M$, alors $r(\o )=0$.) C'est l'\'el\'ement de $\Re (\o )$ qui
est le plus proche de l'origine de $a_0^*$ pour la structure euclidienne d\'efinie sur cet espace.

Pour $\alpha\in\Sigma (T_M)$, $c\in\CC$, d\'esignons par $H_{\alpha ,c}$ l'hyperplan de $a_{M}^*$
donn\'e par $\{\lambda\in a_{M'}^*\vert\langle\lambda ,\alpha ^{\vee }\rangle =c\}$ et muni de
l'orientation donn\'ee par $\alpha $. Remarquons que l'on a, en posant $c=\langle\Re(\sigma '),\alpha
^{\vee } \rangle$, l'\'egalit\'e $\Re (\o _{\sigma ',\alpha })=H_{\alpha ,c}$.

Consid\'erons maintenant un ensemble fini $\S $ d'hyperplans affines radiciels de $\o _{\sigma ,M}$. On
va y associer plusieurs objets: on notera $\A (\S )$ l'ensemble form\'e des sous-espaces affines de $\o
$ qui sont les composantes connexes des intersections de sous-ensembles  finis ou vides de $\S $. On
consid\'erera \'egalement l'ensemble $\H =\H(\S )$ form\'e des parties r\'eelles des  hyperplans dans
$\S $. C'est  un ensemble fini d'hyperplans affines  de l'espace affine $r(\o )+a_M^*$. On y associe
l'ensemble $\L :=\L(\H):=\L (\S )$ dont les \'el\'ements sont les sous-espaces affines de $r(\o)+a_M^*$
obtenus en intersectant les \'el\'ements d'un sous-ensemble (fini ou vide) de $\H $.

Lorsque $A$ est un sous-espace affine radiciel de $\o _{\sigma ,M}$ (resp. $L$ un sous-espace affine
radiciel de $r(\o )+a_M^*$), on pose $\S _A=\{S\in\S\vert S\supseteq A\}$ (resp. $\H _L=\{L\in\L \vert
H\supseteq L\}$). On d\'esignera par $\S (A)$ (resp. $\H (L)$) l'ensemble fini d'hyperplans affines de
$A$ (resp. $L$) dont les \'el\'ements sont les composantes connexes des espaces non vides de la forme
$S\cap A$ avec $S\in \S-\S _A$ (resp. $H\cap L$ avec $H\in \H-\H _L$). On pose $A_{\sS,
\sreg}:=A-\bigcup _{S\in\sS (A)}S$ et $L_{\sS, \sreg}:= L_{\sH, \sreg} := L-\bigcup _{H\in\sH (L)}H$.
Quand cela ne pr\^ete pas \`a  confusion, on \'ecrira plus simplement $A_{\sreg }$ et $L_{\sreg }$.

L'\'el\'ement de $L$ qui est le plus proche de l'origine $r(\o )$ de $\Re (\o )$ (d\'efinie ci-dessus)
pour la structure euclidienne induite par celle sur $a_0^*$ sera not\'e $r(L)$. Si $L$ est la partie
r\'eelle de $A\in\A(\S )$, alors on pose $r(A):=r(L)$.

Si $C$ est une composante connexe de $L$, et si $H\in\H (L)$ est tel que la cl\^oture $\ol {C}$ de $C$
ait une intersection non vide avec $H$, on appelle l'int\'erieur de $H\cap \ol {C}$ une face de $C$.
Cette face est alors \'egale \`a $\ol {C}\cap H_{\sH (H),\sreg }$. Deux composantes connexes $C_1$ et
$C_2$ de $L_{\sreg }$ seront dites adjacentes, si elles ont une face en commun. Cette face est alors
unique. Notons $H$ l'unique \'el\'ement de $\H (L)$ qui contient cette face commune. De tout point de
$C_1$ \`a tout point de $C_2$ il existe un chemin qui ne coupe $\bigcup _{H'\in\sH (L)}H'$ qu'en un
point r\'egulier de $H$.

\null
\null {\bf 1.5} Soit $P=MU$ un sous-groupe parabolique semi-standard de $G$ et $(\pi ,V)$ une
repr\'esentation lisse de $M$. On notera $E_{P,\pi }^G$ l'application lin\'eaire $i_{P\cap
K}^KV\otimes i_{P\cap K}^KV^{\vee } \rightarrow C^{\infty }(G)$ qui associe \`a $v\otimes v^{\vee }$ le
coefficient matriciel $g\mapsto\langle (i_P^G\pi )(g)v,v^{\vee }\rangle$ de la repr\'esentation
$i_P^G\pi $. On omettra l'indice $\pi $ parfois  pour \'ecrire plus simplement $E_P^G$, si cela ne
pr\^ete pas \`a des confusions.

Si $\o $ est l'orbite inertielle d'une repr\'esentation irr\'eductible cuspidale de $M$, et si $\xi $
est une application rationnelle sur $\o $ \`a valeurs dans $i_{P\cap K}^KV\otimes i_{P\cap K}^KV^{\vee
}$, alors, pour tout $g\in G$, la fonction $\sigma\mapsto  E_{P,\sigma }^G(\zeta (\sigma ))(g)$ est
polynomiale sur $\o $, \it i.e. \rm dans $\P (\o,\CC )$. L'assertion analogue vaut pour $\xi $
rationnelle.

Nous consid\'ererons deux classes d'op\'erateurs d'entrelacement: 

Si $w$ est un \'el\'ement de $W$, nous d\'esignerons pour toute repr\'esentation $(\pi ,E_{\pi })$ de
$M$ par $\lambda (w)$ l'isomorphisme $i_P^GE_{\pi }\rightarrow i_{wP}^GwE_{\pi }$ entre les
repr\'esentations $i_P^G\pi$ et
$i_{wP}^Gw\pi$ donn\'e par $v\mapsto v(w^{-1}.)$.

Soit $P'=MU'$ un autre sous-groupe parabolique semi-standard de Levi semi-standard $M$. Soit $\o $
l'orbite inertielle d'une repr\'esentation irr\'eductible de $M$. Pour $\pi\in\o $, notons $v_{\pi }$
l'image d'un \'el\'ement de $i_{P\cap K}^KE_{\so }$ dans l'espace $i_{P}^GE_{\pi }$. On d\'efinit
l'op\'erateur d'entrelacement $J_{P\vert P'}$ dans un c\^one ouvert de $\o $ pour tout $v\in i_{P'\cap
K}^KE_{\so }$, $k\in K$ par l'int\'egrale convergente
$$(J_{P\vert P'}(\pi )v)(k)=\int _{U\cap U'\backslash U'} v_{\pi }(u'k)du'.$$ Par prolongement
analytique, on en d\'eduit un op\'erateur rationnel sur $\o $ \`a valeurs dans $\Hom _{\sCC}
(i_{P'\cap K}^KE_{\so }, i_{P\cap K}^K E_{\so }^{\vee })$ qui d\'efinit en tout point r\'egulier $\pi
$ de $\o $ un homomorphisme \'equivariant entre les repr\'esentations $i_P^G\pi $ et $i_{P'}^G\pi $
(cf. [W] chapitre IV.1). Lorsque $\o $ est l'orbite inertielle d'une repr\'esentation irr\'eductible
temp\'er\'ee, l'int\'egrale converge pour tout $\pi $ avec $\Re (\pi )>_{P'}0$ (cf. [W] chapitre IV.2).

En particulier, la fonction $\o\rightarrow\CC $, $\pi\mapsto E_{P,\pi }^G(J_{P\vert P'}(\pi )\xi (\pi
))(g)$ est rationnelle pour tout $\xi\in\P(\o, i_{P'\cap K}^K E_{\so }\otimes i_{P\cap K}^KE_{\so
}^{\vee })$ et tout $g\in G$.

Pour $\o $ l'orbite inertielle (resp. inertielle unitaire) d'une repr\'esentation irr\'eductible
cuspidale (resp. de carr\'e int\'egrable) d'un sous-groupe de Levi semi-standard de $G$ et $P$, $P'$
deux sous-groupes paraboliques semi-standard de facteur de Levi $M$, on notera $\mu _{P\vert P'}$ la
fonction rationnelle sur $\o $ telle que $\mu _{P\vert P'}(\pi )J_{P\vert P'}(\pi )J_{P'\vert P}(\pi )$
soit l'endomorphisme identit\'e pour tout point r\'egulier
$\pi $ de $\o $. Si $P=\ol {P}$, on \'ecrira $\mu =\mu_{P\vert\ol{P}}$. (Cette fonction ne d\'epend en
effet pas du choix de $P$ (cf. [W] chapitre IV).) La fonction $\mu $ est invariante pour l'action par
$W$ (cf. [W] V.2.1).

\null
\null {\bf 1.6} Pour tout $r\in a_{M'}^*$ et tout sous-groupe de Levi semi-standard $M'\supseteq M$,
posons $\o _{\sigma ,M', r}=\{\sigma '\in\o _{\sigma ,M'}\vert \Re _{M'}(\sigma ')=r \}$. Munissons
$\Xim (T_{M'})$ d'une mesure de Haar de masse totale $1$, $\Xim (M')$ d'une mesure pour laquelle la
restriction $\Xim (M')\rightarrow\Xim (T_{M'})$ conserve localement les mesures et $\o _{\sigma ,M',0}$
d'une mesure telle que l'applicaton $\Xim (M')$ $\rightarrow\o_{\sigma ,M',0}$,
$\chi\mapsto\sigma\otimes\chi $, conserve localement les mesures. Si $\psi $ est une fonction
rationnelle sur $\o _{\sigma ,M'}$ qui est r\'eguli\`ere sur $\o _{\sigma ,M', r}$, on notera $\int
_{\so _{\sigma ,M', r}}\psi (\sigma ') d\Im_{M'}(\sigma ')$ l'int\'egrale $\int _{\so _{\sigma ,M',0}}$
$\psi (\sigma '$ $\otimes \chi _r)d\sigma '$. On \'ecrira \'egalement $\int _{\Re (\sigma ')=r_0}\psi
(\sigma ') d_A\Im (\sigma ')$ avec $A=\o _{\sigma ,M'}$ et $r_0=r(A)+r$.

\null
\null
\null
\vfill
\eject
{\bf 2. Formule de Plancherel et th\'eor\`eme de Paley-Wiener}

\null
\null {\bf 2.1} On notera $\Theta $ l'ensemble des couples $(P,\o )$ form\'es d'un sous-groupe
parabolique semi-standard $P=MU$ et de l'orbite inertielle $\o $ d'une repr\'esentation irr\'eductible
cuspidale de $M$. Si $f$ est un \'el\'ement de l'espace $C_c^{\infty }(G)$ des fonctions lisses \`a
support compact sur $G$, sa transform\'ee de Fourier par rapport \`a $(P,\o )\in\Theta $, $\sigma\in\o
$, not\'ee $\widehat{f}(P,\sigma )$, est un \'el\'ement de $\End _{\sCC}(i_{P\cap K}^KE_{\so })$
d\'efini par
$$\widehat {f}(P,\sigma )v:=(i_P^G\sigma )(f)v=\int _Gf(g)(i_P^G\sigma (g)v)\ dg.$$ On obtient ainsi
une application polynomiale $\widehat {f}(P,\o ):\o\rightarrow i_{P\cap K}^KE_{\so }\otimes i_{P\cap
K}^KE_{\so }^{\vee }$, $\sigma\mapsto \widehat{f}(P,\sigma )$, o\`u on a identifi\'e $i_{P\cap
K}^KE_{\so }\otimes i_{P\cap K}^KE_{\so }^{\vee }$ \`a un sous-espace de $\End _{\sCC}(i_{P\cap
K}^KE_{\so })$.

\null {\bf 2.2} Fixons un \'el\'ement $(P,\o )$ de $\Theta $. On a prouv\'e dans [H1], proposition 2.1,
que, si $\xi $ est une application polynomiale sur $\o $ \`a valeurs dans un sous-espace de dimension
finie de $i_{\ol{P}\cap K}^KE_{\so }\otimes i_{P\cap K}^KE_{\so }^{\vee }$, alors la fonction $f_{\xi
}:G \rightarrow\CC $ d\'efinie avec $\deg (\sigma )$ \'egal au degr\'e formel de $\sigma $ par 
$$f_{\xi }(g):=\gamma (G/M)\int _{\Re (\sigma )=r\gg _P0} \deg(\sigma ) E_{P,\sigma
}^G((J_{\ol{P}\vert P}^{-1}\otimes 1)(\sigma )\xi (\sigma ))(g^{-1})d_{\so } \Im(\sigma )$$  est lisse
\`a support compact. (Ici la notation $r\gg _P0$ signifie \'evidemment que $r$ est tr\`es positif dans
la chambre de Weyl de $P$ dans $a_M^*$.)

Si $(P',\o')\in\Theta $ et si aucun \'el\'ement de $\o $ n'est conjugu\'e \`a un \'el\'ement de $\o '$,
alors $\widehat {f}_{\xi }(P',\o ')=0$  (cf. [H1] proposition 2.4.1). Pour $\sigma\in\o $, on a
l'identit\'e de fonctions rationnelles (cf. [H1] proposition 2.3)
$$\widehat {f}_{\xi }(P,\sigma )=\sum _{w\in W(M,\so )} (J_{P\vert\ol{wP}}(\sigma )\lambda(w)\otimes
J_{P\vert wP}(\sigma ^{\vee })\lambda (w))\xi (w^{-1}\sigma ).$$

\null  {\bf 2.3} L'espace de Paley-Wiener $PW(G)$ est par d\'efinition l'image de $C_c^{\infty }(G)$ 
par transformation de Fourier. Plus pr\'ecis\'ement, un \'el\'ement $\varphi =(\varphi _{P,\so
})_{(P,\so )}$ de
$\displaystyle{\bigoplus _{(P,\so )\in\Theta } \P(\o, }$
$\displaystyle{i_{P\cap K}^KE_{\so }\otimes i_{P\cap K}^KE_{\so }^{\vee })}$ est dans $PW(G)$, si et
seulement s'il existe $f\in C_c^{\infty }(G)$, tel que $\varphi _{P,\so }=\widehat {f}(P,\o )$ pour
tout $(P,\o )\in\Theta $.

Il est prouv\'e dans [H1] (cf. th\'eor\`eme 0.1) que $\varphi =(\varphi _{P,\so })_{(P,\so )}$ est dans
$PW(G)$, si et seulement si

(i) pour tout $(P,\o )\in\Theta $ et tout $w\in W^G$, on a $\lambda (w)\varphi _{P,\so }=\varphi
_{wP,w\so }\lambda (w)$; et si

(ii) pour tout $(P,\o )$, $(P',\o ')\in\Theta $, on a l'identit\'e de fonctions rationnelles
$$J_{P'\vert P}(\sigma )\varphi _{P,\so }(\sigma )=\varphi _{P',\so }(\sigma )J_{P'\vert P}(\sigma ).$$

La preuve de ce th\'eor\`eme utilise les r\'esultats rappel\'es au num\'ero pr\'ec\'edent ainsi que le
fait (cf. [H1] lemme 0.2) que, si la famille $\varphi =(\varphi _{P,\so })_{(P,\so )}$ v\'erifie les
propri\'et\'es (i) et (ii), il existe une application polynomiale $\xi _{\varphi _{P,\sso
}}:\o\rightarrow i_{\ol{P}\cap K}^KE_{\so }\otimes i_{P\cap K}^KE_{\so ^{\vee }}$ telle que
$$\varphi _{P,\so }(\sigma )=\sum _{w\in W(M,\so )} (J_{P\vert\ol{wP}}(\sigma )\lambda (w)\otimes
J_{P\vert wP}(\sigma )\lambda (w))
\xi _{\varphi _{P,\sso }}(w^{-1}\sigma ).$$ L'\'el\'ement de $C_c^{\infty }(G)$ qui correspond \`a
$\varphi $ est d\'etermin\'e de fa\c con unique (cf. [H1] corollaire 3.1.2). On le notera $f_{\varphi
}$.

\null {\bf 2.4} Fixons maintenant $\varphi =(\varphi _{P,\so })_{(P,\so )}\in PW(G)$ et $(P,\o
)\in\Theta $. Soit $\sigma\in\o $ et soit $P'=M'U'$ un sous-groupe parabolique semi-standard de $G$
contenant $P$. Soit $(\pi, V)$ une repr\'esentation de $M'$ qui est un sous-quotient de $i_{P\cap
M'}^{M'}\sigma $. Alors, comme $\varphi _{P,\so }(\sigma )$ laisse - comme transform\'ee de Fourier
d'un \'el\'ement de $C_c^{\infty }(G)$ - invariants les sous-espaces $G$-\'equivariants de $i_P^GE$,
cet endomorphisme d\'efinit par passage au quotient un endomorphisme $i_{P'}^GV\rightarrow i_{P'}^GV$.
Notant $\o _{\pi ,M'}$ l'orbite de $\pi $ sous $M'$, $\varphi _{P,\so }$ induit ainsi une application
polynomiale $\o _{\pi, M'}\rightarrow i_{P'\cap K}^KV \otimes i_{P'\cap K}^KV^{\vee }$ que l'on notera
$\varphi _{P',\so _{\pi, M'}}$. On d\'esignera sa valeur en $\pi'\in\o _{\pi,M'}$ parfois par $\varphi
(P',\pi ')$. Elle est \'egale \`a l'endomorphisme $(i_{P'}^G\pi ')(f_{\varphi })$. 

\null  {\bf 2.5} Notons $\Theta _2$ l'ensemble des couples $(P,\o _0)$ form\'e d'un sous-groupe
parabolique semi-standard $P=MU$ et de l'orbite inertielle unitaire $\o _0$ d'une repr\'esentation
irr\'educti-ble de carr\'e int\'egrable de $M$. Si $(P,\o )\in\Theta $, on note $\Theta _2(\o )$ le
sous-ensemble de $\Theta _2$ form\'e des couples $(P'=M'U',\o _0')\in\Theta _2$ tels que $M'\supseteq
M$ et que le support cuspidal de tout \'el\'ement de $\o _0'$ soit donn\'e par la classe de
$W^{M'}$-conjugaison d'un \'el\'ement de $\o $. 

Soit maintenant $\varphi $ un \'el\'ement de $PW(G)$ tel que $\varphi _{P',\so '}=0$ sauf si $\o '$ est
conjugu\'e \`a $\o $. Appelons deux couples $(P'=M'U',\o_0')$ et $(P''=M''U'',\o_0'')$ de $\Theta _2(\o
)$ \it associ\'es, \rm s'il existe $w\in W$ tel que $wM'=M''$ et $w\o'_0=\o''_0$. Rappelons que
$W(M,\o)$ d\'esigne l'ensemble des $w\in W(M)$ tels que $w\o=\o$. Avec $\gamma (G/M)$ une constante
(d\'efinie dans [W] I.1), la formule de Plancherel de Harish-Chandra donne alors (cf. [W] th\'eor\`eme
VIII.1.1)
$$\displaylinesno{
\qquad f_{\varphi }(g)=\sum _{(P'=M'U',\so _0')\in\Theta _2(\so )/ass}\gamma (G/M')\vert
W(M',\o_0')\vert^{-1}\hfill&(2.5.1)\cr
\hfill\int _{\so _0'}\deg(\pi ')E_{P',\pi '}^G(\varphi (P',\pi '))(g^{-1}) \mu(\pi ') d\pi' .\qquad\cr
}$$   La somme des termes index\'e par des couples $(P',\o _0')$ avec $P'\ne G$ sera appel\'ee la \it
partie continue \rm de la formule de Plancherel.  La somme des autres termes sera la \it partie
discr\`ete.\rm

\null
\null
\null
\vfill
\eject
{\bf 3. Th\'eor\`eme des r\'esidus}

\null
\null {\bf 3.1} Fixons un sous-groupe de Levi semi-standard $M$ de $G$. Soit $\o $ l'orbite inertielle
sous $M$ d'une repr\'esentation irr\'eductible cuspidale d'un sous-groupe de Levi semi-standard de
$M$. Fixons par ailleurs un ensemble fini $\S $ d'hyperplans affines radiciels de $\o $. On
consid\'erera dans la suite de ce paragraphe l'espace $\o $ muni de cet ensemble d'hyperplans.

On note $\R (\o , \S )$ l'ensemble des fonctions rationnelles $\psi $ sur $\o $ qui sont r\'eguli\`eres
sur $\o _{\sreg }:=\o -\bigcup _{S\in\sS } S$. Posons $\A =\A(\S )$, $\H=\H(\S )$ et $\L=\L(\S)$ (cf.
{\bf 1.4} pour la d\'efinition de ceci).

\null {\bf 3.2} Pour donner des \'equations pour les hyperplans affines radiciels, explicitons des
\'el\'ements particuliers de $M$.

\null {\bf Lemme:} \it Soit $\alpha\in\Sigma _{red}(P)$. Il existe une puissance $\ti {h}_{\alpha }$ de
$h_{\alpha }$ par un entier $>0$ tel que $\chi\in\X (M)$ v\'erifie $\chi (\ti {h}_{\alpha })=1$ si et
seulement si $\chi\in\Stab_M(\o )\X (M_{\alpha })$.

\null Preuve: \rm Soit $m>0$ minimal tel que $\chi (h_{\alpha }^m)=1$ pour tout $\chi\in\Stab _M
(\o )$. Montrons que $\ti {h}_{\alpha }:=h_{\alpha }^m$ v\'erifie les propri\'et\'es de
l'\'enonc\'e. Comme $h_{\alpha }\in M_{\alpha }^1$, on a \'evidemment $\chi '(h_{\alpha })=1$ pour
tout $\chi '\in\X (M_{\alpha })$, d'o\`u l'une des deux implications. Soit $\chi\in\X (M)$ avec $\chi
(h_{\alpha }^m)=1$. Ceci implique que
$\chi (h_{\alpha })$ est une racine $m^{\hbox{\srm \`eme }}$ de l'unit\'e. Par minimalit\'e de $m$, il
existe $\chi _1$ dans $\Stab_M (\o )$ avec $\chi _1(h_{\alpha })=\chi (h_{\alpha })$. Il en suit que
$(\chi\chi _1^{-1})(h_{\alpha })=1$. Il en r\'esulte alors par choix de $h_{\alpha }$ et {\bf 1.2} que
$\chi\chi _1^{-1}\in\X (M_{\alpha })$. \hfill {\fin 2}

\null On \'ecrira $\ti{\alpha }$ pour le multiple de $\alpha $ par un r\'eel $>0$ tel que $\langle
H_M(\ti {h}_{\alpha }),\ti {\alpha }\rangle =1$.

\null
\null {\bf Corollaire:} \it Soit $\sigma\in\o $. Alors, pour que $\chi\in\X (M)$ v\'erifie
$\sigma\otimes\chi \in\o_{\sigma, \alpha }$, il faut et il suffit que $\chi(\ti{h}_{\alpha })=1$. 

\null Preuve: \rm On a $\sigma\otimes\chi\in\o _{\sigma, \alpha }$, si et seulement s'il existe
$\chi'\in\X (M_{\alpha })$ tel que
$\sigma\otimes\chi\simeq\sigma\otimes\chi '$. Ceci \'equivaut \`a $\chi\chi'^{-1}\in\Stab(\sigma )$,
d'o\`u
$\chi\in\Stab (\sigma )\X (M_{\alpha })$, ce qui \'equivaut \`a $\chi (\ti{h}_{\alpha })=1$ par le
lemme. \hfill {\fin 2} 

\null
\null {\bf 3.3 Lemme:} \it Soient $A'$ et $A''$ deux sous-espaces affines de $\o $, $A''\not\subseteq
A'$. Alors il existe une fonction polynomiale sur $\o
$ qui est identiquement nulle sur $A'$ et qui n'est pas identiquement nulle sur $A''$.

\null Preuve: \rm Notons $M'$ et $M''$ les sous-groupes de Levi semi-standard de $G$ associ\'es \`a
$A'$ et $A''$ respectivement (cf. {\bf 1.4}). Si $M'\not\subseteq M''$, il existe $\alpha\in\Sigma
_{red}(T_M)$ avec $M_{\alpha }\subseteq M'$ et $M_{\alpha }\not\subseteq M''$. On en d\'eduit
l'existence d'un hyperplan $\o _{\sigma ,\alpha }$ avec $A'\subseteq\o_{\sigma, \alpha }$ et
$A''\not\subseteq \o_{\sigma ,\alpha }$. On peut alors prendre le polyn\^ome d\'ecrivant
$\o_{\sigma ,\alpha }$ explicit\'e dans {\bf 3.2}.

Si $M'\subseteq M''$, on se ram\`ene \`a $M'=M''$, en remplacant $A''$ par un sous-espace
affine de $\o $ plus grand qui n'est pas contenu dans $A'$. Les arguments dans
l'exemple 2, V.1.2 dans [MW1] se g\'en\'eralisent alors et donnent un polyn\^ome ad\'equat sur
$\o $.
\hfill {\fin 2}

\null
\null {\bf 3.4} Pour $\sigma \in \o$, $m\in M$, posons $p_{\sigma ,m}: \o\rightarrow\CC $, $\sigma
\otimes\chi \mapsto (\chi (m)-1)$, lorsque cette fonction est bien d\'efinie. Si $m=\ti{h}_{\alpha }$
avec $\alpha\in\Sigma (T_M)$, on \'ecrira $p_{\sigma ,\alpha }$ ou encore $p_{\so _{\sigma, \alpha }}$.
Il r\'esulte de {\bf 3.2} que $\o _{\sigma ,\alpha }$ est l'ensemble des \'el\'ements de $\o $ en
lesquels $p_{\sigma, \alpha }$ s'annule.

Pour $H=H_{\alpha ,c}$ un hyperplan affine radiciel de $a_M^*$ et $\lambda \in a_M^*$, on pose
$p_H(\lambda )=p_{\alpha ,c}(\lambda )=\langle
\alpha ^{\vee },\lambda \rangle $.

Si $\psi\in\R (\o ,\S )$, alors il existe pour tout $S=\o _{\sigma ,\alpha }\in\S $ un unique entier
$d=d_{\sigma ,\alpha }(\psi )=d_S(\psi )$, tel que $\sigma '\mapsto p_{\sigma ,\alpha }(\sigma ')^d\psi
(\sigma ')$ soit r\'eguli\`ere et non identiquement nulle sur $(\o _{\sigma ,\alpha })_{\sreg }$. Si
$d_{\sigma ,\alpha }(\psi )\geq 1$, on parlera \it d'un p\^ole d'ordre $d_{\sigma ,\alpha }(\psi )$.
\rm Si $d_{\sigma ,\alpha }(\psi )\leq -1$, on parlera d'un \it z\'ero  d'ordre $-d_{\sigma ,\alpha
}(\psi )$. \rm

Si $\S $ contient tous les hyperplans singuliers de $\psi $ et de $1/\psi $, on appelle, pour
$\sigma\in\o $,
$d_{\sigma }(\psi ):=\sum _{S\in\sS,
\sigma\in S}d_S(\psi )$ \it l'ordre de $\psi $ en $\sigma $. \rm Si $d=d_{\sigma }(\psi )>0$, on parle
d'un p\^ole d'ordre $d$, et, si $d<0$, d'un z\'ero d'ordre $-d$.

\null
\it D\'efinition: \rm Pour $\sigma'\in(\o _{\sigma ,\alpha })_{\sreg }$, on pose $(\Res _{\so _{\sigma
,\alpha }}\psi )(\sigma ')=\Res _{z=0}\psi (\sigma '\otimes \chi _{z\ti {\alpha }})$ (en se rappelant
que $\o _{\sigma ,\alpha }$ est muni d'une orientation donn\'ee par $\alpha $). La fonction $\Res _{\so
_{\sigma ,\alpha }}\psi $ est donc d\'efinie sur un ouvert Zariski dense de $\o _{\sigma ,\alpha }$.

\null {\bf 3.5} Pour $\sigma'\in\o $ et $\psi\in\R (\o,\S )$, notons $\psi _{\sigma }$ la fonction
m\'eromorphe sur
$a_{M,\sCC }^*$, donn\'ee par
$\lambda\mapsto \psi _{\sigma '}(\lambda ):=\psi (\sigma'\otimes\chi _{\lambda })$. Notons $D_{\ti
{\alpha }}$ la d\'eriv\'ee en direction de
$\ti{\alpha }$.

\null {\bf Lemme:} \it Supposons que $\psi $ ait un p\^ole d'ordre $d\geq 1$ sur $\o _{\sigma ,\alpha
}$. Alors, pour tout $\sigma '\in (\o _{\sigma ,\alpha })_{\sreg }$, on a
$$(\Res _{\so _{\sigma ,\alpha }}\psi )(\sigma ')={1\over (d-1)!} (D_{\ti{\alpha }}^{d-1}
(\langle\lambda ,H_M(\ti {h}_{\alpha })\rangle ^d\psi _{\sigma '}(\lambda )))_{\vert \lambda =0}.$$ En
particulier, la fonction $\Res _{\so _{\sigma ,\alpha }}\psi $ est rationnelle.

\null Preuve: \rm Par d\'efinition, $(\Res _{\so _{\sigma ,\alpha }}\psi )(\sigma ')=(\Res _{z=0}\psi
)(\sigma '\otimes\chi _{z\ti{\alpha }})$. Comme $\chi _{z\ti {\alpha }}(\ti {h}_{\alpha })-1=q^{-z}-1$
admet un z\'ero d'ordre
$1$ en $z=0$, on trouve $(\Res _{z=0}\psi )(\sigma '\otimes\chi_{z\ti {\alpha }})={1\over (d-1)!}
({d\over dz})^{d-1} (\langle z\ti{\alpha }, H_M(\ti {h}_{\alpha })\rangle ^d\psi _{\sigma
'}(z\ti{\alpha }))_{\vert z=0}$. Ceci est bien
\'egal \`a l'expression dans l'\'enonc\'e. 

\hfill{\fin 2}

\null
\null {\bf 3.6} Soient $C_1$ et $C_2$ deux composantes connexes adjacentes de $\Re (\o )_{\sreg }$.
Notons $H_{\alpha ,c}$ l'hyperplan affine qui les s\'epare (cf. {\bf 1.4}) muni de l'orientation
donn\'ee par $\alpha $. On dira que $C_i$ est \it n\'egatif \rm (resp. \it positif) pour (l'orientation
de) $H_{\alpha ,c}$, \rm si les \'el\'ements de $C_i$ v\'erifient $\langle\alpha ^{\vee
},\lambda\rangle<c$ (resp. $\langle\alpha ^{\vee },\lambda\rangle>c$).

\null  {\bf Proposition:} \it Soit $\psi $ dans $\R (\o ,\S )$. Soient $C_1$ et $C_2$ deux composantes
connexes adjacentes dans $\Re (\o )_{\sreg }$ s\'epar\'ees par $H_{\alpha ,c}$ et supposons $C_1$
n\'egatif et $C_2$ positif pour $H_{\alpha ,c}$. Pour $r_1\in C_1$, $r_2\in C_2$ et $r_{12}$ dans la
face commune  de $C_1$ et $C_2$, on a alors, modulo une constante $\kappa >0$ qui ne d\'epend pas de
$\psi $,
$$\int _{\Re(\sigma )=r_2}\psi (\sigma )d_{\so }\Im (\sigma )-\int _{\Re (\sigma )=r_1} \psi (\sigma
)d_{\so }\Im(\sigma )=\kappa \sum _S \int _{\Re (\sigma )= r_{12}}(\Res _{S}\psi )(\sigma )d_S\Im
(\sigma ),$$  la somme portant sur les hyperplans $S\in \S$ de partie r\'eelle $H _{\alpha ,c}$. 

\null Preuve: \rm Posons $M'=M_{\alpha }$. Comme les int\'egrales ne d\'ependent pas du choix de
$r_j\in C_j$, on peut supposer que $r_2-r_1$ soit un multiple de $\ti {\alpha }$ par un r\'eel $>0$ et
que $r_{12}$ soit sur la droite reliant $r_1$ et $r_2$. Fixons $\sigma_0\in\o $ tel que la projection
de $\Re (\sigma _0)$ sur $a_M^*$ soit $0$. Posons $r_{j,\alpha }=\langle H_M(\ti {h} _{\alpha }),
r_j\rangle $. Par choix de $r_1$, $r_2$ et $r_{12}$, ils ont la m\^eme  projection sur $a_{M'}^*$ que
l'on note $r'$. On obtient alors, avec $\kappa $ une constante positive qui d\'epend uniquement du
choix des mesures que 
$$\eqalign { &\int _{\Re(\sigma )=r_2}\psi (\sigma )d_{\so }\Im (\sigma )-\int _{\Re (\sigma )=r_1}
\psi (\sigma )d_{\so }\Im (\sigma )\cr  =&\kappa ({\log q\over 2\pi }\int _0^{2\pi\over\log q}\int
_{\so _{\sigma _0,M',r'}}\psi (\sigma \otimes\chi _{(r_{2,\alpha }+it)\ti{\alpha }}) d\Im _{M'}(\sigma
)dt\cr &-{\log q\over 2\pi }\int _0^{2\pi\over \log q} 
\int _{\so _{\sigma _0,M',r'}}\psi (\sigma \otimes\chi _{(r_{1,\alpha }+it)\ti {\alpha }}) d\Im
_{M'}(\sigma ) dt)\cr }$$   (cf. {\bf 1.6} pour la d\'efinition de $\o _{\sigma _0,M',r'}$ et celle de
l'int\'egrale int\'erieure).

Choisissons $\theta $ tel que, pour tout $\sigma \in \o _{\sigma _0,M',r'}$, $\sigma\otimes\chi
_{(i\theta+r_{12})\ti{\alpha }}$ ne soit sur aucun hyperplan singulier de partie r\'eelle $H_{\alpha
,c}$. La diff\'erence des int\'egrales est alors \'egale \`a
$$\eqalign {  &{\log q\over 2\pi }\int _{\theta +0}^{\theta +{2\pi\over\log q}}\int _{\so _{\sigma
_0,M',r'}}\psi (\sigma \otimes\chi _{(r_{2,\alpha }+it)\ti {\alpha }})d\Im _{M'}(\sigma )dt\cr  
-&{\log q\over 2\pi}\int _{\theta +0}^{\theta +{2\pi\over\log q}}\int _{\so _{\sigma
_0,M',r'}}\psi(\sigma\otimes\chi _{(r_{1,\alpha }+it)\ti {\alpha }}) d\Im _{M'}(\sigma )dt\cr  }$$  
Posons $f(z)=\int _{\so _{\sigma _0,M',r'}}\psi (\sigma\otimes\chi _{z\ti {\alpha }})d\Im_{M'}(\sigma
)$ pour tout 
$z\in\CC $ pour lequel cette expression fait un sense. Ceci d\'efinit une fonction m\'eromorphe sur
$\CC $.  Comme $f(r+i\theta )=f(r+i(\theta +2\pi/\log q))$ pour 
$r_{1,\alpha }\leq r\leq r_{2,\alpha }$, l'expression ci-dessus est \'egale \`a l'int\'egrale de
${1\over i}f(z)$ le long du contour $\Gamma $  ci-dessous, le facteur ${1\over i}$ venant du fait que
l'on int\`egre ci-dessus le long du chemin $t\mapsto it$.

\vskip 5cm

Pour tout $S\in\S $ de partie r\'eelle $H_{\alpha ,c}$, il existe un unique nombre complexe $c_S$ de
partie imaginaire dans $]\theta, \theta+2\pi [$, tel que $S$ soit l'orbite de $\sigma_0\otimes\chi
_{c_S\ti{\alpha }}$ sous
$M'$. Suite \`a nos hypoth\`eses, $\Re (c_S)=r_{12,\alpha }$.

L'expression ci-dessus devient alors
$$\eqalign {  &{1\over 2\pi i}\int _{\Gamma }(\int _{\so _{\sigma _0,M',r'}}\psi (\sigma\otimes\chi
_{z\ti {\alpha }})\ d\Im _{M'}(\sigma ))\ dz\cr  =&{1\over 2\pi i}\int _{\so _{\sigma _0,M',r'}}(\int
_{\Gamma }\psi (\sigma\otimes\chi _{z\ti {\alpha }})\ dz)\ d\Im _{M'}(\sigma )\cr  =&{1\over 2\pi
i}\int _{\so _{\sigma _0,M',r'}}  2\pi i \sum _S(\Res _{z=c_S}\psi )(\sigma\otimes\chi _{z\ti{\alpha
}})\ d\Im _{M'}(\sigma )\cr  }$$  Pour trouver l'expression de l'\'enonc\'e, il reste \`a remarquer
que  
$$(\Res _{z=c_S}\psi )(\sigma '\otimes\chi _{z\ti{\alpha }})= (\Res _{z=0}\psi) (\sigma '\otimes\chi
_{(z+c_S)\ti{\alpha }})=(\Res _{S}\psi )(\sigma ').$$
\hfill{\fin 2}

\null {\bf 3.7} Soit $A\in \A(\S)$ et soit $\D =(S_0,...,S_h)$ un ensemble ordonn\'e d'\'el\'ements de
$\A (\S )$ avec $S_0=\o $, $S_i$  hyperplan affine radiciel de $S_{i-1}$ pour $i=1,\cdots,h$ et
$S_h=A$, o\`u on suppose par ailleurs implicitement que $S_i$ soit muni d'une orientation en tant
qu'hyperplan affine de $S_{i-1}$. On pose $Res _{\sD }(\psi )=\Res _{S_h}(\Res _{S_{h-1}}(\cdots \Res
_{S_1}\psi ))$. L'op\'erateur $Res _{\sD }$ sera appel\'e un \it op\'erateur r\'esidu en $A$ sur $\o $
relatif \`a $\S $. \rm Une \it donn\'ee de r\'esidu en
$A$ sur $\o $ relative \`a $\S $, \rm not\'ee $Res _A$, sera alors une combinaison lin\'eaire \`a
coefficients r\'eels d'op\'erateurs r\'esidus en
$A$.

\it Remarque: \rm  A l'aide d'une r\'ecurrence, on d\'eduit facilement de la d\'efinition d'un
op\'erateur r\'esidu et du lemme {\bf 3.5} que tout op\'erateur r\'esidu en $A$ d\'efinit une
application lin\'eaire $\R (\o ,\S )$ dans $\R (A, $ $\S (A))$. Il en est alors de m\^eme d'une
donn\'ee de r\'esidu en $A$.

\null
\null Pour $\underline {d}\in\hbox{\ens N}^{\sS}$, notons $\R(\o,\S,\underline {d})$ le sous-espace de
$\R(\o, \S)$ form\'e des applications rationnelles $\psi $ telles que $d_S(\psi )\leq d_S$ pour tout
$S\in\S $ (cf. {\bf 3.4} pour la d\'efinition de $d_S(\psi )$).

\null {\bf Proposition:} \it Soit $A\in\A (\S)$ et soit $Res _A$ une donn\'ee de r\'esidu en $A$ sur
$\o $ relative \`a $\S $. Soit $\underline {d}\in\hbox{\ens N}^{\sS}$ et soit $M'=M_A$ le sous-groupe
de Levi associ\'e \`a $A$. Alors il existe un op\'erateur diff\'erentiel holomorphe $D=D_{\lambda }$
sur $a_{M,\sCC}^{M'}$ \`a coefficients r\'eels constants, tel que pour tout
$\psi\in\R(\o,\S,\underline{d})$ et tout $\sigma\in A_{\sreg }$, on ait
$$(Res_A\psi )(\sigma )=D_{\lambda }((\prod _{S\in \sS_A}p_{\Re (S)}^{d_S} (\lambda ))\psi _{\sigma
}(\lambda ))_{\vert \lambda =0}.$$  (avec $\psi _{\sigma }$ d\'efinie en {\bf 3.5}).

\null  Preuve: \rm  Il suffit de consid\'erer le cas o\`u $Res _A$ est un op\'erateur r\'esidu en $A$.
Il existe donc un ensemble ordonn\'e $\D =(S_0,\cdots, S_h)$ de sous-espaces affines dans $\o $ avec
$S_0=\o $, $S_i$ hyperplan de $S_{i-1}$ pour $i=1,\cdots, h$ et $S_h=A$, tel que $Res _A=Res _{\sD }$.
On va effectuer une r\'ecurrence sur $h$. Si $h=0$, il n'y a rien \`a montrer. Sinon, posons
$A'=S_{h-1}$ et $\D '=(S_0,\cdots ,S_{h-1})$. On a donc $Res _A=Res _{S_h}\circ Res_{\sD '}$. Posons
$M'=M_{A'}$. Soit $\alpha $ la racine positive dans $\Sigma (T_{M'})$ telle que $M'_{\alpha }$ soit le
sous-groupe de Levi associ\'e \`a $A$ et qui refl\`ete l'orientation de $A$ en tant qu'hyperplan de
$A'$. Par hypoth\`ese de r\'ecurrence, il existe un op\'erateur diff\'erentiel holomorphe \`a
coefficients r\'eels constants $D'$ sur $a_{M,\sCC}^{M' *}$, tel que, pour tout $\psi
\in\R(\o,\S,\underline {d})$, on ait
$$(Res_{A'}\psi )(\sigma\otimes\chi_{z\ti{\alpha }})=D'((\prod _{S\in\sS _{A'}}p_{\Re
(S)}^{d_S}(\lambda ))\psi _{\sigma }(\lambda +z\ti{\alpha }))_{\vert \lambda =0}$$  pour $\sigma\in
A_{\sreg }$ et $z\in\CC $ tel que $\sigma\otimes\chi _{z\ti{\alpha }}\in A_{\sreg }'$. Posons
$p(\lambda )=\prod _{S\in \sS_A -\sS_{A'}}p_{\Re(S)}^{d_S}(\lambda )$ pour $\lambda\in a_{M,\sCC }^*$. 

Par la r\`egle de Leibniz, il existe des polyn\^omes $Q_j$ \`a coefficients r\'eels sur $a_M^*$ et des
op\'erateurs diff\'erentiels $D_j$ sur $a_M^{M'*}$ holomorphes \`a coefficients r\'eels constants, tels
que l'on ait pour $z\ne 0$, mais pr\`es de $0$, et toute fonction $f$ m\'eromorphe sur $a_M^{M'*}$ et
holomorphe en $0$, avec $m$ un certain entier $\geq 0$,
$$D'(p(\lambda +z\ti{\alpha })^{-1}f(\lambda ))_{\vert\lambda =0}=p(z\ti {\alpha })^{-m}\sum _j
Q_j(z\ti{\alpha })(D_jf(\lambda ))_{\vert\lambda =0}.$$

On en d\'eduit que
$$(Res _{A'}\psi )(\sigma\otimes\chi _{z\ti{\alpha }})=p(z\ti {\alpha })^{-m} \sum _j Q_j(z\ti{\alpha
})D_j((\prod _{S\in\sS _A} p_{\Re (S)}^{d_S}(\lambda +z\ti{\alpha }))\psi _{\sigma }(\lambda
+z\ti{\alpha }))_{\vert\lambda =0}.$$

On observe alors qu'il existe un entier $m'$ et une constante $c\in \RR$ tels que $p(z\ti{\alpha
})=cz^{m'}$, et on trouve, puisque les fonctions de la forme $D_j(.) _{\vert\lambda =0}$ dans
l'expression ci-dessus sont r\'eguli\`eres en $z=0$,
$$\displaylines {\qquad\Res _{S_h}(Res_{A'}\psi )(\sigma)=\Res _{z=0} (Res_{A'}\psi )
(\sigma\otimes\chi _{z\ti{\alpha }})={1\over (mm'-1)!}\hfill\cr
\hfill [({d\over dz})^{mm'-1}c^{-m}\sum _j Q_j(z\ti{\alpha })D_j((\prod _{S\in\sS _A} p_{\Re
(S)}^{d_S}(\lambda +z\ti{\alpha }))\psi _{\sigma }(\lambda +z\ti{\alpha })) _{\vert\lambda =0}]_{\vert
z=0}.\cr }$$

En appliquant de nouveau la r\`egle de Leibniz, on obtient une expression du type voulu, en se
rappelant que les polyn\^omes $Q_j$ sont \`a coefficients r\'eels. \hfill {\fin 2}

\null
\null {\bf 3.8} Les r\'esultats techniques de ce num\'ero seront utiles avec ceux de {\bf 3.7} pour
\'etablir des propri\'et\'es d'unicit\'e des donn\'ees de r\'esidu. 

\null Par ce qui a \'et\'e expliqu\'e en {\bf 1.2}, on peut introduire des variables $z_1,\cdots ,z_n$
sur $\X (M)$ d\'eployant la suite exacte  $1\rightarrow\X (M')\rightarrow\X(M)$ $\rightarrow X(M\cap
M'^1/M^1)\rightarrow 1$, ainsi qu'une base de $a_{M,\sCC}^{M'*}$ form\'e d'\'el\'ements de $a_M^{M'*}$,
telles que

(i) les fonctions polynomiales sur $\X (M)$ sont les polyn\^omes de Laurent en $z_1,\cdots ,z_n$, les
fonctions polynomiales sur $X(M\cap M'^1/M^1)$ correspondant aux variables $z_1,\cdots
,z_m$ pour un certain entier $m$, $0\leq m\leq n$;

(ii) pour tout vecteur $\ul {k}=(k_1,\cdots ,k_m)$ form\'e d'entiers rationnels, les fonctions sur
$a_{M,\sCC}^{M'*}$ qui sont les compos\'ees de mon\^omes de Laurent sur $X(M'^1\cap M/M^1)$ de degr\'e
partiel donn\'e par le vecteur $\ul {k}$ avec l'application $\lambda\mapsto\chi _{\lambda }$ sont des
produits d'expressions de la forme $q^{k_j\lambda _j}$ par une constante, les $\lambda _j$ d\'esignant
les coordonn\'ees de $\lambda $.

\null La notion de degr\'e et de mon\^ome utilis\'ee dans le reste de cette section
sera celle d\'etermin\'ee par le choix de base et de variables ci-dessus.

\null  {\bf Lemme:} \it Soit $D$ un op\'erateur diff\'erentiel holomorphe sur $a_{M,\sCC}^{M'*}$ \`a
coefficients constants donn\'e par un mon\^ome unitaire non nul. Soit $\sigma\in\o $. Il existe une
fonction polynomiale $\psi $ sur $\o $ telle que $(D_{\lambda }\psi _{\sigma }(\lambda ))_{\vert\lambda
=0}\ne 0$, mais $(D'_{\lambda }\psi _{\sigma }(\lambda ))_{\vert\lambda =0}=0$ pour tout op\'erateur
diff\'erentiel holomorphe \`a coefficients constants $\ne D$ venant d'un mon\^ome unitaire de degr\'e
total inf\'erieur ou \'egal \`a celui de $D$.

\null Preuve: \rm Soit $\psi $ une fonction polynomiale sur $\o $, invariante pour l'action par $\X
(M')$. Il existe des fonctions $p_{\ul{k}}$ sur $a_{M,\sCC}^{M'*}$ qui sont les compos\'ees de
mon\^omes de Laurent sur $X(M'^1\cap M/M^1)$ de degr\'es partiels donn\'es par le vecteur $\ul{k}$ par
l'application $\lambda\mapsto\chi _{\lambda }$, tels que, pour tout $\lambda\in a_M^{M'*}$, on ait
$\psi (\sigma\otimes \chi_{\lambda })=\sum _{\ul{k}}p_{\ul{k}} (\lambda )$. Par nos choix, si $\lambda
_1,\cdots, \lambda _m$ et $k_1,\cdots ,k_m$ sont les coordonn\'ees de $\lambda $ et $\ul{k}$
respectivement, $p_{\ul{k}}$ est un produit d'expressions de la forme $q^{ k_j\lambda _j}$ par une
constante.

A tout op\'erateur diff\'erentiel holomorphe \`a coefficients constants $D$ correspond donc un
polyn\^ome $Q_D$ sur $\hbox{\ens R}^{\dim (a_M^{M'*})}$, tel que $(D_{\lambda }\psi(\sigma\otimes\chi
_{\lambda }))_{\vert\lambda =0}=\sum _{\ul{k}}Q_D(\log q \ul{k})p_{\ul{k}}(0)$. L'application
$D\rightarrow Q_D$ est lin\'eaire, et, si $D$ vient d'un mon\^ome, $Q_D$ est un mon\^ome du m\^eme
degr\'e.

Tout revient donc \`a choisir une famille $(a_{\ul{k}})$ de nombres complexes dont presque tous - sauf
un nombre fini - sont nuls, telle que $(a_{\ul{k}})$ soit orthogonale \`a $(Q_{D'}(\log q\ul{k}))$ pour
tout $D'\ne D$ de degr\'e total inf\'erieur ou \'egal \`a celui de $D$, et non orthogonal \`a
$(Q_D(\log q\ul{k}))$. Or, ceci est toujours possible, puisque des mon\^omes de degr\'es partiels
distincts sont lin\'eairement ind\'ependants. 

\hfill{\fin 2}

\null  {\bf Proposition:} \it Soit $D$ un op\'erateur diff\'erentiel holomorphe sur $a_{M,\sCC}^{M'*}$
\`a coefficients constants. Soit $\sigma\in\o $, soit $p$ un polyn\^ome sur $\o $ qui n'est pas
identiquement nul sur $A:=\o _{\sigma ,M'}$, et soit $\ti{p}$ une fonction m\'eromorphe sur
$a_{M,\sCC}^{M'*}$ qui est r\'eguli\`ere et non nulle dans l'origine de $a_{M,\sCC}^{M'*}$. Supposons
$$\int _{\Re(\sigma ')=\Re(\sigma )}D_{\lambda }(\ti{p}(\lambda )p_{\sigma ' } (\lambda)\psi _{\sigma
'}(\lambda ))_{\vert \lambda =0}d_A\Im(\sigma ')=0$$  pour toute fonction polyn\^omiale $\psi $ sur $\o
$.

Alors $D=0$. \rm

\null Montrons d'abord le lemme suivant:

\null {\bf Lemme:} \it Sous les hypoth\`eses de la proposition, on a $D_{\lambda }(\ti{p}(\lambda
)p_{\sigma '}(\lambda )\psi _{\sigma '}(\lambda )) _{\vert\lambda =0}=0$ pour tout $\sigma'\in\o $
avec $\Re (\sigma ')=\Re (\sigma )$. 

\null Preuve: \rm On peut \'ecrire
$$p_{\sigma '}(\lambda )\psi_{\sigma '}(\lambda )=\sum _{\ul{k}} m_{\ul{k}}(\sigma
')\psi_{\ul{k}}(\lambda ),$$  o\`u $m_{\ul{k}}(\sigma ')$ est un polyn\^ome sur $\o_{\sigma ,M'}$ dont
les degr\'es partiels par rapport aux variables de $\X (M')$ sont donn\'es par le vecteur $\ul{k}$, et
o\`u $\psi_{\ul{k}}$ est une fonction sur $a_{M,\sCC}^{M'*}$ qui est le compos\'e d'un polyn\^ome sur
$X(M\cap M'^1/M^1)$ avec l'application $\lambda\mapsto\chi _{\lambda }$. Il suffit de prouver que
$D_{\lambda }(\ti{p}(\lambda )\psi_{\ul{k}}(\lambda ))_{\vert\lambda =0}=0$ pour tout $\ul{k}$ avec
$m_{\ul{k}}\not\equiv 0$.

Or, par lin\'earit\'e de $D_{\lambda }$, on trouve sous les hypoth\`eses de la proposition par le
th\'eor\`eme de Cauchy avec $\kappa $ une constante non nulle et $\ul{0}$ le vecteur dont toutes les
coordonn\'ees sont nulles,
$$0=\int _{\Re(\sigma ')=\Re (\sigma )} D_{\lambda }(\ti {p}(\lambda )p_{\sigma '}(\lambda
)\psi_{\sigma '}(\lambda ))_{\vert\lambda =0}d_{\so _{\sigma ,M'}}\Im(\sigma ')=\kappa\ m_{\ul{0}}
D_{\lambda }(\ti {p}(\lambda )\psi_{\ul{0}} (\lambda ))_{\vert\lambda =0}.$$  Comme $m_{\ul{0}}$ est
constant, $m_{\ul{0}}\not\equiv 0$ implique $(D_{\lambda }(\ti{p}(\lambda ) \psi_{\ul{0}}(\lambda
)))_{\vert\lambda =0}=0$. Le r\'esultat pour les autres $\ul {k}$ se prouvent par le m\^eme argument,
en remplacant la fonction $\psi $ par $m_{\ul{k}}^{-1}\psi $ pour tout $\ul{k}$ tel que
$m_{\ul{k}}\not\equiv 0$.\hfill {\fin 2}

\null
\it Preuve: \rm (de la proposition) Supposons par absurde $D\ne 0$. Notons $D_{\ul{k}}$ un mon\^ome non
nul de degr\'e total maximal dans $D$. Par hypoth\`ese, il existe $\sigma'\in\o_{\sigma ,M'}$ tel que
$p_{\sigma '}(0)\ne 0$. Choisissons $\psi $ satisfaisant aux conclusions du premier lemme relatif \`a
$D_{\ul{k}}$ et $\sigma '$. Alors $(D_{\ul{k}}\psi_{\sigma '}(\lambda ))_{\vert\lambda =0}\ti
{p}(0)p_{\sigma '} (0)\ne 0$, et, par les r\`egles de Leibniz, $D_{\lambda }(\ti{p}(\lambda )p_{\sigma
'} (\lambda )\psi_{\sigma '}(\lambda ))_{\vert\lambda =0}=\ti {p}(0)p_{\sigma '}(0)$ 
$(D_{\ul{k}} \psi_{\sigma '}(\lambda ))_{\vert\lambda =0}\ne 0$. Ceci contredit les conclusions du
deuxi\`eme lemme, d'o\`u une contradiction, l'espace des fonctions holomorphes sur $a_{M,\sCC }^{M'*}$
\'etant invariant par translations. \hfill{\fin 2}

\null
\null {\bf 3.9 Proposition:} \it Soient $Res _A^{(1)}$ et $Res _A^{(2)}$ deux donn\'ees de r\'esidus en
$A$ telles que $Res _A^{(1)}\psi $ $=Res _A^{(2)}\psi $ pour tout $\psi\in\R (\o, \S _A)$. Alors $Res
_A^{(1)}=Res _A^{(2)}$.

\null Preuve: \rm Soit $\ul {d}\in\hbox{\ens N}^{\sS}$. Posons $M'=M_A$ et notons $D^{(1)}$ (resp.
$D^{(2)}$) un op\'erateur diff\'erentiel holomorphe \`a coefficients constants qui v\'erifie les
conclusions de la proposition {\bf 3.7} par rapport \`a $Res _A^{(1)}$ (resp. $Res _A^{(2)}$) et
$\ul{d}$. Il suffit de prouver que $D:=D^{(1)}-D^{(2)}=0$. 

Soit $\sigma\in A_{\sreg }$. Par hypoth\`ese, on a, avec $\ti{p}(\lambda )=\prod _{S\in\sS_A} (p_{\Re
(S)}^{d_S}(\lambda )p_S^{-d_S}(\sigma\otimes\chi  _{\lambda }))$ (cf. {\bf 3.4}), $D_{\lambda
}(\ti{p}(\lambda )\psi (\sigma\otimes \chi_{\lambda })) _{\vert\lambda =0}=0$ pour tout
$\psi\in\P(\o,\CC )$. Par les arguments de la preuve de la proposition {\bf 3.8}, ceci implique $D=0$.
\hfill{\fin 2}

\null
\null {\bf 3.10}   Supposons $\S $ stable par passage \`a la partie imaginaire, i.e. que $\o _{\sigma
,\alpha }\in\S $ implique $\o _{\Im _M(\sigma ),\alpha }\in\S $ (chaque \'el\'ement de $\H=\H(\S )$ est
alors le translat\'e d'un hyperplan affine de $\Re(\o )$ passant par l'origine $r(\o )$ de $\Re(\o )$
et appartenant \`a $\H$).

Avant d'\'enoncer le th\'eor\`eme des r\'esidus, nous allons choisir selon le proc\'ed\'e suivant, pour
tout $L\in \L=\L(\S )$, un ensemble de points de $L_{\sreg }$ qui est contenu dans la r\'eunion des
composantes connexes de $L_{\sreg }$ dont la cl\^oture contient l'origine $r(L)$ de $L$: Remarquons
d'abord que, l'ensemble $\L $ \'etant fini, il existe $\epsilon >0$ tel que pour tout $L\in \L$ la
boule $B(r(L),\epsilon )$ dans $r(\o )+a_M^*$ de centre $r(L)$ et de rayon $\epsilon $ n'ait une
intersection non vide avec un hyperplan affine $H\in\H$ que si celui-ci contient $r(L)$.

Fixons maintenant $L\in\L $  et notons $M_L$ le sous-groupe de Levi semi-standard qui lui est
associ\'e. Par nos hypoth\`eses, $r(\o )+a_{M_L}^* \in\L $. D\'esignons par $\H _0(r(\o )+a_{M_L}^*)$
l'ensemble form\'e des $H\in\H (r(\o )+a_{M_L}^*)$ avec $H$ passant par $r(\o )$. Notons $\P _{\sS}
(M_L)$ l'ensemble des composantes connexes de $a_{M_L,0}^*:=a_{M_L}^*-\bigcup _{H\in\sH _0(r(\so
)+a_{M_L}^*)}(-r(\o )+H)$.

Par choix de $\epsilon $, la translation de vecteur $r(L)$ dans $a_0^*$ envoie  $B(0,\epsilon )\cap
a_{M_L,0}^*$ dans $B(r(L),\epsilon )\cap L_{\sreg }$. Choisissons pour tout $Q\in\P_{\sS}(M_L)$ un
\'el\'ement $\epsilon _Q\in Q\cap B(0,\epsilon )$. Alors les points $r(L)+\epsilon _Q$, $Q\in
\P_{\sS}(M_L)$, se trouvent dans des composantes connexes de $L_{\sreg }$ dont la cl\^oture contient
$r(L)$.

\null   {\bf Th\'eor\`eme:} \it Soit $C$ une composante connexe de $\Re (\o )_{\sreg }$. Il existe pour
tout $A\in \A(\S )$ une donn\'ee de r\'esidu $\Res _A=\Res _A^C$ qui ne d\'epend que du choix de $C$,
telle que, pour toute fonction $\psi $ dans $\R (\o,\S )$ et tout $r\in C$, on ait
$$\displaylines {
\qquad\int _{\Re (\sigma )=r}\psi (\sigma )d_{\so }\Im (\sigma )=\hfill\cr
\hfill\sum _{L\in\sL }\sum _{Q\in\sP _{\ssS}(M_L)}\vert\P_{\sS} (M_L)\vert ^{-1} \sum _{A\in \sA(\sS
), \Re(A)=L}\int_{\Re (\sigma )=r(L)+\epsilon _Q}(\Res _A^C\psi )(\sigma )d_A\Im(\sigma ).\qquad\cr 
}$$  Par ailleurs, $\Res _{\so }=\id $ et, pour $S\in\S $, $\Res ^C_S$ est un multiple de $\Res _S$ par
un r\'eel (qui peut \^etre nul).

\null
\it Remarque: \rm On aurait \'evidemment pu prendre ci-dessus tout de suite la somme sur tous les
\'el\'ements $A$ de $\A (\S )$ (ce que l'on fera plus tard), mais \`a ce stade il nous a sembl\'e plus
utile d'ordonner les \'el\'ements de $\A(\S )$ par leur partie r\'eelle.

\null  
\it Preuve: \rm La preuve va suivre les lignes de celle du th\'eor\`eme 1.13 dans [BS], en utilisant
les r\'esultats pr\'eliminaires \'etablis dans les num\'eros pr\'ec\'edents. Rappelons que l'article
[BS] utilise une nouvelle approche inspir\'ee de [HO1] et cite [A], [MW1] et [L], la preuve de 
l'unicit\'e des donn\'ees de r\'esidu \'etant plus indirecte dans le travail de Langlands 
(cf. [MW1] V.2.2 Remarques (2) et V.3.13 Remarques (ii)).

Pour montrer l'existence, on va effectuer une r\'ecurrence sur la dimension $m$ de $a_M^*$. Si $m=0$,
il n'y a rien \`a montrer. Soit $m>0$ et supposons  l'existence \'etablie lorsque $\dim (a_M^*)<m$. On
peut \'ecrire
$$\displaylines {
\qquad\int _{\Re (\sigma )=r}\psi (\sigma )d_{\so }(\sigma )=\sum _{Q\in \sP_{\ssS}(M_{\sso })} \vert
\P_{\sS}(M_{\so })\vert ^{-1}\int _{\Re (\sigma )=r(\so )+\epsilon _Q}(\Res _{\so }\psi )(\sigma
)d_{\so }\Im (\sigma )\hfill\cr
\hfill +\sum _{Q\in\sP_{\ssS}(M_{\sso })}\vert\P_{\sS}(M_{\so })\vert ^{-1} (\int _{\Re (\sigma )=r}
\psi (\sigma )d_{\so }\Im (\sigma ) -\int _{\Re (\sigma )=r(\so )+\epsilon _Q}\psi (\sigma )d_{\so }\Im
(\sigma )). \cr  }$$  La premi\`ere partie est identique \`a la somme des termes de l'\'enonc\'e
correspondant \`a $L=\Re (\o )$. Les termes entre parenth\`eses (...) se d\'ecomposent en somme de
termes de la forme 
$$\int _{\Re (\sigma )=r_2}\psi (\sigma )d_{\so }\Im (\sigma )-\int _{\Re (\sigma )=r_1}\psi (\sigma
)d_{\so }\Im(\sigma )$$   avec $r_1$ et $r_2$ dans des composantes connexes adjacentes $C_1$ et $C_2$
de $\Re (\o )_ {\sreg }$.

Notons $H_{\alpha ,c}$ l'hyperplan de $\Re (\o )$ qui s\'epare $C_1$ et $C_2$. On applique la
proposition {\bf 3.6} \`a la diff\'erence de ces int\'egrales et on obtient que celle-ci est modulo une
constante r\'eelle \'egale \`a
$$\sum _S \int _{\Re (\sigma )=r_{12}} (\Res _S\psi )(\sigma ) d_S\Im(\sigma ),$$  la somme portant sur
les hyperplans $S\in\S $ de partie r\'eelle
$H_{\alpha ,c}$. Fixons un tel hyperplan $S$. Le sous-groupe de Levi associ\'e \`a $H_{\alpha ,c}$
\'etant $M_{\alpha }$, l'hypoth\`ese de r\'ecurrence s'applique \`a l'int\'egrale $\int _{\Re (\sigma )
=r_{12}}(\Res _S\psi )(\sigma ) d_S\Im(\sigma )$, et on obtient 
$$\displaylines {
\qquad\int _{\Re (\sigma )=r}\psi (\sigma )d_{\so }(\sigma )=\hfill\cr
\hfill\sum _{L\in\sL }\sum _{Q\in\sP _{\ssS}(L)}\sum _{A\in \sA(\sS ), \Re(A)=L}c(A) \int _{\Re
(\sigma )=r(L)+\epsilon _Q}(Res _A\psi )(\sigma )d_A\Im(\sigma ), \qquad\cr }$$  o\`u $c(A)$ d\'esigne
une constante r\'eelle qui peut \'eventuellement \^etre nulle. Mais peu importe, en multipliant $Res
_A$ par une constante convenable qui peut \^etre nulle, on obtient une expression du type voulu.

Prouvons maintenant l'unicit\'e. Il suffit de montrer que, si une famille $(Res_A)_{A\in\sA(\sS)}$ de
donn\'ees de r\'esidu v\'erifie
$$\sum _{L\in\sL }\sum _{Q\in\sP _{\ssS}(M_L)} \vert\P_{\sS}(M_L)\vert ^{-1} \sum _{A\in\sA(\sS);\Re
(A)=L}\int _{\Re (\sigma )=r(L)+\epsilon _Q}(Res _A\psi) (\sigma )d_A\Im(\sigma )=0$$ pour toute $\psi
\in \R(\o ,\S )$, alors $Res _A\equiv 0$ pour tout $A$.

On va effectuer une r\'ecurrence sur la dimension de $L=\Re (A)$. Soit $0\leq m\leq \dim (a_M^*)$ et
supposons $Res _A\equiv 0$ pour tout $A$ avec $\Re (A)$ de dimension $>m$. (Pour $m=\dim (a_M^*)$,
ceci est certainement vrai.) On a alors
$$0=\sum _{L\in\sL;\dim (L)\leq m}\sum _{Q\in\sP_{\ssS}(M_L)} \vert\P_{\sS}(M_L)\vert^{-1}\sum
_{A\in\sA(\sS ),\Re (A)=L} 
\int _{\Re (\sigma )=r(L)+\epsilon _Q} (Res _A\psi )(\sigma )d_A\Im(\sigma )$$ pour toute
$\psi\in\R(\o, \S)$. Fixons
$A\in\A(\S )$ avec $L:=\Re (A)$ de dimension $m$. D'apr\`es le lemme {\bf 3.3}, il existe pour tout
$A'\in\A (\S )$, $A'\ne A$, dont la partie r\'eelle est de dimension $\leq m$ un polyn\^ome $p_{A'}$
sur $\o $ qui est identiquement nul sur $A'$, mais qui n'est pas identiquement nul sur $A$. 

Par suite, il existe, pour tout $\ul{d}\in\hbox{\ens N}^{\sS}$, un polyn\^ome $p$ sur $\o $ qui n'est
pas identiquement nul sur $A$, telle que $Res_{A'}(p\psi )\equiv 0$ pour tout $A'\ne A$ et toute
$\psi\in\R(\o, \S,\underline{d})$. On en d\'eduit
$$0=\sum _{Q\in\sP_{\ssS}(M_L)}\vert\P_{\sS}(M_L)\vert ^{-1}\int _{\Re (\sigma ) =r(A)+\epsilon
_Q}(Res _A p\psi )(\sigma ) d_A\Im (\sigma ).$$ Ceci est en particulier v\'erifi\'e pour tout
$\psi\in\R (\o,\S_A, \ul{d}_{\vert\sS_{A }})$ (ici $\ul{d}_{\vert\sS_{A}}$ d\'esigne l'\'el\'ement de
$\hbox{\ens N}^{\sS _A}$ de composante $d_S$ en tout
$S\in\S _A$). Pour de telles $\psi $, $Res _A(p\psi )$ est une fonction r\'eguli\`ere sur
$A$, en sorte que l'on peut remplacer $r(L)+\epsilon _Q$ par $r(L)$  par le th\'eor\`eme de Cauchy pour
tout $Q\in\P _{\sS} (M_L)$. On obtient donc
$$0=\int _{\Re (\sigma )=r(A)} (Res_Ap\psi )(\sigma )d_A\Im (\sigma ).$$ Posons $M'=M_A$. D'apr\`es la
proposition {\bf 3.7}, il existe un op\'erateur diff\'erentiel holomorphe $D=D_{\lambda }$ sur
$a_{M,\sCC}^{M'*}$ \`a coefficients r\'eels constants, tel que 
$$(Res _A (p\psi ))(\sigma )=D_{\lambda }((\prod _{S\in\sS_A}p_{\Re (S)}^{d_S} (\lambda ))p_{\sigma
}(\lambda )\psi _{\sigma }(\lambda ))_{\vert\lambda =0}$$  pour tout $\psi \in
\R(\o,\S_{A},\ul{d}_{\vert \sS_{A}})$, d'o\`u, avec
$\ti{p}(\lambda )=\prod _{S\in\sS_A} (p_{\Re (S)}^{d_S}(\lambda )p_S^{-d_S} (\sigma\otimes\chi 
_{\lambda }))$,
$$0=\int _{\Re (\sigma )=r(A)} D_{\lambda }(\ti{p}(\lambda )p_{\sigma }(\lambda )\psi _{\sigma
}(\lambda ))_{\vert\lambda =0}d_A\Im(\sigma )$$  pour toute fonction polynomiale $\psi$ sur $\o $.

Comme $\ti{p}$ est une fonction m\'eromorphe sur $a_{M,\sCC}^{M'*}$ qui est r\'eguli\`ere et non nulle
en $0$, on peut appliquer la proposition {\bf 3.9} qui montre que $D=0$, d'o\`u $\Res _A\equiv
0$.\hfill {\fin 2}
 
\null
\null
\null
\vfill
\eject
{\bf 4. La fonction $\mu $ de Harish-Chandra, son syst\`eme de racines singuli\`eres et son
groupe de Weyl}

\null
\null  On se fixe dans cette section un sous-groupe parabolique semi-standard $P=MU$, et $\o $
d\'esignera l'orbite inertielle d'une repr\'esentation irr\'eductible cuspidale de $M$.

\null {\bf 4.1 Proposition:} \it i) Supposons $M$ maximal. Si $\mu $ n'est pas constante sur $\o $, il
existe une repr\'esentation unitaire $\sigma\in\o$ telle que $\mu (\sigma )=0$. Posons $\Sigma
_{red}(T_M)=\{\alpha, -\alpha \}$. Alors il existe une constante $c_{\alpha }> 0$ et des nombres
r\'eels strictement positifs $k:=k_{\alpha }$ et $l:=l_{\alpha }$ tels que l'on ait, avec
$\lambda\in\CC $ (et $\ti {\alpha }$ d\'efini dans ({\bf 3.2})), l'identit\'e de fonctions rationnelles
$$\mu (\sigma\otimes\chi_{\lambda\ti{\alpha }})=c_{\alpha }\ \mu_{\Re} (\sigma\otimes
\chi_{\lambda\ti{\alpha }})\mu _{\Im}(\sigma\otimes\chi_{\lambda\ti{\alpha }})$$
$$\mu_{\Re}(\sigma\otimes\chi_{\lambda\ti{\alpha }})={(1-q^{\lambda }) (1-q^{-\lambda })\over
(1-q^{k-\lambda })(1-q^{k+\lambda })}
\leqno{\hbox{\it avec}}$$
$$\mu_{\Im}(\sigma\otimes\chi_{({\pi i\over\log q}+\lambda )\ti{\alpha }})=
\cases{{(1-q^{\lambda })(1-q^{-\lambda })\over (1-q^{l-\lambda }) (1-q^{l+\lambda })}, & si
$\mu(\sigma\otimes\chi_{{\pi i\over\log q} \ti{\alpha }})=0$\cr 1& sinon.\cr}
\leqno{\hbox{\it et}}$$

ii) Soit $\sigma\in\o $ une repr\'esentation unitaire. L'ensemble des racines
$\alpha\in\Sigma_{red}(T_M)$ telles que
$\mu ^{M_{\alpha }}(\sigma )=0$ forme un syst\`eme de racines $\Sigma _{sp }(\sigma )$. Si $\alpha $,
$\beta $ sont des racines dans $\Sigma _{sp}(\sigma )$ de m\^eme longueur, alors $k_{\alpha }=k_{\beta
}$ et $l_{\alpha }=l_{\beta }$ (i.e. si l'un des deux nombres est d\'efini, l'autre l'est aussi, et les
deux nombres sont \'egaux).
\rm

\null
\it Preuve: \rm L'assertion (i) est le th\'eor\`eme 1.6 de [S2]. (Remarquons que la preuve de ce
th\'eor\`eme donn\'ee dans [S2] repose sur la simplicit\'e des p\^oles de $\mu ^{M_{\alpha }}$ dont la
preuve dans [S1] est erronn\'ee, comme cela nous a \'et\'e rapport\'e par A. Silberger. L'auteur
conna\^\i t une autre preuve de ce r\'esultat qui lui a \'et\'e communiqu\'e par J.-L. Waldspurger.)
Concernant l'assertion (ii), le fait que $\Sigma _{sp}(\sigma )$ est un syst\`eme de racines est
indiqu\'e dans [S3] 3.5. (La preuve est analogue \`a celle de la proposition {\bf 4.2} ci-apr\`es.) Les
autres assertions de (ii) r\'esultent alors du fait que le groupe de Weyl de $\Sigma _{sp}(\sigma )$
est inclus dans $W(M,\o)$ (cf. [S3] 3.5) et que la fonction $\mu $ est invariante par $W(M,\o )$ (cf.
{\bf 1.5}). \hfill{\fin 2}

\null {\it D\'efinition:} Notons $\Sigma _{sp}^+(\sigma )$ un ensemble de racines positives dans
$\Sigma _{sp}(\sigma )$ pour un certain ordre. On d\'efinit
$\mu _{sp,\sigma }=\prod _{\alpha\in\Sigma _{sp}^+(\sigma )}\mu_{\Re }^{M_{\alpha }}$ et $\mu
_{nsp,\sigma }=\mu\
\mu_{sp,\sigma }^{-1}$.

\null Remarquons que $\Im (\sigma )$ est un point r\'egulier pour $\mu _{nsp,\sigma }$ et que $\mu
_{nsp,\sigma }(\Im (\sigma ))\ne 0$.

\null
\null {\bf 4.2} Notons $\Sigma _{\so }$ l'ensemble des $\alpha \in \Sigma _{red}(T_M)$  tels qu'il
existe $w_{\alpha }\in W^{M_{\alpha }}(M,\o )$ v\'erifiant $w_{\alpha }(P\cap M_{\alpha })=\ol{P}\cap
M_{\alpha }$. (L'\'el\'ement $w_{\alpha }$ est n\'ecessairement unique (cf. [Cs] 7.1).)

\null {\bf Lemme:} \it L'automorphisme induit par $w_{\alpha }$ sur $a_M^*$ est une r\'eflexion
v\'erifiant
$w_{\alpha }(\alpha )=-\alpha $.

\null Preuve: \rm Comme $w_{\alpha }(P\cap M_{\alpha })w_{\alpha }^{-1}=\ol{P}\cap M_{\alpha }$, il est
clair que $w_{\alpha }(\alpha )=-\alpha $. L'espace $a_M^{M_{\alpha }*}$ \'etant de dimension $1$,
l'automorphisme induit par $w_{\alpha }$ sur $a_M^{M_{\alpha }*}$ est une r\'eflexion. Pour que
l'automorphisme induit par $w_{\alpha }$ sur $a_M^*$ soit une r\'eflexion, il reste \`a v\'erifier que
$w_{\alpha }$ agit trivialement dans
$a_{M_{\alpha }}^*$. Or, ceci est imm\'ediat, puisque $w_{\alpha }\in M_{\alpha }$ et que $a_{M_{\alpha
}}^*$ a une base form\'ee de caract\`eres de $M_{\alpha }$.\hfill {\fin 2}

\null
\null Notons $W_{\so }$ le sous-groupe de $W(M,\o )$ engendr\'e par les r\'eflexions $w_{\alpha }$ avec
$\alpha\in\Sigma _{\so }$. La proposition suivante est analogue \`a celle sur les racines $\omega
$-sp\'eciales dans le num\'ero 3.5 dans [S3].

\null {\bf Proposition:} \it L'ensemble $\Sigma _{\so }$ est un syst\`eme de racines dans un sous-espace
$a_{M,\so}^*$ de $a_M^*$. Il est invariant par $W(M,\o )$. Son groupe de Weyl est engendr\'e par les
restrictions des r\'eflexions $w_{\alpha }$ de $a_M^*$ \`a ce sous-espace. Il est isomorphe \`a $W_{\so
}$ au moyen de cette restriction. Le groupe $W_{\so }$ est distingu\'e dans $W(M,\o )$.

\null Preuve: \rm On va appliquer les propositions 8 et 9 du paragraphe 2, chapitre VI de [B]. Notons
$a_{M,\so}^*$ le sous-espace de $a_M^*$ engendr\'e par $\Sigma _{\so }$. Remarquons d'abord que $\Sigma
_{\so }$ est stable pour l'action par $W(M,\o)$, puisque $w^{-1}(P\cap M_{\alpha })w=wPw^{-1}\cap
M_{w\alpha }$ et que $ww_{\alpha }w^{-1}(wPw^{-1}\cap M_{w\alpha })= (w\ol{P}w^{-1}\cap M_{\alpha })$.
Notons $W_{\so }'$ l'image de $W_{\so }$ dans $\GL (a_{M,\so }^*)$ et $w'$ l'image d'un \'el\'ement
$w\in W_{\so }$ par cette application.

Par ce qui pr\'ec\`ede, le groupe $W_{\so }'$ est engendr\'e par des r\'eflexions. Il est essentiel,
i.e. aucun
\'el\'ement $\ne 0$ de $a_{M,\so }^*$ n'est invariant par $W_{\so }'$: comme $\Sigma _{\so }$ engendre
$a_{M,\so }^*$ et que $v$ invariant par $W_{\so }$ implique $\sum _{w\in W_{\sso }} wv\in \hbox{\ens
R}^*v$, il suffit de v\'erifier que $\Sigma _{w\in W_{\sso }}w\alpha =0$ pour tout $\alpha\in\Sigma
_{\so }$. Or, on a
$\sum _{w\in W_{\sso }}w\alpha =\sum _{w\in W_{\sso }}w(w_{\alpha }\alpha )=-\sum _{w\in W_{\sso
}}w\alpha $, d'o\`u
$\sum _{w\in W_{\sso }}w\alpha =0$.

Notons $L_{M,\so }$ le $\hbox{\ens Z}$-sous-module de $a_{M,\so }^*$ engendr\'e par $\Sigma _{\so }$.
Il est sans torsion et de type fini, donc libre de m\^eme rang que $a_{M,\so }^*$. On a vu que
$L_{M,\so }$ est invariant par $W_{\so }'$. La proposition 9 (loc. cit.) prouve donc que $W_{\so }'$
est le groupe de Weyl d'un syst\`eme de racines dans $a_{M,\so }^*$. Il r\'esulte alors des
propositions 8 et de la preuve de l'assertion ((iv) implique (ii)) de la proposition 9 (loc. cit.) que
$\Sigma _{\so }$ est un syst\`eme de racines de groupe de Weyl $W_{\so }'$.

Il reste \`a voir que la restriction $W_{\so }\rightarrow W_{\so }'$ est injective. Soit $w\in W_{\so
}$ d'image triviale dans $W_{\so }'$. Alors $w\alpha =\alpha $ pour tout $\alpha\in\Sigma _{\so }$.
Soit $M'$ un sous-groupe de Levi minimal qui contient $\Sigma _{\so }$. Alors, $w\in M'$ et $\Sigma
_{\so }$ engendre $a_M^{M'*}$. Comme $w$ op\`ere trivialement sur $\Sigma _{\so }$, il s'en suit que
$w=1$. \hfill{\fin 2}

\null
\null {\bf 4.3} Notons $\S _{\mu }$ le plus petit ensemble qui contient les hyperplans affines de $\o $
qui sont singuliers pour $\mu $ et qui est stable par passage \`a la partie imaginaire, i.e. que $\o
_{\sigma ,\alpha }\in\S_{\mu }$ implique $\o _{\Im(\sigma ),\alpha }\in\S _{\mu }$.

\null {\bf Proposition:} \it L'ensemble $\S _{\mu }$ est $W(M,\o )$ invariant. L'ensemble des
$\alpha\in\Sigma (T_M)$ tels qu'il existe $\sigma\in\o $ avec $\o _{\sigma ,\alpha }\in\S _{\mu }$ est
contenu dans $\Sigma _{\so }$. 

\null Preuve: \rm La premi\`ere assertion r\'esulte de l'invariance de la fonction $\mu $ par $W(M,\o
)$. La deuxi\`eme assertion est due \`a Harish-Chandra (cf. [S1] corollaire 5.4.2.3 et 5.4.2.2).

\null  {\bf 4.4} Pour all\'eger les notations, on \'ecrira dans la suite de ce paragraphe plus
simplement $\S $ au lieu de $\S _{\mu }$.

\null {\bf Lemme:} \it Soit $\S _0$ le sous-ensemble de $\S $ form\'e des hyperplans affines dont la
partie r\'eelle contient l'origine de $a_M^*$. L'op\'erateur rationnel $\sigma\mapsto\mu (\sigma
)J_{P\vert\ol{P}}(\sigma )$ est r\'egulier sur $\o -\bigcup _{S\in\sS -\sS_0} S$.

\null Preuve: \rm Par la formule du produit pour la fonction $\mu $ et pour les op\'erateurs
d'entrelace-ment, on se ram\`ene au cas d'un sous-groupe parabolique maximal. Si $\Re (\sigma )>_P0$ ou
$-\Re (\sigma )>_P0$, alors $J_{P\vert\ol{P}}$ est r\'egulier en la repr\'esentation cuspidale $\sigma
$ par les r\'esultats de Harish-Chandra (cf. [S1] Th. 5.4.2.1). Dans les autres cas, la r\'egularit\'e
en $\sigma $ de l'application dans l'\'enonc\'e r\'esulte de [W] Corollaire V.2.3.\hfill{\fin 2}

\null
\null {\bf 4.5}  Notons $\Delta _{\so }$ la base de $\Sigma _{\so }$ qui est form\'ee de racines qui
sont positives pour $P$. Si $\Omega\subseteq\Delta _{\so }$, d\'esignons par $M_{\Omega }$ le plus
petit sous-groupe de Levi semi-standard $M'$ contenant $M$ tel que $\Omega\subseteq a_M^{M'*}$. Comme
$M_{\Omega }$ peut \^etre obtenu en adjoignant successivement des racines de $\Omega $ \`a $M$, $\Omega
$ forme une base de $a_M^{M'*}$. Un sous-groupe de Levi $M'$ de $G$ qui contient $M$ sera dit \it
$(P,\S)$-standard, \rm si $M'=M_{\Omega }$ pour un $\Omega\subseteq \Delta _{\so }$. On \'ecrira
$\Omega =
\Delta _{\so }^{M'}$.

Pour un sous-groupe de Levi $(P,\S )$-standard $M'$ on pose alors

$$\eqalign {  W_{\so ,M'}^+=\ &\{ w\in W_{\so }\vert w\ \hbox{\rm de longueur minimale dans $wW_{\so
}^{M'}$ } \},\cr W _{\Delta _{\sso }}=\ &\{w\in W(M,\o )\vert w\Delta _{\so }=\Delta _{\so }\},\cr W
_{\Delta _{\sso },M'}=\ &\hbox{\rm syst\`eme de repr\'esentants de $W _{\Delta _{\sso }}W^{M'}(M,\o
)/W^{M'}(M,\o )$, form\'e d'\'el\'e- }\cr &\hbox{\rm ments de $W _{\Delta _{\sso }}$ et contenant
l'identit\'e de $W(M,\o )$},\cr  W_{M'}^+(M,\o )=\ &W _{\Delta _{\sso },M'}W_{\so ,M'}^+.\cr   }$$

\null
\it Remarque: \rm $W _{\Delta _{\sso }}$ est un sous-groupe de $W(M,\o )$ qui v\'erifie $W_{\so }\cap W
_{\Delta _{\sso }}=\{1\}$. Il existe bien des cas o\`u $W _{\Delta _{\sso }}\ne\{1 \}$: les syst\`emes
de racines de type $A_n$, $D_n$ et $E_6$ poss\`ede par exemple des automorphismes ne venant pas d'un
\'el\'ement de leur groupe de Weyl (cf. [B]) et qui peuvent bien \^etre r\'ealis\'es par un \'el\'ement
du groupe de Weyl d'un syst\`eme de racines plus grand de m\^eme rang. 

\null  {\bf Lemme:} \it L'ensemble $W_{M'}^+(M,\o )$ est un syst\`eme de repr\'esentants de  $W(M,\o
)/W^{M'}$ $(M,\o )$ dans $W(M,\o )$.

\null Preuve: \rm Soit $w\in W(M,\o )$. Alors il existe $w_{\so }\in W_{\so }$, tel que $ww_{\so
}^{-1}$ v\'erifie $ww_{\so }^{-1}\Delta _{\so }\subseteq\Sigma (P)$. Comme $\Sigma _{\so }$ est
invariant par $W(M,\o )$, il en suit que $w _{\Delta }:=ww_{\so }^{-1}\in W _{\Delta _{\sso }}$.
\'Ecrivons $w _{\Delta }=w _{\Delta ,M'}w _{\Delta }^{M'}$ avec $w _{\Delta ,M'}\in W _{\Delta _{\sso
} ,M'}$ et $w _{\Delta }^{M'}\in W^{M'}(M,\o )$. Alors $w=w _{\Delta ,M'}w _{\Delta }^{M'}w_{\so }=w
_{\Delta ,M'}(w _{\Delta }^{M'}w_{\so }w _{\Delta } ^{M'\ -1})w _{\Delta }^{M'}$. Comme $w _{\Delta
}^{M'}w_{\so }w _{\Delta }^{M'\ -1}\in W_{\so }$, il est clair que tout \'el\'ement de $W(M,\o
)/W^{M'}(M,\o )$ est repr\'esent\'e par un \'el\'ement de $W_{M'}^+(M,\o )$.

Prouvons l'unicit\'e de cet \'el\'ement: soient $w _{\Delta }, w _{\Delta }'\in W _{\Delta _{\sso }
,M'}$ et $w, w'\in W_{\so ,M'}^+$ et supposons que $w _{\Delta }w$ et $w _{\Delta }'w'$ soient dans la
m\^eme classe \`a droite modulo $W^{M'}(M,\o )$. Il existe donc $w''\in W^{M'}(M,\o )$ tel que
$w_{\Delta }w=w_{\Delta }'w'w''$. Alors $w _{\Delta }'w'w''w^{-1}\in W_{\Delta _{\sso }}$, ce qui
implique $w'w''w^{-1}\in W_{\Delta _{\sso }}$. Comme $w$ et $w'$ sont de longueur minimale dans leur
classe \`a droite modulo $W^{M'}_{\so }$, $w''$ doit   v\'erifier $w''\Delta _{\so }^{M'}=\Delta _{\so
}^{M'}$ et donc \^etre un \'el\'ement de $W_{\Delta _{\sso }}$. On en d\'eduit $w''^{-1}w'w''\in W_{\so
,M'}^+$, $W_{\so }$ \'etant distingu\'e dans $W(M,\o )$, et $(w''^{-1}w'w'')w^{-1}\in W_{\Delta _{\sso
}}$. Comme d'autre part $(w''^{-1}w'w'') w^{-1}\in W_{\so }$, c'est l'\'el\'ement neutre de ce groupe,
d'o\`u $w_{\Delta }=w_{\Delta }'w''$. Mais, comme $w_{\Delta }$ et $w_{\Delta }'$ sont des
repr\'esentants d'une m\^eme classe \`a droite modulo $W^{M'}(M,\o )$, il en r\'esulte $w''=1$, d'o\`u
l'\'egalit\'e recherch\'ee. \hfill{\fin 2}

\null
\null {\bf 4.6} Continuons \`a supposer $M'$ $(P,\S )$-standard. Il existe alors $A\in \A(\S )$ tel que
$M'=M_A$. La composante connexe $Q$ de $\P_{\sS}(M_A)$ qui contient les \'el\'ements $\lambda $ de
$a_{M_A}^*$ avec $\langle\alpha^{\vee },\lambda \rangle>0$ pour tout $\alpha\in\Delta _{\so } -\Omega
$ sera dite \it $(P,\S )$-standard. \rm On la notera $P_A$. 

\null  {\bf Lemme:} \it Soit $A\in\A (\S)$. Pour tout \'el\'ement $w$ de $W_{\so }$, l'action par $w$
sur $a_M^*$ induit une application de $\P_{\sS }(M_A)$ sur $\P_{\sS }(M_{wA})$. Soit
$Q\in\P_{\sS}(M_A)$. Il existe un unique sous-groupe de Levi $(P,\S)$-standard $M'$ de $G$ tel que $Q$
soit conjugu\'e \`a l'\'el\'ement $(P,\S)$-standard de $\P_{\sS}(M')$ par un \'el\'ement de $W_{\so }$.
La classe \`a droite modulo $W_{\so }^{M_A}$ de cet \'el\'ement est d\'etermin\'ee de fa\c con unique.

\null  Preuve: \rm Soit $M_1$ le sous-groupe de Levi contenant $M$, tel que $\Delta _{\so }$ engendre
$a_M^{M_1*}$, i.e. $a_M^{M_1*}=a_{M,\so }^*$. On a $M_A\subseteq M_1$ et tout \'el\'ement de $\P
_{\sS}(M_A)$ est uniquement d\'etermin\'e par son intersection avec $a_M^{M_1*}$. L'intersection de
$a_{M_A}^*$ avec $a_M^{M_1*}$ \'etant $a_{M_A}^{M_1*}$, on se retrouve dans le contexte des chambres
relevant d'un syst\`eme de racines, o\`u le r\'esultat est connu (cf. [Ca] Th. 2.5.8 et prop. 2.6.1,
2.6.3). \hfill{\fin 2}

\null
\null {\bf 4.7} Posons, pour $M'$ un sous-groupe de Levi de $G$ contenant $M$, $\W (M',\o)=\{w\in
W(M,\o )\vert wM'=M'\}$ et notons $W(M', \o)$ le groupe quotient $\W (M',\o )/W^{M'}(M,\o )$.

\null {\bf Lemme:} \it Soit $M'$ un sous-groupe de Levi $(P,\S )$-standard. Le nombre de sous-groupes
de Levi $(P,\S )$-standard de $G$ conjugu\'es \`a $M'$ par un \'el\'ement de $W(M,\o )$ est \'egal \`a
${\vert\sP_{\ssS}(M')\vert\ \vert W _{\Delta _{\sso } ,M'}\vert \over\vert W(M',\so)\vert}$.

\null  Preuve: \rm   Soit $A\in \A (\S )$ avec $M'=M_A$, et notons $[A]_{P,\sS}$ l'ensemble des
$A'\in\A(\S )$ qui sont conjugu\'es \`a $A$ par un \'el\'ement de $W(M,\o )$ et qui sont tels que
$M_{A'}$ soit $(P,\S )$- standard. On voit que $A'\in [A]_{P,\sS}$, si et seulement si $M_{A'}$ est
$(P,\S )$-standard et qu'il existe $w\in W_{M_A}^+(M,\o )$ tel que $A'=wA$. \'Ecrivons $w=w_{\Delta
}w_{\so }$ avec $w_{\Delta }\in W _{\Delta _{\sso },M_A}$ et $w_{\so }\in W_{\so ,M_A}^+$. Comme $W
_{\Delta _{\sso }}$ envoie tout sous-groupe de Levi $(P,\S )$-standard sur un sous-groupe de Levi $(\P
,\S )$-standard, il faut que $w_{\so }A\in [A]_{P,\sS}$. Les \'el\'ements de $[A]_{P,\sS}$ sont donc de
la forme $w_{\Delta }w_{\so }A$ avec $w_{\Delta }\in W _{\Delta _{\sso },M_A}$ et $w_{\so }\in W_{\so
,M_A}^+$ tel que $w_{\so } A\in [A]_{P,\sS }$. Notons $W_{\so ,M_A}^+(P, \S)$ le sous-ensemble de
$W_{\so ,M_A}^+$ d\'ecrit par cette propri\'et\'e.

Par le lemme {\bf 4.6}, on d\'efinit une application de $W_{\so ,M_A}^+(P,\S )$ dans $\P_{\sS }(M_A)$,
en associant \`a un \'el\'ement $w$ l'\'el\'ement $w^{-1}P_{wA}$ de $\P _{\sS }(M_A )$, o\`u $P_{wA}$
d\'esigne l'\'el\'ement $(P,\S )$-standard de $\P _{\sS }(M_{wA})$. Le lemme {\bf 4.6} prouve
\'egalement que cette application est surjective et qu'elle est injective. Par suite, $\vert W_{\so
,M_A}^+(P,\S )\vert =\vert \P_{\sS }(M_A)\vert $. On en d\'eduit que le sous-ensemble $W_{M_A}^+(M,\o
)_{P, \sS}$ de $W_{M_A}^+(M,\o )$ form\'e des \'el\'ements $w$ tels que $wA\in [A]_{P,\sS }$ est de
cardinal $\vert \P_{\sS }(M_A)\vert\ \vert W _{\Delta _{\sso },M_A }\vert $. Il reste \`a voir que les
fibres de l'application $W_{M_A}^+(M,\o )\rightarrow \{M_{A'}\vert A'\in [A]_{P,\sS }\}$, $w\mapsto
wM_A$, sont tous en bijection avec  $W(M_A, \o)$.

Fixons $w\in W_{M_A}^+(M,\o )$ et consid\'erons la fibre au-dessu de $wM_A$. Si $w'M_A=wM_A$, alors il
existe un unique \'el\'ement $w''\in W(M,\o )$ avec $w''M_A=M_A$ tel que $ww''=w'$, d'o\`u une
application de la fibre au-dessus de $wM_A$ dans $W(M_A, \o)$. Inversement, soit $w'\in W(M_A,\o )$.
Alors il existe un unique repr\'esentant $w''$ de $w'$ dans $W(M,\o )$ tel que $ww''\in W_{M_A}^+(M,\o
)$, et on a $ww''M_A=wM_A$. Par construction, ces deux applications sont inverses l'une de l'autre.
\hfill{\fin 2}

\null
\null {\bf 4.8} On va maintenant appliquer le th\'eor\`eme de r\'esidus du paragraphe 3 au syst\`eme
fini d'hyperplans affines $\S=\S _{\mu }$ . Rappelons (cf. {\bf 1.4}) que l'on a not\'e, pour $A\in
\A(\S )$, par $M_A$ le sous-groupe de Levi semi-standard de $G$ associ\'e \`a $A$. Il contient $M$.
Fixons $\psi\in\R(\o ,\S )$.

\null {\bf Lemme:} \it Pour tout $w\in W(M,\o )$, on a $(Res _{wA}^{wC}\psi)\circ w=Res_A^C(\psi\circ
w)$. 

\null Preuve: \rm Soit $r$ un point de $C$. Alors $wr\in wC$. Comme $w$ est une transformation qui
pr\'eserve les mesures, on a 
$$\int _{\Re (\sigma )=r}(\psi\circ w)(\sigma )d_{\so }\Im (\sigma )=\int _{\Re (\sigma )=wr}\psi
(\sigma )d_{\so }\Im (\sigma ).$$  L'ensemble $\S $ (et donc $\A (\S )$) \'etant invariant par
$w^{-1}$, on trouve par le th\'eor\`eme {\bf 3.10}
$$\displaylines {
\qquad\int _{\Re (\sigma )=r}(\psi\circ w)(\sigma )d_{\so }\Im (\sigma )=\sum _{L\in\sL }\sum
_{Q\in\sP _{\ssS}(M_L)}\vert\P _{\sS}(M_L)\vert ^{-1}\hfill\cr
\sum _{A\in \sA(\sS _{\mu }), \Re(A)=L}\int _{\Re (\sigma )=w^{-1}(r(L)+\epsilon _Q)}(Res
_{w^{-1}A}^C(\psi\circ w))(\sigma )d_{w^{-1}A}\Im(\sigma ).\cr  }$$  et
$$\displaylines {
\qquad\int _{\Re (\sigma )=wr}\psi(\sigma )d_{\so }\Im (\sigma )=\hfill \cr
\hfill\sum _{L\in\sL }\sum _{Q\in\sP_{\ssS}(L)}\vert\P _{\sS}(M_L)\vert ^{-1}\sum _{A\in \sA,
\Re(A)=L}\int _{\Re (\sigma )=r(L)+\epsilon _Q}(Res _{A}^{wC}\psi )(\sigma ) d_{A}\Im(\sigma
).\qquad\cr  }$$  Comme
$$\displaylines{
\qquad\int _{\Re (\sigma )=r(L)+\epsilon _Q}(Res_A^{wC}\psi )(\sigma )d_A\Im (\sigma )=\hfill\cr
\hfill\int _{\Re (\sigma )=w^{-1}(r(L)+\epsilon _Q)} ((Res _A^{wC}\psi )(w\sigma )d_{w^{-1}A}\Im
(\sigma ),\qquad\cr }$$  le lemme r\'esulte de l'unicit\'e des donn\'ees de r\'esidu relatives \`a $C$
et de l'invariance de $\A(\S )$ par les \'el\'ements de $W(M,\o )$ (cf. {\bf 4.3}). \hfill {\fin 2}

\null
\null {\bf 4.9} Rappelons que $\Re (\o )=a_M^*$. Notons $\Res _A^{P}$ la donn\'ee de r\'esidu en $A$
correspondant \`a la composante connexe de $\R (\o )_{\sreg }$ (cet espace a \'et\'e d\'efini en {\bf
1.4}) contenant les $\lambda \in a_M^*$ avec $\lambda\gg _P0$. Rappelons (cf. {\bf 4.3}) que $\S _{\mu
^{M_A}}$ d\'esigne le plus petit ensemble qui contient les hyperplans affines de $\o $ qui sont
singuliers pour la fonction $\mu ^{M_A}$ d\'efinie relative \`a $M_A$ et qui est stable par passage \`a
la partie imaginaire.

\null {\bf Lemme:} \it Soit $A\in \A (\S )$. La donn\'ee de r\'esidu $\Res _A^{P\cap M_A}$ d\'efinie
sur $\R (\o,\S _{\mu ^{M_A}})$ se prolonge \`a $\R (\o,\S )$ et y est \'egale \`a $\Res _A^P$.

\null Preuve: \rm Il est clair par d\'efinition d'une donn\'ee de r\'esidu que $\Res _A^{P\cap M_A}$ se
prolonge
\`a $\R (\o,\S )$. L'\'egalit\'e avec $\Res _A^P$ r\'esulte alors de l'unicit\'e des donn\'ees de
r\'esidus et de {\bf 3.9}: les deux donn\'ees de r\'esidus co\"\i ncident pour $\psi\in\R(\o,\S _A)$,
puisque $\int _{\Re (\sigma )=r\gg _P0}\psi (\sigma )d_A\Im(\sigma ) =\int _{\Re (\sigma )=r_1\gg
_{P\cap M_A}0}\psi (\sigma )d_A\Im(\sigma )$ pour de tels $\psi $.
\hfill {\fin 2}

\null
\null {\bf 4.10 Lemme:} \it Soit $A\in\A (\S )$ et posons $\Sigma _{\so }^{M_A} (P\cap M_A)= \Sigma
_{\so }^{M_A}\cap\Sigma ^{M_A}(P\cap M_A)$. Si $w\in W(M,\o )$ v\'erifie $w\Sigma _{\so}^{M_A}(P\cap
M_A)\subseteq\Sigma (P)$, alors $Res _A^{P}(\psi )=Res_A^{w^{-1}P}(\psi )$.

\null  Preuve: \rm Comme une donn\'ee de r\'esidu en $A$ est uniquement d\'etermin\'ee par ses valeurs
sur $\R (\o,\S _A)$ (cf. {\bf 3.9}), il suffit de montrer que, si $r \in a_M^*$ v\'erifie $r\gg _P0$,
alors $w^{-1}r$ est dans la m\^eme composante connexe de $\Re (\o )_{\sS _A,\sreg }$ que $r$. Or, si
$r\gg _P0$ et $w\Sigma _{\so } ^{M_A}(P\cap M_A)\subseteq \Sigma(P)$, on a bien $\langle(w\alpha
)^{\vee },r\rangle\gg _P0$ pour tout $\alpha\in\Sigma _{\so }^{M_A} (P\cap M_A)$. \hfill{\fin 2}

\null
\null  {\bf 4.11} Le r\'esultat suivant est une g\'en\'eralisation \`a notre situation du r\'esultat
principal de [BS].

\null {\bf Proposition:} \it Notons $\Delta _{\so }$ la base de $\Sigma _{\so }$ qui est form\'ee de
racines qui sont positives pour $P$. Alors on a
$$\displaylines {
\qquad\int _{\Re (\sigma )=r\gg _P0}\psi(\sigma )d_{\so}\Im (\sigma )=\sum _{\Omega\subseteq\Delta
_{\sso }}\sum _{A\in\sA(\sS ), M_A=M_{\Omega }}\hfill\cr
\hfill\sum _{w\in W_{M_A}^+(M,\so)} \vert W _{\Delta _{\sso },M_A}\vert^{-1} \vert \P_{\sS }(M_A)\vert
^{-1}\int _{\Re(\sigma )=r(A)+\epsilon _{P_A}} Res^P_A(\psi\circ w)(\sigma )d_A\Im (\sigma ).\qquad\cr
}$$  

\null  Preuve: \rm  Par {\bf 4.6}, il existe pour tout $A\in\A (\S )$ et tout $Q$ dans $\P_{\sS }
(M_A)$ un unique sous-groupe de Levi $(P,\S )$-standard $M'$ de $G$ tel que $Q$ soit conjugu\'e \`a
l'\'el\'ement $(P,\S )$-standard de $\P_{\sS }(M')$ par un \'el\'ement $w$ de $W_{\so }$. La classe
modulo $W_{\so }^{M_A}$ de $w$ est d\'etermin\'ee de fa\c con unique. Remarquons que $\A (\S )$ est
$W(M,\o )$-invariant (cf. {\bf 4.3}). On d\'efinit donc une bijection entre l'ensemble des couples
$(A,Q)$ form\'e d'un \'el\'ement $A$ de $\A (\S )$ et d'un \'el\'ement $Q$ de $\P _{\sS }(M_A)$ avec
l'ensemble des couples $(A,w)$ form\'e d'un \'el\'ement $w$ de $W_{\so ,M_A}^+$ et d'un \'el\'ement $A$
de $\A (\S )$ avec $M_A$ $(P,\S )$-standard, en associant \`a $(A,w)$ le couple $(wA,wP_A)$, $P_A$
d\'esignant la composante $(P,\S )$-standard de $\P _{\sS }(M_A)$. De plus, on voit que l'on peut
choisir les $\epsilon _Q$ dans la proposition {\bf 3.10} tels que $w\epsilon _{P_A}=\epsilon _{wP_A}$
pour de tels $w$.

Avant d'appliquer ceci \`a la proposition {\bf 3.10}, ajoutons que, pour $w\in W_{\Delta _{\sso },M_A}$
et $A\in \A(\S )$ avec $M_A$ $(P,\S )$-standard, le sous-groupe de Levi $wM_A=M_{wA}$ est \'egalement
$(P,\S )$-standard. Comme $wP_A=P_{wA}$ pour de tels $w$, on peut en outre choisirs les $\epsilon
_{P_A}$ tels que $w\epsilon _{P_A}=\epsilon _{P_{wA}}$. Il r\'esulte alors de ces remarques et de la
proposition {\bf 3.10}, ayant $W_{M_A}^+(M,\o )=W _{\Delta _{\sso }, M_A}W_{\so ,M_A}^+$ par {\bf 4.5},
que
$$\displaylines {
\qquad\int _{\Re (\sigma )=r\gg _P0}\psi (\sigma )d_{\so}\Im (\sigma )=\sum _{\Omega\subseteq\Delta
_{\sso }}\sum _{A\in\sA(\sS), M_A=M_{\Omega }}\hfill\cr
\hfill\sum _{w\in W_{M_A}^+(M,\so )}\vert W _{\Delta _{\sso },M_A}\vert^{-1}\ \vert
\P_{\sS}(M_A)\vert^{-1}\ \int _{\Re(\sigma )=w(r(A)+\epsilon _{P_A})}(Res ^P_{wA} \psi )(\sigma
)d_{wA}\Im (\sigma ).\qquad\cr  }$$  Il ne reste alors qu'\`a remarquer que 
$$\int _{\Re (\sigma )=w(r(A)+\epsilon _{P_A})} (Res^P_{wA}\psi )(\sigma )d_{wA} \Im (\sigma )=\int
_{Re (\sigma )=r(A)+\epsilon _{P_A}}(Res_{wA}^P\psi)(w\sigma )d_A\Im(\sigma ),$$  et que, pour $w\in
W_{M_A}^+(M,\o )$, $(Res^P_{wA}\psi )\circ w=Res ^{(w^{-1}P)} _A(\psi\circ w)=Res^P_A(\psi\circ w)$ par
les lemmes de {\bf 4.8} et de {\bf 4.10}, pour en d\'eduire l'expression dans l'\'enonc\'e de la
proposition. \hfill{\fin 2}

\null
\null {\bf 4.12} En g\'en\'eralisant les m\'ethodes de [BS] on pourra probablement montrer que $\Res
_A^P= 0$, si
$r(A)$ n'est pas une combinaison lin\'eaire d'\'el\'ements de $\Delta _{\so }^{M_A}:=\Delta _{\so
}\cap\Sigma^{M_A}(A_M)$ \`a coefficients $\geq 0$.

\null
\null
\null
\vfill
\eject
{\bf 5. Les d\'eriv\'ees d'un coefficient matriciel}

\null
\null  {\bf Proposition:} \it Soit $P=MU$ un sous-groupe parabolique semi-standard de $G$ et
$D=D_{\lambda }$ un op\'erateur diff\'erentiel holomorphe \`a coefficients constants sur
$a_{M,\sCC}^{G*}$. Soit $(\sigma ,E)$ une repr\'esentation irr\'eductible lisse de $M$ et $\tau $ dans
$i_{P\cap K}^KE\otimes i_{P\cap K}^KE^{\vee }$.

Alors il existe une famille finie $\{(\ti{\sigma }_j,\ti {E}_j)\}_{j\in J}$ de repr\'esentations
admissibles lisses de $M$ de longueur finie et de sous-quotients irr\'eductibles isomorphes \`a
$(\sigma ,E)$ et une famille
$\{\ti {\tau }_j\}_{j\in J}$ d'\'el\'ements de $i_{P\cap K}^K\ti{E}_j\otimes i_{P\cap K}^K\ti{E}^{\vee
}_j$, telles que, pour $\lambda _0\in a_{M,\sCC}^{G*}$,
$$(D_{\lambda }E_{P,\sigma\otimes\chi_{\lambda}}^G(\tau))_{\vert\lambda =\lambda_0}  =\sum _{j\in
J}E_{P,\ti{\sigma}_j\otimes\chi _{\lambda _0}}^G(\ti {\tau }_j).$$  Par ailleurs, les restrictions \`a
$T_G$ des repr\'esentations $i_P^G\ti{\sigma }_j$, $j\in J$, et $i_P^G\sigma $ agissent sur leurs
espaces respectifs par le m\^eme caract\`ere qui est \'egal \`a la restriction \`a $T_G$ du caract\`ere
central de $\sigma $. \rm

\null  Avant de prouver le th\'eor\`eme, \'etablissons le lemme suivant:  

\null  {\bf Lemme:} \it Soit $(\sigma ,E)$ une repr\'esentation irr\'eductible lisse de $M$, $p_1$ un
polyn\^ome sur $a_{M,\sCC}^G$ et $e_1\otimes e_1^{\vee }$ dans $E\otimes E^{\vee }$. Alors
l'application $m\mapsto p_1(H_M^G(m))
\langle\sigma (m)e_1,e_1^{\vee }\rangle $ est le coefficient matriciel d'une repr\'esentation
admissible de longueur finie $\ti{\sigma }$ de $M$ dont les sous-quotients sont isomorphes \`a $\sigma
$.

\null  Preuve: \rm Notons $\ti {E}$ l'espace vectoriel engendr\'e par les applications lisses de $M$
dans $E$ de la forme $\ti {e}_p: m\mapsto p(H_M^G(m)) \sigma (m)e$ avec $e\in E$ et $p$ polyn\^ome sur
$a_{M,\sCC}^G$. Faisons agir $M$ sur $\ti {E}$ par translation \`a droite et notons $\ti{\sigma }$ la
repr\'esentation lisse de $M$ qui correspond \`a cette action. 

Notons $l$ la dimension de $a_{M,\sCC}^G$. Munissons $\hbox{\ens N}^l$ de l'ordre lexicographique.
Fixons une base de $a_{M,\sCC}^G$. On a alors une notion de mon\^ome sur $a_{M,\sCC}^G$ correspondant
aux coordonn\'ees par rapport \`a la base choisie. Lorsque $p$ est un tel mon\^ome, notons
$\ul{\deg}(p)$ le vecteur dans $\hbox{\ens N}^l$ dont les composantes sont les degr\'es partiels de $p$
par rapport aux diff\'erentes coordonn\'ees. On l'appellera le degr\'e de $p$. Lorsque $p$ est un
polyn\^ome, on d\'esignera par $\ul{\deg }(p)$ le degr\'e du mon\^ome non nul dans $p$ de degr\'e
maximal pour l'ordre lexicographique.

Pour $\ul {d}\in\hbox{\ens N}^l$, soit $\ti {E}_{\leq\ul{d}}$ (resp. $\ti {E}_{<\ul{d}}$) le
sous-espace de $\ti {E}$ engendr\'e par les applications $\ti {e}_p$ avec $\ul{\deg}(p)\leq \ul{d}$
(resp. $\ul{\deg}(p)<\ul{d}$). Ce sous-espace est stable pour l'action par $\ti {\sigma }$. Notons
$\ti {\sigma }_{\leq \ul{d}}$ (resp. $\ti {\sigma }_{<\ul{d}}$) la repr\'esentation de $M$ dans cet
espace d\'eduite de $\ti{\sigma }$. Les repr\'esentations $\ti {\sigma }_{\leq \ul{d}}$ et $\ti {\sigma
}_{<\ul{d}}$ sont admissibles, puisque $\sigma $ l'est. 

Soit $p$ un mon\^ome avec $\ul{\deg}(p)=\ul {d}$. Tout \'el\'ement de $\ti
{E}_{\leq\ul{d}}/\ti{E}_{<\ul{d}}$ est l'image d'une application de la forme $\ti {e}_p$ avec $e\in E$.
Par ailleurs, on constate que l'application
$m\mapsto (p(H_M^G(mm_1))-p(H_M^G(m)))\sigma(m)(\sigma (m_1)e)$ est dans $\ti {E}_{<\ul {d}}$, puisque
$\ul{\deg }(p(.+H_M^G(m_1))-p(.))<\ul{d}$. On en d\'eduit que l'application de $E$ dans
$\ti{E}_{\leq\ul {d}}/\ti{E}_{<\ul {d}}$ qui associe \`a $e$ la classe de $\ti {e}_p$ d\'efinit un
isomorphisme entre $\sigma $ et la repr\'esentation quotient $\ti {\sigma }_{\leq\ul {d}}/\ti {\sigma
}_{<\ul{d}}$. Les repr\'esenta-tions $\ti{\sigma } _{\leq\ul {d}}$, $\ul {d}\in\hbox{\ens N}^l$, sont
donc admissibles de longeur finie et leurs sous-quotients irr\'eductibles sont isomorphes \`a $\sigma
$.  

Faisons correspondre \`a $e^{\vee }\in E^{\vee }$ l'application $\ti {e^{\vee
}}:\ti{E}_{\leq\ul{d}}\rightarrow\CC$, $\ti {e}\mapsto \langle\ti{e}(1), e^{\vee }\rangle $. C'est un
\'el\'ement du dual lisse de $\ti{E}_{\leq\ul{d}}$: il faut v\'erifier l'existence d'un sous-groupe
ouvert $H$ qui laisse $\ti {e^{\vee }}$ invariant. Comme $e^{\vee }$ est dans le dual lisse de $E$, il
existe un sous-groupe ouvert $H$ contenu dans $K$ qui laisse $e^{\vee }$ invariant. On a alors pour
tout $e\in E$, tout polyn\^ome $p$ sur $a_{M,\sCC}^G$, et tout
$h\in H$, 
$$\eqalign {\langle\ti{\sigma }(h)\ti {e}_p,\ti {e^{\vee }}\rangle&=\langle\ti {e}_p(h),e^{\vee
}\rangle=p(H_M^G(h))
\langle\sigma(h)e,e^{\vee }\rangle\cr  &=p(H_M^G(1))\langle e,e^{\vee }\rangle=\langle\ti{e}_p(1),
e^{\vee }\rangle\cr &=\langle\ti{e}_p,\ti{e^{\vee }}\rangle.\cr  }$$  Le lemme r\'esulte alors du fait
que (avec les notations de l'\'enonc\'e) l'application $m\mapsto p_1(H_M(m))\langle\sigma (m)e_1,$
$e_1^{\vee }\rangle
$ est \'egale au coefficient matriciel $m\mapsto\langle\ti{\sigma }_{\leq\ul{d}}(m)
\ti{e}_1,\ti{e_1^{\vee }}\rangle$ de $\ti{\sigma }_{\ul{d}}$ avec $\ul{d}=\ul{\deg}(p_1)$. \hfill{\fin
2} 

\null
\it Preuve: \rm (de la proposition)   Il suffit de consid\'erer le cas o\`u $\tau=v\otimes v^{\vee }$.  

Pour $g\in G$, on a 
$$\eqalign{
\langle i_P^G(\sigma\otimes\chi )(g)v,v^{\vee }\rangle&=\int _K\langle v_{\chi }(kg),v^{\vee
}(k)\rangle dk\cr  &=\int _K\delta _P^{1/2}(m_P(kg))\chi (m_P(kg))
\langle\sigma(m_P(kg))v(k_P(kg)),v^{\vee}(k)\rangle dk,\cr  }$$  o\`u $m_P(kg)$ et $k_P(kg)$
d\'esignent les composantes de $kg$ selon la d\'ecomposition $G=MUK$.

Les applications \'etant lisses, l'int\'egrale est en fait une somme finie. L'op\'erateur
diff\'erentiel commute donc avec l'int\'egrale et on obtient
$$\displaylines{
\qquad D_{\lambda }(\langle i_P^G(\sigma\otimes\chi _{\lambda })(g)v,v^{\vee }\rangle)_{\vert\lambda
=\lambda _0}=\hfill\cr
\hfill\int _K D_{\lambda }(\delta _P^{1/2}(m_P(kg))\chi_{\lambda }(m_P(kg))
\langle\sigma(m_P(kg))v(k_P(kg)),v^{\vee}(k)\rangle)_{\vert\lambda =\lambda_0}dk.\qquad\cr  }$$    Le
terme sous l'int\'egrale est une expression de la forme  $p_1(H_M^G(m_P(kg)) \delta _P^{1/2}(m_P(kg))$
$\chi_{\lambda }(m_P(kg)) \langle\sigma(m_P(kg))v(k_P(kg)),v^{\vee}(k)\rangle$ avec $p_1$ polyn\^ome
sur
$a_{M,\sCC}^G$ et $k\in K$.

D'apr\`es le lemme pr\'ec\'edent, il existe une repr\'esentation admissible lisse de longueur finie
$(\ti {\sigma },\ti {E})$ et des \'el\'ements
$\ti {v}(k_P(kg))$, $\ti {v^{\vee }}(k)$ de $\ti {E}$ et de $\ti {E}^{\vee }$, tels que cette
expression vaille $\delta _P^{1/2}(m_P(kg))\chi_{\lambda _0} (m_P(kg)) \langle\ti{\sigma
}(m_P(kg))\ti{v}(k_P(kg)),\ti{v^{\vee }} (k)\rangle $. En prenant pour $(\ti {\sigma },\ti {E})$ la
repr\'esentation explicit\'ee dans la preuve du lemme, on voit par ailleurs que 
$$(\ti{v}(k_P(kg)))(m)=p_1(H_M^G(m))\sigma (m) v(k_P(kg))\qquad\hbox{\rm pour $m\in M$}$$  et
$\langle\ti{e},\ti{v^{\vee }}(k)\rangle=\langle\ti{e}(1),v^{\vee }(k)\rangle$ pour $\ti{e}\in\ti{E}$.

On d\'efinit donc ainsi des applications $k\mapsto \ti{v}(k)$ et $k\mapsto
\ti{v^{\vee }}(k)$. Montrons qu'elles sont \`a valeurs dans $i_{P\cap K}^K\ti{E}$ et $i_{P\cap
K}^K\ti{E}^{\vee }$ respectivement:  elles sont toutes les deux invariantes \`a droite par un
sous-groupe ouvert de $K$, puisque $i_P^G\sigma $ est lisse. Par ailleurs, on trouve pour $k_M\in K\cap
M$,
$k_U\in K\cap U$ et $k\in K$
$$\eqalign {
\ti {v}(k_Mk_Uk)(m)&=p_1(H_M^G(m))\sigma (m) v(k_Mk_Uk)\cr &=p_1(H_M^G(m))\sigma (m)\sigma (k_M)v(k)\cr
&=p_1(H_M^G(mk_M))\sigma (mk_M)v(k)\cr &=(\ti{\sigma }(k_M)\ti v(k))(m)\cr  }$$  ainsi que pour tout
$\ti{e}\in\ti{E}$
$$\eqalign{
\langle\ti{e},\ti{v^{\vee}}(k_Mk_Uk)\rangle&=\langle\ti{e}(1),v^{\vee}(k_Mk_Uk)
\rangle\cr  &=\langle\ti{e}(1),\sigma^{\vee}(k_M)v^{\vee}(k)\rangle\cr &=\langle\sigma
(k_M)^{-1}\ti{e}(1), v^{\vee }(k)\rangle\cr  &=\langle\ti{\sigma
}(k_M)^{-1}\ti{e},\ti{v^{\vee}}(k)\rangle, \cr }$$  la derni\`ere \'egalit\'e r\'esultant du fait que
$(\ti{\sigma }(k_M^{-1}) \ti{e})(1)=\sigma(k_M^{-1})\ti{e}(1)$, puisque $H_M^G(k_M^{-1})=H_M^G(1)$.

On d\'eduit de ceci que $D_{\lambda }(\langle i_P^G(\sigma\otimes\chi _{\lambda })(g)v,v^{\vee
}\rangle)_{\vert \lambda=\lambda _0}$ est une expression de la forme
$$\eqalign{ &\int _K \delta _P^{1/2}(m_P(kg))\chi _{\lambda _0}(m_P(kg))
\langle\ti{\sigma }(m_P(kg))\ti{v}(k_P(kg)),\ti{v^{\vee }}(k)\rangle dk\cr =&\int _K\langle
(i_P^G(\ti{\sigma }\otimes\chi _{\lambda _0} )(g)\ti{v})(k), \ti{v^{\vee }}(k)\rangle dk\cr  =&\langle
i_P^G(\ti{\sigma }\otimes\chi _{\lambda _0})(g)\ti{v},\ti{v^{\vee }} \rangle,\cr  }$$  d'o\`u la
proposition, puisque l'on a remarqu\'e que $\ti{v}\otimes\ti{v^{\vee}} \in i_{P\cap K}^K\ti{E}\otimes
i_{P\cap K}^K\ti{E}^{\vee }$.

La derni\`ere assertion de la proposition concernant la restriction des repr\'esentation \`a $T_G$ se
voit par une v\'erification directe, en remarquant que $H_M^G(t)=0$ pour $t\in T_G$.
\hfill{\fin 2}

\null
\null
\null
\vfill
\eject
{\bf 6. Prolongement rationnel d'une identit\'e polynomiale}

\null
\null {\bf 6.1} Soit $P=MU$ un sous-groupe parabolique semi-standard et soit $\o $ l'orbite inertielle
d'une repr\'esentation irr\'eductible cuspidale de $M$. Soit $M'$ un sous-groupe de Levi semi-standard
de $G$ contenant $M$. Lorsque $\ti {\xi }^{M'}$ est une application rationnelle sur $\o $ \`a valeurs
dans un sous-espace de dimension finie de $i_{\ol{P}\cap K\cap M'}^{K\cap M'}E_{\so } \otimes i_{P\cap
K\cap M'}^{K\cap M'}$ $E_{\so }^{\vee }$, notons $\varphi _{\ti{\xi }^{M'}}$ l'application rationnelle
sur $\o $ \`a valeurs dans $i_{P\cap K\cap M'}^{K\cap M'}E_{\so }\otimes i_{P\cap K\cap M'}^{K\cap
M'}E_{\so }^{\vee }$ d\'efinie par
$$\varphi _{\ti {\xi }^{M'}}(\sigma )=\sum _{w\in W^{M'}(M,\so)} (J^{M'}_{P\cap M'\vert\overline
{wP}\cap M'}(\sigma )\lambda (w)\otimes J_{P\cap M'\vert wP\cap M'}^{M'}(\sigma ^{\vee })\lambda (w))\
\ti{\xi }^{M'}(w^{-1}\sigma ).$$   Si $\ti {\xi}^{M'}$ est polynomiale, $\varphi _{\ti {\xi }^{M'}}$
est la composante en $(P\cap M',\o )$ d'un \'el\'ement de l'espace de Paley-Wiener $PW(M')$ pour $M'$
(cf. {\bf 2.2}). Par abus de notation, cet \'el\'ement de l'espace de Paley-Wiener $PW(M')$ sera
\'egalement not\'e $\varphi _{\ti {\xi }^{M'}}$, ce qui permettra par exemple d'\'ecrire $\varphi _{\ti
{\xi }^{M'}} (P_1\cap M',\sigma _1)$, si $P_1$ est un deuxi\`eme parabolique semi-standard de
$G$ de Levi semi-standard $M$ et $\sigma _1\in \o$. 

On notera par ailleurs, pour $A\in\A(\S )$, $\Stab (A)=\{w\in W^{M_A}(M,\o ) \vert wA=A\}$, et, pour
$\sigma\in\o $, $\E_2(M',\sigma )$ l'ensemble des repr\'esentations de carr\'e int\'egrable de $M'$ de
support cuspidal \'egal \`a la classe de $W^{M'}$-conjugaison de $\sigma $. (Une condition
n\'ecessaire, pour que $\sigma $ soit dans le support cuspidal d'une repr\'esentation de carr\'e
int\'egrable de $M'$ est \'evidemment que la restriction de $\sigma $ \`a $T_{M'}$ agit par un
caract\`ere unitaire.) 

Pour simplifier, on \'ecrira parfois $E_P^G(J_{\ol {P}\vert P}^{-1},\xi,\sigma )$ \`a la place de
$E_P^G((J_{\ol{P}\vert P}^{-1}(\sigma )\otimes 1)\xi(\sigma ))$ etc., et $Res _A(E_P^G(J_{\ol{P}\vert
P}^{-1},\xi ,\sigma )(g))$ d\'esignera
$Res_A (E_P^G(J_{\ol{P}\vert P}^{-1},\xi ,.)(g))(\sigma )$ etc.

L'objet de ce paragraphe est la preuve de la proposition suivante: 

\null  {\bf Proposition:} \it Soit $P'=M'U'$ un sous-groupe parabolique semi-standard contenant $P$
avec $M'=M_A$ pour un certain \'el\'ement $A$ de $\A (\S )$. Fixons $\sigma\in A$ avec $\Re (\sigma
)=r(A)$. 

Supposons qu'il existe une donn\'ee de r\'esidu $Res _A$ en $A$, telle que, pour toute application
polynomiale $\xi ^{M'}$ sur $\o $ \`a valeurs dans un sous-espace de dimension finie de $i_{\ol{P}\cap
K\cap M'}^{K\cap M'}$ $E_{\so } \otimes i_{P\cap K\cap M'}^{K\cap M'}E_{\so }^{\vee }$ et tout $g\in
G$, on ait l'identit\'e   
$$\displaylinesno {
\qquad\sum _{w\in W^{M'}(M,\so )}\gamma(M'/M)\deg(\sigma )Res _{wA}(E_{P\cap M'}^{M'}(J_{\ol {P}\cap
M'\vert P\cap M'}^{M'\ -1},\xi ^{M'}, w(\sigma\otimes\chi'))(g))\hfill&(6.1.1)\cr
\hfill =\vert\Stab (A)\vert \sum_{\pi \in\sE_2(M',\sigma )} \deg(\pi )
E_{M'}^{M'}(\varphi_{\xi^{M'}}(M',\pi\otimes\chi '))(g).\qquad\cr }$$  de fonctions rationnelles sur
$\X (M')$.

Soit $\ti {\xi }^{M'}$ une application rationnelle sur $\o $ \`a valeurs dans un sous-espace de
dimension finie de $i_{\ol{P}\cap K\cap M'}^{K\cap M'}E_{\so }\otimes i_{P\cap K\cap M'}^{K\cap
M'}E_{\so }^{\vee }$ telle que  $\varphi _{\ti {\xi}^{M'}}$ soit le produit de la composante en $(P\cap
M',\o )$ d'un \'el\'ement de l'espace de Paley-Wiener $PW(M')$ de $M'$ par une fonction rationnelle
$\ti {\mu }$. Supposons qu'il existe un ensemble fini $\ti{\S }$ d'hyperplans radiciels de $\o $ de la
forme $\o_{\sigma ',\alpha }$ avec $\alpha\in\Sigma (P)-\Sigma ^{M'}(P\cap M')$ tel que $\ti {\xi
}^{M'}$ et $\ti {\mu }$ soient r\'eguli\`eres sur $\o -\bigcup_{S\in\ti {\sS }}S$. Alors l'identit\'e
ci-dessus vaut pour $\ti {\xi }^{M'}$. 

\null  Preuve: \rm  Remarquons que, pour $\pi\in\E _2(M',\sigma )$, $\varphi _{\ti {\xi }^{M'}}
(M',\pi\otimes\chi')$ est bien d\'efinie pour tout $\chi'\in\X (M')$ avec $\sigma\otimes\chi '$
r\'egulier pour $\ti {\mu}$. Le terme \`a  droite de l'identit\'e (6.1.1) est donc une fonction
rationnelle sur $\X (M')$. La fonction \`a gauche est rationnelle par d\'efinition des donn\'ees de
r\'esidu. Par prolongement analytique, il suffit donc de montrer qu'il existe $r_{M'}\in a_{M'}^*$, tel
que l'identit\'e (6.1.1) vaut en tout $\chi '\in\X (M')$ de partie r\'eelle $>_{P'}r_{M'}$. Remarquons
que $\Re (\sigma )_{M'}=0$.

L'ensemble des sous-espaces affines radiciels $A$ de $\o $ dans $\A (\S )$ avec $M_A=M'$ \'etant fini,
on peut choisir $r^{M'}$ dans la chambre de Weyl de $P\cap M'$ dans $a_M^{M'*}$, tel que tout point
$\sigma '$ d'un tel sous-espace v\'erifie $-r^{M'}<_{P\cap M'}\Re (\sigma ')^{M'}<_{P\cap M'}r^{M'}$.
On va alors prendre pour $r_{M'}$ un \'el\'ement de la chambre de Weyl de $P'$ dans $a_{M'}^*$ tel que,
pour tout $H_{\alpha ,c}\in\H(\ti{\S })$ avec $\alpha\in\Sigma (P)-\Sigma ^{M'}(P\cap M')$, et tout
$\lambda ^{M'}\in a_M^{M'*}$ avec $-r^{M'}\leq _{P\cap M'}\lambda ^{M'}\leq _{P\cap M'}r^{M'}$, on ait
$\langle\alpha^{\vee }, r_{M'}\rangle >\Re(c)-\langle\alpha^{\vee },\lambda ^{M'}\rangle $. (Ceci est
possible, puisque les formes lin\'eaires $\langle\alpha ^{\vee },.\rangle$, $\alpha\in\Sigma
(P)-\Sigma ^{M'}(P\cap M')$, sont non nulles sur $a_{M'}^*$, alors que $\lambda ^{M'}$ parcourt un
compact.)   Notons $U$ le domaine des \'el\'ements de $\o $ dont la partie r\'eelle
$\lambda $ v\'erifie $\lambda _{M'}>_{P'}r_{M'}$ et $-r^{M'}<_{P\cap M'}\lambda ^{M'}<_{P\cap
M'}r^{M'}$. Les fonctions $\ti {\xi }^{M'}$ et $\ti {\mu }$ sont r\'eguli\`eres sur $U$: en effet, par
choix de $r_{M'}$, on a, pour
$H_{\alpha ,c}\in\H(\ti{\S})$ avec $\alpha\in\Sigma (P)$, et $\sigma'\in U$,
$\langle\alpha ^{\vee },\Re (\sigma )_{M'}\rangle >\langle\alpha ^{\vee },r_{M'}\rangle >\Re
(c)-\langle\alpha ^{\vee },\Re (\sigma ')^{M'}\rangle
$, et donc $\langle\alpha ^{\vee },\Re (\sigma )\rangle\ne\Re (c)$.  

En particulier, la fonction $\chi\mapsto \ti{\xi }^{M'}(\sigma\otimes\chi )$ est pour tout
$r>_{P'}r_{M'}$ analytique sur le polydisque $U_r'$ de
$\X (M)$ form\'e des caract\`eres $\chi $ v\'erifiant $-r^{M'}<_{P\cap M'} \Re(\sigma\otimes\chi
)<_{P\cap M'}r^{M'}$ et $r_{M'}<_{P'}\Re (\chi )<_{P'}r$. Par la th\'eorie des fonctions analytiques en
plusieurs variables, $\ti {\xi }^{M'}(\sigma\otimes .)$ peut s'\'ecrire comme limite d'une suite
d'applications polynomiales $(\xi _n)$ convergeant uniform\'ement sur tout compact de $U_r'$ et prenant
ses valeurs dans un m\^eme sous-espace de dimension finie de
$i_{\ol{P}\cap K\cap M'}^{K\cap M'}E_{\so }\otimes i_{P\cap K\cap M'}^{K\cap M'}E_{\so }^{\vee }$ que
$\ti {\xi }^{M'}$. Pour obtenir des applications polynomiales sur $\o $, on les rend invariantes par
$\Stab (\o )$ (cf. {\bf 1.3}). Ce proc\'ed\'e ne change rien \`a la convergence uniforme sur tout
compact de $U_r'$ vers $\ti {\xi }^{M'}$, puisque $\ti {\xi }^{M'}$ est invariante. 

Notons $U_r$ l'ensemble des $\sigma\otimes\chi $ avec $\chi\in U_r'$. Cet ensemble est \'egal au
sous-ensemble de $U$ form\'e des $\sigma '$ avec
$\Re (\sigma ')_{P'}<_{P'}r$. On a donc montr\'e l'existence d'une suite $(\xi _n^{M'})$
d'applications polynomiales sur $\o $ qui converge uniform\'ement sur tout compact de $U_r$ vers $\ti
{\xi }^{M'}$. Remarquons que, comme $\sigma $ a \'et\'e choisie dans $A$, on a, par choix de $r^{M'}$
et gr\^ace \`a l'invariance de $\A (\S )$ par les \'el\'ements de $W^{M'}(M,\o )$, $w(\sigma\otimes\chi
')\in U_r$ pour tout $w\in W^{M'}(M,\o )$ et tout $\chi'\in\X (M')$, $r_{M'}<_{P'}\Re (\chi ')<_{P'}r$.

Suite aux hypoth\`eses de la proposition, l'identit\'e (6.1.1) est v\'erifi\'ee pour $\xi _n^{M'}$ en
tout $\chi'\in\X (M')$, et donc en particulier si $r_{M'}<_{P'}\Re (\chi ')<_{P'}r$. Il reste \`a voir,
si la suite des valeurs en $\sigma\otimes\chi '$ converge pour $n\rightarrow\infty$ respectivement vers 
$$\sum_{w\in W^{M'}(M,\so )}\gamma (M'/M)\deg(\sigma )Res _{wA}(E_{P\cap M'}^{M'} (J_{\ol{P} \cap
M'\vert P\cap M'}^{M'\ -1}, \ti {\xi }^{M'},w(\sigma\otimes\chi '))(g))$$
$$\vert\Stab (A)\vert\sum _{\pi\in\sE_2(M',\sigma )}\deg(\pi ) E_{M'}^{M'}(\varphi _{\ti{\xi
}^{M'}}(M',\pi\otimes\chi '))(g).\leqno {\hbox{et}}$$   L'identit\'e (6.1.1) en $\chi'$ en r\'esulte.

\'Etudions les deux cas s\'epar\'ement. Comme le calcul de la donn\'ee de r\'esidu $Res _A$ ne fait
intervenir que des combinaisons lin\'eares r\'eelles de d\'eriv\'ees en direction d'\'el\'ements de
$a_{M,\sCC}^{M'*}$ et qu'une suite de d\'eriv\'ees de $\xi _n^{M'}$ converge par le th\'eor\`eme de la
convergence uniforme des fonctions holomorphes  uniform\'ement sur tout compact de $U_r$ vers la
d\'eriv\'ee de $\ti {\xi }^{M'}$, on trouve pour tout
$w\in W^{M'}(M,\o )$, ayant $w(\sigma\otimes\chi')\in U_r$, 
$$\displaylines{
\qquad\lim _{n\rightarrow\infty} Res_{wA}(E_{P\cap M'}^{M'}(J_{\ol{P}\cap M'\vert P\cap M'}^{M'\
-1},\xi _n^{M'}, w(\sigma\otimes\chi '))(g))=\hfill\cr
\hfill Res _{wA}(E_{P\cap M'}^{M'}(J_{\ol{P}\cap M'\vert P\cap M'}^{M'\ -1},\ti {\xi
}^{M'},w(\sigma\otimes\chi '))(g)).\qquad\cr }$$ 

Concernant le deuxi\`eme terme, il suffit de prouver que $\varphi _{\xi _n^{M'}}(M',\sigma\otimes\chi
')$ converge vers $\varphi _{\ti{\xi }^{M'}} (M',\sigma\otimes \chi')$. Notant $\mu ^+$ un
d\'enominateur pour la fonction rationnelle $\mu ^{M'}$ qui est invariant par $W^{M'}(M,\o )$, on a une
identit\'e  
$$\displaylinesno {
\qquad\varphi _{\mu^+\xi _n^{M'}}(\sigma ')=\sum _{w\in W^{M'}(M,\so )}\mu ^+(\sigma'
)\hfill\qquad&(6.1.2)\cr
\hfill(J_{P\cap M'\vert\ol{wP}\cap M'}^{M'}(\sigma ')\lambda (w)\otimes J_{P\cap M'\vert wP\cap
M'}^{M'}(\sigma '^{\vee }) \lambda (w))\xi _n^{M'}(w^{-1}\sigma ') \qquad\cr   }$$   de fonctions
rationnelles sur $\o $. Choisissons un voisinage relativement compact et connexe $U'$ de
$\sigma\otimes\chi'$ dans $U_r$ tel que $wU'\subseteq U_r$ pour tout $w\in W^{M'}(M,\o )$. Comme
l'op\'erateur $\sigma'\mapsto\mu^+(\sigma ') (J_{P\cap M'\vert\ol{wP}\cap M'}^{M'}(\sigma ')\lambda
(w)\otimes  J_{P\cap M'\vert wP\cap M'}^{M'}(\sigma '^{\vee })\lambda (w))$ est r\'egulier sur $\o $
(cf. {\bf 4.4}) et que $\mu ^+$ est invariant par $W^{M'}(M,\o)$, (6.1.2) converge uniform\'ement sur
tout compact de $U'$ vers $\varphi _{\mu ^+ \ti {\xi }^{M'}}(\sigma ')$. Par suite, 
$\lim _{n\rightarrow\infty }$ $\varphi _{\xi _n^{M'}}(\sigma ')=
\lim _{n\rightarrow\infty } {\varphi _{\mu ^+\xi _n^{M'}}(\sigma ')\over\mu ^+(\sigma ')}= {\varphi
_{\mu ^+\ti {\xi }^{M'}}(\sigma ')\over\mu ^+(\sigma ')}= 
\varphi _{\ti {\xi }^{M'}}(\sigma ')$  en tout point $\sigma '$ de $U'$ avec $\mu ^+(\sigma ')\ne 0$.
Comme l'op\'erateur $\varphi _{\ti {\xi }^{M'}}$  est r\'egulier en $\sigma\otimes\chi'$ et que la
limite d\'efinit un op\'erateur m\'eromorphe sur $U'$, on en d\'eduit l'\'egalit\'e en
$\sigma\otimes\chi '$ par prolongement analytique.\hfill{\fin 2}  

\null
\null {\bf 6.2} La proposition {\bf 6.1} sera utilis\'ee dans la suite pour les applications
rationnelles d\'ecrites dans le lemme suivant:

\null {\bf Lemme:} \it Soit $\xi $ une application polynomiale sur $\o $ \`a valeurs dans un
sous-espace de dimension finie de $i_{\ol{P}\cap K}^KE_{\so }\otimes i_{P\cap K}^KE_{\so }^{\vee }$.
Soit $P'=M'U'$ un sous-groupe parabolique $(P,\S)$-standard, $P_1=(P\cap M')U'$ et $\ti {P_1}=(P\cap
M')\ol{U'}$. Identifions $i_{\ol{\ti{P}}_1}^GE_{\so }\otimes i_{P_1}^GE_{\so }^{\vee }$ avec
$i_{P'}^Gi_{\ol{P}\cap M'}^{M'}E_{\so }\otimes i_{P'}^Gi_{P\cap M'}^{M'}E_{\so }^{\vee }$. Soient
$g',g\in G$. L'application rationnelle $\ti{\xi}_{g'g,g'}^{M'}:\o\rightarrow i_{\ol{P}_1\cap K\cap
M'}^{K\cap M'}E_{\so }\otimes i_{P_1\cap K\cap M'}^{K\cap M'}E_{\so }^{\vee }$,
$$\displaylines {\qquad\sigma\mapsto \mu _{P_1\vert\ti {P_1}}(\sigma )[\sum _{w\in W_{M'}^+(M,\so )}
\mu_{w^{-1}P\cap M'\vert P\cap M'}^{M'}(\sigma )\hfill\cr
\hfill(J_{\ol{\ti{P}}_1\vert\ol{w^{-1}P}}(\sigma )\lambda (w^{-1})\otimes J_{P_1\vert w^{-1}P}(\sigma
^{\vee })\lambda (w^{-1}))\xi (w\sigma )](g'g,g')\qquad\cr }$$  est r\'eguli\`ere en dehors d'un
ensemble fini $\ti {\S }$ d'hyperplans radiciels de $\o $ de la forme $\o _{\sigma ,\alpha }$ avec
$\alpha\in\Sigma (P_1)-\Sigma ^{M'}(P\cap M')$. Par ailleurs, $\varphi ^{M'}_{\ti {\xi
}_{g'g,g'}^{M'}}(\sigma )=\mu _{P_1\vert\ti {P}_1}(\sigma )(\varphi _{\xi } (P_1,\sigma )(g'g,g'))$. 

\null
\it Remarque: \rm L'application polynomiale $\sigma \mapsto \varphi _{\xi }(P_1,\sigma )(g'g,g')$ du
lemme est bien la composante en $P_1\cap M'$ d'un \'el\'ement de l'espace de Paley-Wiener de $M'$,
puisque les hypoth\`eses de {\bf 2.3} valent relatives \`a cette application (ce qui est une
v\'erification directe en utilisant la fonctorialit\'e des op\'erateurs d'entrelacement vis-\`a-vis de
l'induction parabolique).

\null
\null 
\it Preuve: \rm La fonction $\mu _{w^{-1}P\cap M'\vert P\cap M'}^{M'}$ est polynomiale sur $\o $, si
$w\in W_{M'}^+(M,\o )$: pour que cette fonction ait un p\^ole, il faut par la formule du produit pour
la fonction $\mu $ et {\bf 4.3} qu'il existe une racine $\alpha\in\Sigma _{\so }^{M'}\cap\Sigma (P)$
qui est n\'egative pour $w^{-1}P\cap M'$. Or, comme $w$ est le produit d'un \'el\'ement de $W$ qui
envoie $\Delta _{\so }$ dans lui-m\^eme par un \'el\'ement de $W_{\so }$ de longueur minimale dans sa
classe \`a droite modulo $W_{\so }^{M'}$, $\alpha\in \Sigma _{\so }^{M'} \cap\Sigma (P)$ implique
$w\alpha\in\Sigma (P)$, i.e. $\alpha $ est n\'ecessairement positif pour $w^{-1}P\cap M'$.

Par la formule de d\'ecomposition pour les op\'erateurs d'entrelacement, on sait de m\^eme qu'un
hyperplan de la forme $\o _{\sigma ,\alpha }$ avec $\alpha\in\Sigma ^{M'}(P\cap M')$ est singulier pour
$J_{(M'\cap\ol{P})U' \vert\ol{w^{-1}P}}$ (resp. $J_{P_1\vert w^{-1}P}$), seulement s'il existe
$\alpha\in\Sigma _{\so }^{M'}\cap \Sigma (P)$ qui est n\'egative pour
$w^{-1}P\cap M'$. Or, on vient de voir qu'un tel $\alpha $ n'existe pas. Il est clair que les
hyperplans de cette forme ne sont pas non plus singuliers pour $\mu _{P_1\vert\ti{P}_1}$. Ceci prouve
la premi\`ere assertion du lemme.

D\'eterminons maintenant $\varphi _{\ti {\xi }_{g'g,g'}^{M'}}^{M'}$ en un point g\'en\'erique $\sigma
$. Comme $\mu _{P_1\vert\ti {P_1}}$ est invariante par
$W^{M'}(M,\o)$ en tant que quotient des fonctions $W^{M'}(M,\o)$-invariantes $\mu $ et $\mu ^{M'}$, le
produit de $\mu _{P_1\vert\ti {P}_1}(\sigma )^{-1}$ avec
$$\sum _{w'\in W^{M'}(M,\so )}(J_{P_1\cap M'\vert\ol{w'P_1}\cap M'}^{M'}(\sigma )\lambda (w')\otimes
J_{P_1\cap M'\vert w'P_1\cap M'}^{M'}(\sigma )\lambda (w')) \ti{\xi }_{g'g,g'}^{M'}(w'^{-1}\sigma )$$ 
est l'\'evaluation en $(g'g,g')$ de
$$\displaylines {
\ \sum _{w'\in W^{M'}(M,\so )}\sum _{w\in W_{M'}^+(M,\so )} \mu^{M'}_{w^{-1}P\cap M'\vert P\cap
M'}(w'^{-1}\sigma )\hfill\cr
\qquad (J_{P_1\vert (\ol{w'P}\cap M')U'}(\sigma )\lambda (w')J_{(\ol{P}\cap M')
U'\vert\ol{w^{-1}P}}(w^{-1}\sigma )\lambda (w^{-1})\otimes\hfill\cr
\hfill J_{P_1\vert w'P_1}(\sigma ^{\vee })\lambda (w')J_{P_1\vert w^{-1}P} (w'^{-1}\sigma ^{\vee
})\lambda (w^{-1}))\xi (ww'^{-1}\sigma )\qquad\cr =\sum _{w'\in W^{M'}(M,\so )}\sum _{w\in
W_{M'}^+(M,\so )} \mu ^{M'}_{w^{-1}P\cap M'\vert P\cap M'}(w'^{-1}\sigma )\hfill\cr
\qquad (J_{P_1\vert(\ol{w'P}\cap M')U'}(\sigma )J_{(\ol{w'P}\cap M')U'\vert \ol{w'w^{-1}P}}(\sigma
)\lambda (w'w^{-1})\otimes \hfill\cr
\hfill J_{P_1\vert w'P_1}(\sigma ^{\vee })J_{w'P_1\vert w'w^{-1}P}(\sigma ^{\vee })\lambda
(w'w^{-1}))\xi (ww'^{-1}\sigma )\qquad\cr  =\sum _{w'\in W^{M'}(M,\so )}\sum _{w\in W_{M'}^+(M,\so )}
\mu ^{M'}_{w^{-1}P\cap M'\vert P\cap M'}(w'^{-1}\sigma )j_{w'w^{-1}P\cap M'\vert w'P\cap
M'}^{M'}(\sigma )\hfill\cr
\hfill (J_{P_1\vert\ol{w'w^{-1}P}}(\sigma )\lambda (w'w^{-1})\otimes J_{P_1\vert w'w^{-1}P}(\sigma
^{\vee })\lambda (w'w^{-1}))\xi (ww'^{-1}\sigma ),\qquad\cr }$$  o\`u on a utilis\'e la formule du
produit pour les op\'erateurs d'entrelacement. Comme
$$\eqalign { j_{w'w^{-1}P\cap M'\vert w'P\cap M'}^{M'}(\sigma )\ id &=\lambda
(w'^{-1})J^{M'}_{w'w^{-1}P\cap M'\vert w'P\cap M'}(\sigma ) J^{M'}_{wP\cap M'\vert w'w^{-1}P\cap
M'}(\sigma )\lambda (w')\cr &=J^{M'}_{w^{-1}P\cap M'\vert P\cap M'}(w'^{-1}\sigma )J^{M'}_{P\cap
M'\vert w^{-1}P\cap M'}(w'^{-1}\sigma )\cr  &=j ^{M'}_{w^{-1}P\cap M'\vert P\cap M'}(w'^{-1}\sigma )\
id,\cr  }$$  on a $\mu ^{M'}_{w^{-1}P\cap M'\vert P\cap M'}(w'^{-1}\sigma )j_{w'w^{-1}P\cap M'\vert
w'P\cap M'}^{M'}(\sigma )=1$.  Remarquons que $W(M,\o )=W^{M'}(M,\o )(W_{M'}^+(M,\o ))^{-1}$ par {\bf
4.5}. L'expression que l'on a obtenu est donc \'egale \`a $(J_{P_1\vert P}(\sigma )\otimes J_{P\vert
P_1}(\sigma ^{\vee })^{-1})\varphi _{\xi }(\sigma )$, ce qui vaut $\varphi _{\xi }(P_1,\sigma )$,
d'o\`u le lemme. \hfill {\fin 2}

\null
\null
\null
\vfill
\eject
{\bf 7. L'\'el\'ement $\varphi _{P,\so }^H$ de l'espace de Paley-Wiener}

\null
\null  {\bf 7.1} Fixons un sous-groupe parabolique semi-standard $P=MU$ de $G$ et l'orbite inertielle
$\o $ d'une repr\'esentation irr\'eductible cuspidale de $M$. Fixons en outre un sous-groupe ouvert
compact distingu\'e $H$ de $K$ admettant une d\'ecomposition d'Iwahori par rapport \`a tout couple
parabolique semi-standard de $G$ et tel que $E_{\so }^{H\cap M}\ne 0$. On a d\'efini dans [H1] B.4.
l'application polynomiale $\varphi _{P,\so }^H$ suivante sur $\o $ \`a valeurs dans $i_P^GE_{\so
}\otimes i_P^GE_{\so }^{\vee }$:

D\'esignons pour une repr\'esentation lisse $(\pi ,V)$ de $G$ par $V_P$ son module de Jacquet qui est
une repr\'esentation lisse de $M$. L'espace $V_P$ est somme directe de deux sous-espace
$M$-\'equivariante $V_P(\o )$ et $V_P(hors \o)$ tels que tout sous-quotient irr\'eductible de la
repr\'esentation de $M$ dans $V_P(\o )$ soit dans $\o $ et qu'aucun sous-quotient de la
repr\'esentation de $M$ dans $V_P(hors \o)$ ne le soit.

Pour $H$ un sous-groupe ouvert compact de $G$ admettons une d\'ecomposition d'Iwahori par rapport \`a
$P=MU$, le foncteur de Jacquet $V^H\rightarrow V_P^{H\cap M}$ d\'efini relatif aux invariants par $H$
et par $H\cap M$ admet une section canonique sur un sous-espace not\'e $S_P^H(V)$ de $V^H$.

La valeur de $\varphi _{P,\so }^H$ en un point $\sigma $ de $\o $ est alors \'egale \`a la projection
de $i_{P\cap K}^KE_{\sigma }$ sur $S_P^H(i_{P\cap K}^K E_{\sigma })(\o )$ de noyau \'egal \`a
l'intersection des \'el\'ements de $S_{\ol {P}}^H(i_{P\cap K}^KE_{\sigma }^{\vee })(\o ^{\vee })$.

On a montr\'e qu'il existe un \'el\'ement $\varphi _{\so }$ de l'espace de Paley-Wiener de $PW(G)$ tel
que $(\varphi _{\so })_{P,\so }=\varphi _{P,\so }^H$ et que $(\varphi _{\so })_{P',\so '}=0$, si $\o'$
n'est pas conjugu\'e \`a $\o $.

Choisissons des \'el\'ements $v_i\in (i_{\ol {P}\cap K}^KE_{\so })^H$ (resp. $v_i^{\vee }\in (i_{P\cap
K}^KE_{\so }^{\vee })^H$), $i\in I$, \`a support dans $(\ol {P}\cap K)H$ (resp. $(P\cap K)H)$, tels que
$e_i:=v_i(1)$ (resp.
$e_i^{\vee }:=v_i^{\vee }(1)$) forment des bases de $E_{\so }^{M\cap H}$ (resp. $(E_{\so ^{\vee
}}^{\vee })^{M\cap H}$) qui sont duales. Posons $c=\mes _{U_P}(H\cap U_P) \mes _{\ol{U}_P}(H\cap\ol
{U}_P)$.  On a montr\'e dans [H1] B.5 qu'une solution polynomiale de la relation (\#) de l'introduction
relative \`a $\varphi _{P,\so }^H$ est donn\'ee par
$$\xi _{\varphi _{P,\sso }^H}:=c^{-1}\sum _{i\in I} v_i\otimes v_i^{\vee }.$$

Notons $f_{\varphi _{\sso }}$ l'\'el\'ement de $C_c^{\infty }(G)$ qui correspond \`a $\varphi _{\so }$
(cf. {\bf 2.3}). 

\null
\null   {\bf Lemme:} \it Soit $\pi $ une repr\'esentation de carr\'e int\'egrable de $G$ dont le
support cuspidal est contenu dans la classe de
$W$-conjugaison de $\o $. Alors $\tr (\pi (f_{\varphi _{\sso }}))$ est un entier strictement positif. 

\null   Preuve: \rm  Notons $V_{\pi }$ l'espace de $\pi $. D'apr\`es [H1] B.3.4.1, on a $(V_{\pi
})_P(\o )\ne 0$. Par r\'eciprocit\'e de Frob\'enius, il existe donc $\sigma\in \o $, tel que $\pi
\hookrightarrow i_P^G\sigma $. Il en r\'esulte en outre que $S_P^H(i_P^G\sigma )(\o )\cap V_{\pi
}\supseteq S_P^H(V_{\pi })(\o )\ne 0$ (dans les notations de [H1] B.2) et que $\pi(f_{\varphi _{\sso
}})=\varphi _{P,\so }^H (\sigma )\mid_{V_{\pi }}$. L'endomorphisme $\pi (f_{\varphi })$ est donc comme
restriction d'une projection une projection ou bien identiquement 0. Comme la trace d'une projection
est toujours un entier strictement positif, il reste \`a montrer que $\pi (f_{\varphi })\ne 0$. Or,
d'apr\`es ce qui pr\'ec\'edait, il existe $v\in S_P^H(i_P^G\sigma )(\o )\cap V_{\pi }$, $v\ne 0$, et on
a $\pi (f_{\varphi })(v)=\varphi _{P,\so }^H(\sigma )(v)=v\ne 0$, d'o\`u le lemme. 

\hfill {\fin 2} 

\null
\null  {\bf 7.2 Lemme:} \it L'application rationnelle
$$\o\rightarrow \CC,\qquad \sigma\mapsto E_{P,\sigma }^G((J_{P\vert\ol{P}}(\sigma )\otimes 1)\xi
_{\varphi _{P,\sso}^H}(\sigma ))(1)$$  est constante et $>0$.

\null  Preuve: \rm Remarquons d'abord que, pour $\sigma $ r\'egulier,
$$\eqalign { E_{P,\sigma }^G((J_{P\vert\ol{P}}(\sigma )\otimes 1)\xi _{\varphi _{P,\sso}^H} (\sigma
))(1)&=c^{-1}\sum _{i\in I}E_{P,\sigma }^G(J_{P\vert\ol{P}}(\sigma )v_i \otimes v_i^{\vee })(1)\cr 
&=c^{-1}\sum _{i\in I} \langle J_{P\vert\ol{P}}(\sigma )v_i,v_i^{\vee }\rangle\cr  &=c^{-1}\sum _{i\in
I} \int_K \langle (J_{P\vert\ol{P}}(\sigma )v_i)(k), v_i^{\vee }(k)\rangle dk.\cr   }$$   Soit $p\in
P\cap K$ et $h\in H$. Notons $m_P(p)u_P(p)$ la d\'ecomposition de $p$ selon la d\'ecomposition $P=MU$.
On a 
$$\eqalign{ (J_{P\vert\ol{P}}(\sigma )v_i)(ph)&=((i_P^G\sigma )(h)J_{P\vert\ol{P}}(\sigma
)v_i)(p)=(J_{P\vert\ol{P}}(\sigma )(i_{\ol {P}}^G \sigma )(h)v_i)(p)\cr  &=(J_{P\vert\ol{P}}(\sigma
)v_i)(p)=\sigma (m_P(p))(J_{P\vert\ol{P}}(\sigma )v_i)(1),\cr   }$$   avec, d'apr\`es le lemme B.5.2.1
de [H1],
$$(J_{P\vert\ol{P}}(\sigma )v_i)(1)=\mes _{U_P}(H\cap U_P)e_i,$$    alors que $v_i^{\vee }(ph)=\sigma
^{\vee }(m_P(p)) v_i^{\vee }(1)=\sigma ^{\vee }(m_P(p))e_i^{\vee }$.

Comme $\supp(v_i^{\vee })\subseteq (P\cap K)H$ et que $\langle e_i,e_i^{\vee }\rangle =1$, on en
d\'eduit que
$$\int_K \langle (J_{P\vert\ol{P}}(\sigma )v_i)(k), v_i^{\vee }(k)\rangle dk =\mes((P\cap K)H)\mes
_{U_P}(H\cap U_P).$$   Par suite, $E_{P,\sigma }^G((J_{P\vert\ol{P}}(\sigma )\otimes 1)\xi _{\varphi
_{P,\sso}^H}(\sigma ))(1)$ est constant et $>0$ pour $\sigma $ dans un ouvert de Zariski de $\o $,
d'o\`u le lemme par prolongement analytique. \hfill{\fin 2} 

\null
\null 
\null
\vfill
\eject
{\bf 8. La d\'ecomposition spectrale}

\null Dans tout ce paragraphe, on fixe un sous-groupe parabolique semi-standard $P=MU$ de $G$ et
l'orbite inertielle $\o $ d'une repr\'esentation irr\'eductible cuspidale de $M$. Par ailleurs, le
symbole $\S $ d\'esignera dans tout ce paragraphe l'ensemble fini $\S _{\mu }$ d'hyperplans affines de
$\o $ associ\'e \`a la fonction $\mu $ de Harish-Chandra d\'efinie relative \`a $\o $ (cf. {\bf 4.3}).

\null  {\bf 8.1} Appelons un sous-espace affine radiciel $A$ de $\o $ \it r\'esiduel pour $\mu $, \rm
s'il existe une fonction polynomiale $p$ sur $\o $ et une donn\'ee de r\'esidu $Res _A$ en $A$, tel que
$Res _A(p\mu)\ne 0$. Notons $\A _{\mu }(\S )$ le sous-ensemble de $\A (\S )$ form\'e des sous-espaces
affines qui sont r\'esiduels pour $\mu $. 

Notons $\rg _{ss}(H)$ le rang semi-simple d'un groupe r\'eductif $H$. Si $M'$ est un sous-groupe de
Levi de $G$, on appelle $\rg _{ss}(G)-\rg _{ss}(M')$ le \it rang parabolique de $M'$ (relatif \`a $G$).
\rm

Le lemme suivant est une reformulation de r\'esultats combinatoires de Heckman-Opdam et Opdam:

\null {\bf Lemme:} \it Soit $\Sigma '$ un syst\`eme de racines dans $a_M^{G*}$ de dimension $m$ de
groupe de Weyl $W'$, soit $q_1$ un r\'eel $>1$ et soit $\underline {k}=(k_{\alpha })_{{\alpha
}\in\Sigma '}$ une famille de nombres $>0$ tels que $k_{w\alpha }=k_{\alpha }$ pour tout $w\in W'$.
Posons pour $\lambda\in a_{M,\sCC }^{G*}$
$$\mu (\lambda ,\ul {k})=\prod _{\alpha\in\Sigma '} {1-q_1^{\langle\alpha ^{\vee },\lambda\rangle
}\over  1-q_1^{k_{\alpha }+\langle\alpha^{\vee },\lambda \rangle}}$$

Pour qu'il existe une suite strictement croissante de sous-groupes de Levi $M\subset M_1\subset
\cdots\subset M_{m-1}\subset M_m=G$ et une fonction $p$ m\'eromorphe sur $a_{M,\sCC }^{G*}$ et holomorphe 
en un voisinage de $\lambda $ telle que les r\'esidus successifs de $p\mu (.,\ul {k})$ par rapport aux
sous-espaces $\lambda +a_{M_i,\sCC }^{G*}$ soient non nuls, il faut et il suffit que $\lambda $ soit un
p\^ole d'ordre $m$ de $\mu (.,\ul {k})$ (l'ordre d'un p\^ole \'etant d\'efini comme dans ({\bf 3.4})).
Ces p\^oles de $\mu (.,\ul {k})$ sont d'ordre maximal. 

De plus, si $\lambda $ est un p\^ole d'ordre maximal $m$ de $\mu (.,\ul {k})$, alors les r\'esidus
successifs de $\mu (.,\ul {k})$ par rapport aux sous-espaces $\lambda +a_{M_i,\sCC }^{G*}$ sont des
fonctions non nulles qui ont des p\^oles simples par rapport aux sous-espaces $\lambda +a_{M_{i+1},\sCC
}^{G*}$.

\null Preuve: \rm Il r\'esulte de [HO1] th\'eor\`eme 3.9 et identit\'e (3.17) que la condition est
suffisante, pour que $\lambda $ soit un p\^ole d'ordre $m$, et que les p\^oles des r\'esidus successifs
par rapport aux $a_{M_i,\sCC }^{G*}$ sont alors simples. La r\'eciproque et la maximalit\'e de l'ordre
de ces p\^oles r\'esultent de [O2] corollaire 7.10 et remarque 7.11. \hfill {\fin 2} 

\null
\null  {\bf Proposition:} \it Pour qu'un sous-espace affine radiciel $A$ de $\o $ soit dans $\A _{\mu
}(\S)$, il faut et il suffit qu'un (ou tout) \'el\'ement de $A_{\sreg }$ soit un p\^ole de $\mu^{M_A}$
d'ordre \'egal au rang parabolique de $M$ d\'efini relatif \`a $M_A$. Ces p\^oles sont d'ordre maximal
pour $\mu ^{M_A}$.

\null
\it Preuve: \rm Soit $\sigma\in A_{\sreg }$. Consid\'erons d'abord le cas $M_A=G$. On a alors
$A=A_{\sreg }$. En \'echangeant $P$ le cas \'ech\'eant, on peut supposer $\Re (\sigma )\geq _P0$. Si
$\o_{\sigma, \alpha}$ est un hyperplan singulier pour $\mu $, alors $\mu^{M_{\alpha }}(\Im(\sigma ))=0$
par {\bf 4.1}. Si $\mu ^{M_{\alpha }}$ s'annule sur $\o_{\sigma, \alpha }$, alors la projection $\Re
(\sigma )^{M_{\alpha }}$ de $\Re (\sigma )$ sur
$a_M^{M_{\alpha }*}$ est nulle (cf. {\bf 4.1}). Notons $\Sigma _{sp}^+(\sigma )$ le sous-ensemble de
l'ensemble $\Sigma _{sp}(\sigma )$ (d\'efini en {\bf 4.1}), form\'e des racines qui sont positives pour
$P$. \'Ecrivons $\mu=\mu _{nsp,\sigma }\ \mu_{sp, \sigma }$ avec les notations de {\bf 4.1} et posons
$\sigma _0:=\Im (\sigma)$. La fonction $\lambda\mapsto\mu_{nsp,\sigma }(\sigma _0\otimes\chi _{\lambda
})$ est alors r\'eguli\`ere et non nulle dans la partie r\'eelle de tout hyperplan affine qui
appartient \`a $\S $ et qui contient $A$ (i.e. qui est un \'el\'ement de $\S _A$). On peut \'ecrire
$$\mu _{sp,\sigma }(\sigma _0\otimes\chi_{\lambda })=\prod _{\alpha\in\Sigma _{sp}^+(\sigma
)}{(1-q^{-m_{\alpha }\langle\alpha ^{\vee },\lambda\rangle}) (1-q^{m_{\alpha }\langle\alpha ^{\vee
},\lambda\rangle})
\over (1-q^{k_{\alpha }-m_{\alpha }\langle\alpha ^{\vee },\lambda\rangle}) (1-q^{k_{\alpha }+m_{\alpha
}\langle\alpha ^{\vee },\lambda\rangle})}
\eqno (8.1.1)$$  avec des nombres rationnels $m_{\alpha }>0$ qui v\'erifient $m_{w\alpha }= m_{\alpha
}$ pour $w\in W(M,\o )$ (en fait $m_{\alpha }=\langle\alpha^{\vee },\ti{\alpha }\rangle^{-1}$, ce qui
est clairement invariant par $W(M,\o )$).

Posons $m=\rg _{ss}(G)-\rg _{ss}(M)$, $q_1=q^{m_{\alpha }}$ et $k_{\alpha }'=k_{\alpha }/m_{\alpha }$.
Comme, par ce qui pr\'ec\`ede, $k_{\alpha }'=k_{w\alpha }'$ pour $w\in W(M,\o )$, on peut appliquer le
lemme \`a l'expression 8.1.1.

Pour que $A\in\A_{\mu }(\S )$, il faut qu'il existe une suite strictement croissante $M\subset
M_1\subset \cdots\subset M_{m-1}\subset M_m=G$ de sous-groupes de Levi de $G$ et un polyn\^ome $p$ sur
$\o $ tel que les r\'esidus successifs de $\lambda\mapsto p\mu _{sp, \sigma }(\sigma_0\otimes\chi
_{\lambda })$ par rapport aux sous-espaces $\Re (\sigma )+a_{M_i,\sCC }^{G*}$ soient non nuls. On
d\'eduit alors du lemme que ceci implique que $\sigma $ est un p\^ole d'ordre $m$ de $\mu $. 

R\'eciproquement, le lemme montre que, si $\sigma $ est un p\^ole d'ordre $m$ de $\mu $, alors il
existe une  suite strictement croissante $M\subset M_1\subset \cdots\subset M_{m-1}\subset M_m=G$
telle que les r\'esidus successifs de $\lambda\mapsto\mu _{sp, \sigma }(\sigma_0\otimes\chi _{\lambda
})$ par rapport aux sous-espaces $\Re (\sigma )+a_{M_i,\sCC }^{G*}$ soient des fonctions non nulles qui
ont des p\^oles simples par rapport aux sous-espaces $\Re (\sigma )+a_{M_{i+1},\sCC }^{G*}$. Comme
$\lambda\mapsto\mu_{nsp,\sigma }(\sigma_0 \otimes\chi _{\lambda })$ est r\'eguli\`ere et non nulle dans
la partie r\'eelle de tout hyperplan affine de $\S $ qui contient $A$, ceci implique bien que $A\in
\A_{\mu }(\S )$.

On d\'eduit du lemme \'egalement que ces p\^oles sont d'ordre maximal.

Soit maintenant $A\in \A(\S )$, $M_A\ne G$. Tout \'el\'ement $A'$ de $\A (\S )$ qui contient $A$ est
alors dans $\A (\S _{\mu ^{M_A}})$, puisque les hyperplans affines dans $\S \setminus\S _{\mu ^{M_A}}$
sont d\'efinis \`a partir de racines dans $\Sigma (P)\setminus\Sigma (P\cap M_A)$ et ils ne
contiennent donc pas $A$. De ce que l'on vient de montrer, on d\'eduit qu'il faut et qu'il suffit que
$A\in\A_{\mu ^{M_A}}(\S _{\mu^{M_A}})$, pour que tout \'el\'ement de $A_{reg}$ soit un p\^ole de $\mu
^{M_A}$ d'ordre \'egal au rang parabolique de $M$ relatif \`a $M_A$. On est donc r\'eduit \`a \'etablir
que $A\in\A_{\mu^{M_A}}(\S _{\mu ^{M_A}})$ \'equivaut \`a $A\in\A _{\mu }(\S )$. Soit $p$ une fonction
polynomiale sur $\o $, et soit $\D $ une suite d\'ecroissante de sous-espaces affines de $\o $
contenant $A$ et d\'efinissant un op\'erateur r\'esidu en $A$ (cf. {\bf 3.7}). Soit $\sigma\in
A_{reg}$. Comme $\sigma $ est alors r\'egulier pour $\mu(\mu^{M_A})^{-1}$, on d\'eduit du lemme
appliqu\'e \`a $\mu ^{M_A}$ et $a_M^{M_A*}$ que $Res _{\sD}(p\mu )(\sigma )\ne 0$ ou $Res _{\sD }(p\mu
^{M_A})(\sigma )\ne 0$ implique que les r\'esidus successifs de $\mu ^{M_A}$ par rapport aux
sous-espaces afines formant $\D $ ont des p\^oles simples par rapport \`a ces sous-espaces affines,
d'o\`u l'\'egalit\'e $Res _{\sD }(p\mu )(\sigma )=(\mu (\mu ^{M_A})^{-1}))(\sigma ) Res _{\sD}(p\mu
^{M_A})(\sigma )$. Comme  $(\mu (\mu ^{M_A})^{-1}))$ ne s'annule par d\'efinition de $\S $ pas sur
$A_{reg }$, il en suit que $Res _{\sD}(p\mu )(\sigma )\ne 0$ \'equivaut \`a $Res _{\sD }(p\mu
^{M_A})(\sigma )\ne 0$, d'o\`u l'\'equivalence recherch\'ee entre $A\in \A_{\mu^{M_A}}(\S _{\mu
^{M_A}})$ et $A\in\A _{\mu }(\S )$.
\hfill {\fin 2}

\null
\it Remarque: \rm Pour tous $A\in\A _{\mu }(\S )$, $\sigma\in A_{\sreg }$ et toute fonction rationnelle
$\psi $ sur $\o $ r\'eguli\`ere sur $A_{\sreg }$, nous avons montr\'e au cours de la preuve de la
proposition l'identit\'e 
$$(\Res ^P_A\psi\mu )(\sigma )=\psi(\sigma )(\Res ^P_A\mu )(\sigma )\eqno (8.1.2).$$

\null
\null  Nous avons \'egalement \'etabli \`a la fin de la preuve de la proposition que

\null {\bf Corollaire:} \it Soit $M'$ un sous-groupe de Levi semi-standard de $G$, $M'\supseteq M$.
Soit $A\in\A (\S _{\mu ^{M'}})$.

Pour que $A$ soit dans $\A _{\mu ^{M'}}(\S _{\mu ^{M'}})$, il faut et il suffit que $A$ appartienne \`a
$\A _{\mu }(\S )$. \rm

\null
\null {\bf 8.2} On est maintenant en mesure d'\'enoncer et de prouver notre r\'esultat principal.
Appelons sous-groupe parabolique \it $(P,\S)$-standard de $G$ \rm tout sous-groupe parabolique
semi-standard de la forme
$P'=M'U'$, o\`u $M'=M_A$ est un sous-groupe de Levi $(P,\S)$-standard, tel que $P_A$ contienne la
chambre de Weyl associ\'ee \`a $P'$ dans $a_{M_A}^*$ (cf. {\bf 4.6}). Rappelons que l'on a not\'e (cf.
{\bf 6.1}) $\Stab (A)=\{w\in W^{M_A}(M,\o )\vert wA=A\}$ pour $A\in\A(\S )$. On d\'esignera par $[A]$
la classe de $W^{M_A}(M,\o )$-conjugaison de $A$. On fixe un ensemble de repr\'esentants $[\A(\S )]$
dans $\A(\S )$ de ces classes. On pose $[\A_{\mu }(\S )]:=[\A(\S )] \cap\A _{\mu}(\S )$.

Rappelons \'egalement que $\E _2(M',\sigma )$ d\'esigne l'ensemble des repr\'esentations de carr\'e
int\'egrable de $M'$ dont le support cuspidal contient $\sigma $ (cf. {\bf 6.1}) et que, pour
$\Omega\subseteq\Delta _{\so }$, $M_{\Omega }$ d\'esigne le plus petit sous-groupe de Levi $M'$ de $G$
contenant $M$ et tel que $\Omega\subseteq a_M^{M'*}$ (cf. {\bf 4.6}).

\null  {\bf Th\'eor\`eme:} \it Pour toute application polynomiale $\xi :\o\rightarrow i_{\ol{P}\cap
K}^KE_{\so }\otimes i_{P\cap K}^KE_{\so }^{\vee }$ \`a valeurs dans un sous-espace de dimension finie,
on a
$$\displaylinesno {
\qquad\gamma (G/M)\int _{\Re (\sigma )=r\gg _P0} \deg(\sigma ) E_P^G((J_{\ol {P}\vert P}^{-1}(\sigma )\otimes 1)\xi (\sigma ))(g)\ d_{\so }\Im (\sigma )\hfill& (8.2.1)\cr =\sum _{\Omega\subseteq\Delta _{\sso }}\sum
_{A\in[\sA _{\mu }(\sS )], M_A=M_{\Omega }}\ \vert W _{\Delta _{\sso }, M_A}\vert ^{-1}
\vert \P_{\sS }(M_A)\vert ^{-1}\gamma (G/M)\vert\Stab (A)\vert^{-1}\int _{\Re (\sigma )=r(A)}\deg (\sigma )\hfill\cr
\hfill\sum_{w'\in W^{M_A}(M,\so )}\sum _{w\in W_{M_A}^+(M,\so )}\Res
_{w'A}^P(E_P^G((J_{\ol{P}\vert P}^{-1}(ww'\sigma )\otimes 1)\xi (ww'\sigma ))(g))d_A\Im (\sigma ).\qquad\cr }$$ 
Par ailleurs, pour tout $A\in [\A_{\mu }(\S )]$ et tout $\sigma\in A$ avec $\Re (\sigma )=r(A)$, on a,
avec $M'=M_A$ et $P'=M'U'$ un sous-groupe parabolique $(P,\S )$-standard, l'identit\'e suivante de
fonctions rationnelles en $\chi '\in\X (M')$:
$$\displaylinesno {
\qquad\gamma (G/M)\deg(\sigma )\vert\Stab (A)\vert^{-1}\sum _{w'\in W^{M'}(M,\so )}\sum _{w\in
W_{M'}^+(M,\so)} 
\hfill& (8.2.2)\cr
\hfill\Res_{w'A}^P(E_P^G((J_{\overline {P}\vert P}^{-1}(ww'(\sigma\otimes\chi '))\otimes 1)\xi (ww'(\sigma\otimes\chi ')))(g))\qquad\cr 
\qquad=\gamma (G/M')\sum _{\pi\in\sE_2(M',\sigma )}\deg (\pi )\mu(\pi\otimes\chi ')E_{P'}^G(\varphi _{\xi }(P',\pi\otimes\chi '))(g).\hfill\cr 
}$$  
De plus, pour que $\sigma \in\o $ v\'erifie $\E _2(G,\sigma )\ne\emptyset $, il faut qu'il appartienne \`a un \'el\'ement $A$ de $\A _{\mu }(\S )$ qui v\'erifie $M_A=G$.

\null  
Preuve: \rm 
Pour simplifier les notations (et \'economiser de la place), on va identifier $i_{P'\cap K}^KE_{\so }\otimes i_{P\cap K}^KE_{\so }^{\vee }$ \`a un sous-espace de $\End (i_{P\cap K}^KE_{\so }, i_{P'\cap K}^KE_{\so })$ etc. via le plongement usuel. Ainsi $(J_{\ol {P}\vert P}^{-1}(\sigma )\otimes 1)\xi (\sigma )$ s'identifie \`a $J_{\ol {P}\vert P}^{-1}(\sigma )\xi (\sigma )$.

Par le th\'eor\`eme des r\'esidus {\bf 4.11},
$$\displaylinesno {
\qquad\gamma (G/M)\int _{\Re (\sigma )=r\gg _P0}\deg(\sigma ) E_P^G(J_{\overline {P}\vert
P}^{-1}(\sigma )\xi (\sigma ))(g) d_{\so }\Im (\sigma )\hfill\cr  =\sum _{\Omega\subseteq \Delta
_{\sso }}\sum _{A\in[\sA (\sS )],M_A=M_{\Omega}} \vert W _{\Delta _{\sso ,M_A}}\vert ^{-1}\
\vert\P_{\sS}(M_A)\vert^{-1}\gamma (G/M)\ \int _{\Re (\sigma )=r(A)+\epsilon _A}\deg(\sigma )\hfill&
(8.2.3)\cr
\hfill\vert\Stab (A)\vert ^{-1}\sum _{w'\in W^{M_A}(M,\so )}\sum _{w\in W_{M_A}^+(M,\so)}\Res
_{w'A}^P(E_P^G(J_{\ol{P}\vert P}^{-1}(ww'\sigma )\xi (ww'\sigma ))(g))d_A\Im (\sigma ).\qquad\cr  }$$ 
On va montrer que les contributions des termes venant d'un \'el\'ement de $[\A(\S )]-[\A_{\mu }(\S )]$
sont nulles, que les identit\'es (8.2.2) sont v\'erifi\'ees pour tout $A\in [\A(\S )]$ et que la somme
sous l'int\'egrale (8.2.3) est r\'eguli\`ere en tout $\sigma $ avec $\Re (\sigma )=r(A)$. A l'aide de
la formule de Plancherel, ceci prouvera le th\'eor\`eme.

On va effectuer une r\'ecurrence sur le rang semi-simple de $G$. Les diff\'erentes \'etapes de la
preuve (traitement de la partie "continue" de (8.2.3) \`a l'aide de l'hypoth\`ese de r\'ecurrence,
identification de la partie continue de (8.2.3) avec celle de la formule de Plancherel, traitement de
la partie "discr\`ete")  sont reparties sur les num\'eros {\bf 8.3}, {\bf 8.4} et {\bf 8.5} ci-apr\`es.

\null
\null  {\bf 8.3} On va effectuer une r\'ecurrence sur le rang semi-simple de $G$. Si $rg_{ss}(G)=0$, il
n'y a rien \`a montrer. Supposons le th\'eor\`eme prouv\'e pour tout groupe de rang semi-simple plus
petit que celui de $G$. Soit $\Omega \subseteq\Delta _{\so }$ tel que $M':=M_{\Omega }\ne G$. Soit
$A\in[\A (\S )]$ avec $M'=M_{A}$ et notons $P'$ un sous-groupe parabolique $(P,\S )$-standard de Levi
$M'$. Posons $P_1=(M'\cap P)U'$. C'est un sous-groupe parabolique semi-standard contenu dans $P'$ de
sous-groupe de Levi $M$. Comme $\Delta _{\so }\subseteq\Sigma (P)$ et que $P'$ est $(P,\S)$-standard,
tout \'el\'ement $\alpha\in\Delta _{\so }$ v\'erifie ou $\alpha\in \Sigma (P\cap M')$ ou
$\alpha_{\vert M'}\in\Sigma (P')$. On en d\'eduit $\Delta _{\so }\subseteq \Sigma (P_1)$.

On a alors pour $\sigma\in \o $ en position g\'en\'erique, en utilisant les r\`egles de composition,
d'adjonction et de compatibilit\'e vis-\`a-vis de l'induction parabolique pour les op\'erateurs
d'entrelacement (o\`u on identifie $i_{P_1}^G$ avec $i_{P'}^Gi_{P_1\cap M'}^{M'}$) (cf. [W] chapitre
IV)
$$\displaylines {
\ \ \sum _{w\in W_{M'}^+(M,\so)} E_P^G(J_{P\vert\ol{P}}(w\sigma )\xi(w\sigma))(g)\hfill\cr  =\sum
_{w\in W_{M'}^+(M,\so)} E_P^G(\lambda (w)\lambda (w^{-1})J_{P\vert\ol{P}}(w\sigma )\xi (w\sigma
))(g)\hfill\cr =\sum _{w\in W_{M'}^+(M,\so)} E_P^G(\lambda (w)J_{w^{-1}P\vert\ol{w^{-1}P}}(\sigma
)\lambda (w^{-1})\xi (w\sigma ))(g)\hfill\cr  =\sum _{w\in W_{M'}^+(M,\so)}\mu ^{M'}_{w^{-1}P\cap
M'\vert P\cap M'}(\sigma ) E_{P_1}^G(\lambda(w)J_{w^{-1}P\vert P_1}(\sigma )J_{P_1\vert
(M'\cap\ol{P})U'}(\sigma )\hfill\cr
\hfill J_{(M'\cap\ol{P})U'\vert\ol{w^{-1}P}}(\sigma )\lambda (w^{-1})\xi (w\sigma ))(g)\qquad\cr 
=\sum _{w\in W_{M'}^+(M,\so)} E_{P_1}^G((J_{P_1\vert(M'\cap\ol{P})U'}(\sigma ) \otimes 1)[\mu
^{M'}_{w^{-1}P\cap M'\vert P\cap M'}(\sigma )\hfill\cr
\qquad\hfill (J_{(M'\cap\ol{P})U'\vert\ol{w^{-1}P}}(\sigma )\lambda (w^{-1})\otimes J_{P_1\vert
w^{-1}P}(\sigma ^{\vee })\lambda (w^{-1}))\xi (w\sigma )])(g)\cr  =\sum _{w\in W_{M'}^+(M,\so)}
E_{P'}^G(i_{P'}^G(J^{M'}_{M'\cap P\vert M'\cap\ol{P}}(\sigma )\otimes 1)[\mu ^{M'}_{w^{-1}P\cap
M'\vert P\cap M'}(\sigma )\hfill\cr
\qquad\hfill (J_{(M'\cap\ol{P})U'\vert\ol{w^{-1}P}}(\sigma )\lambda (w^{-1})\otimes J_{P_1\vert
w^{-1}P}(\sigma ^{\vee })\lambda (w^{-1}))\xi (w\sigma )])(g)\cr  =\int _{P'\backslash G} E_{P\cap
M'}^{M'}(J_{P\cap M'\vert\ol{P}\cap M'}^{M'}(\sigma )[\sum _{w\in W_{M'}^+(M,\so)}\mu
^{M'}_{w^{-1}P\cap M'\vert P\cap M'}(\sigma )\hfill\cr
\qquad\hfill((J_{(M'\cap\ol{P})U'\vert\ol{w^{-1}P}} (\sigma )\lambda (w^{-1})\otimes J_{P_1\vert
w^{-1}P}(\sigma ^{\vee })\lambda (w^{-1}))\xi (w\sigma ))(\ol{g'}g,\ol{g'})])(1)d\ol{g'} \qquad\cr }$$

Posons $\ti {P}_1=(M'\cap P)\ol{U'}$. Comme l'int\'egrale est en fait une somme finie, il s'en suit
que, pour $A'$ dans $[A]$ et $\sigma\in A'$ en position g\'en\'erique 
$$\displaylinesno {
\qquad\sum _{w\in W_{M'}^+(M,\so)}\Res_{A'}(E_P^G(J_{\ol{P}\vert P}^{-1}(w\sigma )\xi (w\sigma
))(g))=\hfill&(8.3.1)\cr
\qquad\int _{P'\backslash G}\Res_{A'}(E_{P\cap M'}^{M'}(J_{\ol{P}\cap M'\vert P\cap M'}^{M'\
-1}(\sigma )[\mu _{P_1\vert\ti{P}_1}(\sigma )\sum _{w\in W_{M'}^+(M,\so)}\mu ^{M'}_{w^{-1}P\cap
M'\vert P\cap M'}(\sigma )\hfill\cr
\hfill(J_{(\ol{P}\cap M')U'\vert\ol{w^{-1}P}}(\sigma )\lambda (w^{-1})\otimes J_{P_1\vert
w^{-1}P}(\sigma ^{\vee })\lambda (w^{-1}))\xi (w\sigma )](\ol{g'}g,\ol{g'}))(1)) d\ol{g'}\qquad\cr }$$
L'application rationnelle $\ti {\xi }^{M'}_{g'g,g'}:\o\rightarrow i_{\ol{P}\cap K\cap M'}^{K\cap
M'}E_{\so }\otimes i_{P\cap K\cap M'}^{K\cap M'}E_{\so }^{\vee }$, 
$$\displaylines{
\qquad\sigma\mapsto\mu_{P_1\vert\ti{P_1}}(\sigma )[\sum _{w\in W_{M'}^+(M,\so)} \mu ^{M'}_{w^{-1}P\cap
M'\vert P\cap M'}(\sigma)\hfill\cr
\hfill ((J_{(\ol{P}\cap M')U'\vert\ol{w^{-1}P}}(\sigma )\lambda (w^{-1})\otimes J_{P_1\vert
w^{-1}P}(\sigma ^{\vee })\lambda (w^{-1}))\xi (w\sigma )](g'g,g') \qquad  }$$  v\'erifie d'apr\`es
{\bf 6.2} les hypoth\`eses de la proposition {\bf 6.1}. Le lemme {\bf 6.2} montre par ailleurs que
$\varphi _{\ti {\xi }^{M'}_{g'g,g'}} (\sigma )=\mu _{P_1\vert\ti {P_1}}(\sigma )(\varphi _{\xi
}(P_1,\sigma ) (g'g,g'))$ (dans les notations de {\bf 6.2}). 

Fixons $\sigma\in A$ avec $\Re (\sigma )=r(A)$. On trouve alors par la proposition {\bf 6.1} avec
l'hypoth\`ese de r\'ecurrence, en utilisant l'\'egalit\'e $\gamma (G/M)=\gamma (G/M')\gamma (M'/M)$, et
en posant $\sigma '=\sigma\otimes\chi '$ l'\'egalit\'e 
$$\displaylinesno {
\qquad\vert\Stab (A)\vert^{-1}\sum _{w'\in W^{M'}(M,\so )}\gamma (G/M)\deg (\sigma )\hfill&(8.3.2)\cr
\hfill\Res_{w'A}^{P\cap M'}(E_{P\cap M'}^{M'}(J_{P\cap M'\vert\ol{P}\cap M'}^{M'\ -1}(w'\sigma ')[\mu
_{P_1\vert\ti{P_1}}(w'\sigma ')\sum _{w\in W_{M'}^+(M,\so )} \mu ^{M'}_{w^{-1}P\cap M'\vert P\cap
M'}(w'\sigma ')\qquad\cr
\hfill (J_{(\ol{P}\cap M')U'\vert\ol{w^{-1}P}}(w'\sigma ')\lambda (w^{-1})\otimes J_{P_1\vert
w^{-1}P}(w'\sigma '^{\vee })\lambda (w^{-1}))\xi (ww'\sigma ')] (g'g,g'))(1))\cr
\qquad =\gamma (G/M')\sum _{\pi\in\sE_2(M',\sigma )}\deg (\pi)\mu _{P_1\vert\ti {P_1}}(\sigma
')E_{M'}^{M'}(\varphi _{\xi }(P_1,\pi \otimes\chi ')(g'g,g'))(1). \hfill\cr  }$$  de fonctions
rationnelles en $\chi '\in\X (M')$. Remarquons que, par hypoth\`ese de r\'ecurrence, $\E _2(M',\sigma
)$ est vide, si $A\not\in\A_{\mu ^{M'}} (\S_{\mu ^{M'}})$, ce qui \'equivaut par le corollaire {\bf
8.1} \`a $A\not\in\A_{\mu }(\S)$. La somme (8.3.2) est donc nulle, si $A\in [\A (\S )]-[\A _{\mu }(\S
)]$.

En appliquant {\bf 4.9}, en int\'egrant l'\'egalit\'e ci-dessus sur $P'\backslash G$ par rapport \`a
$\ol{g'}$, en ins\'erant le r\'esultat dans l'identit\'e (8.3.1) et en utilisant l'identit\'e $\mu
_{P_1\vert\ti{P_1}} (\sigma\otimes\chi ')=\mu (\pi\otimes\chi' )$, on en d\'eduit l'identit\'e 
$$\displaylinesno {
\qquad\gamma (G/M)\deg(\sigma )\vert\Stab (A)\vert^{-1}\sum _{w'\in W^{M'}(M,\so )}\sum _{w\in
W_{M'}^+(M,\so)}
\hfill&(8.3.3) \cr
\hfill\Res_{w'A}^P(E_P^G(J_{\ol{P}\vert P}^{-1}(ww'(\sigma\otimes\chi '))
\xi (ww'(\sigma\otimes\chi ')))(g))\qquad\cr
\hfill=\gamma (G/M')\sum _{\pi\in\sE_2(M',\sigma )}\mu (\pi\otimes\chi' )\deg (\pi )E_{P'}^G(\varphi
_{\xi }(P',\pi\otimes\chi '))(g),\qquad\cr  }$$  de fonctions rationnelles sur $\X (M')$, o\`u
$\E_2(M',\sigma )=\emptyset $, si $A\not\in\A_{\mu }(\S )$. Ceci implique bien les identit\'es (8.2.2)
du th\'eor\`eme pour tout $\Omega\subseteq\Delta_{\so }$, $M_{\Omega }\ne G$ et que les termes dans
(8.2.3) index\'es par $[\A (\S )]-[\A _{\mu }(\S )]$ sont nuls. Comme la fonction $\mu $ est
r\'eguli\`ere en une repr\'esentation de carr\'e int\'egrable (cf. [W] IV.3(6)), la fonction d\'efinie
par (8.3.3) est r\'eguli\`ere en tout $\chi '$ avec $\Re (\chi ')=0$. On peut donc poser $\epsilon
_A=0$ dans (8.2.3).

\null
\null {\bf 8.4} On va maintenant montrer que la somme des termes dans (8.2.3) index\'es par les
sous-ensembles $\Omega $ de $\Delta _{\so }$, $M_{\Omega }\ne G$, s'identifie \`a la partie continue de
la formule de Plancherel. 

Soit $\Omega\subseteq\Delta _{\so }$ tel que $M_{\Omega }\ne G$, et supposons qu'il existe $A\in[\A (\S
)]$ tel que $M_A=M_{\Omega }$. D\'esignons par $[A]_{P,\sS}$ l'ensemble des $A'\in [\A (\S )]$ qui sont
conjugu\'es \`a $A$ par un \'el\'ement de $W(M,\o )$ et qui sont tels que $M_{A'}$ soit $(P,\S
)$-standard. (L'ensemble $[A]_{P,\sS}$ n'est qu'un sous-ensemble de celui d\'esign\'e par le m\^eme
symbole \`a l'int\'erieur de la preuve du lemme {\bf 4.7}.) Fixons pour tout $A'$ de $[A]_{P,\sS}$ un
\'el\'ement $\sigma _{A'}$ de $A'$ avec $\Re (\sigma _{A'})=r(A')$, et un sous-groupe parabolique
$(P,\S )$-standard  de Levi $M_{A'}$ que l'on notera (abusivement) $P_{A'}$. On peut supposer que les
repr\'esentations $\sigma _{A'}$ sont mutuellement conjugu\'ees. Comme les identit\'es (8.2.2) ont
\'et\'e \'etablies pour les \'el\'ements de $[A]_{P,\sS}$, la somme des termes dans (8.2.3), index\'es
par les \'el\'ements de $[A]_{P,\sS}$, est \'egale \`a
$$\displaylinesno {
\qquad\sum _{A'\in [A]_{P,\ssS}} \vert W _{\Delta _{\sso},M_{A'}}\vert^{-1}\vert
\P_{\sS}(M_{A'})\vert^{-1}\gamma (G/M_{A'})\hfill&(8.4.1)\cr
\hfill\sum_{\pi\in\sE_2(M_{A'},\sigma _{A'})}\int _{\so _{\pi ,0}}\deg (\pi')  E_{P_{A'}}^G(\varphi
_{\xi }(P_{A'},\pi'))(g)\mu (\pi ')d\pi'.\qquad\cr  }$$  Remarquons que les constantes $\vert W
_{\Delta _{\sso},M_{A'}}\vert$, $\vert\P_{\sS}(M_{A'})\vert$ et $\deg(\sigma _{A'})$ ne d\'ependent pas
du choix de $A'\in [A]_{P,\sS}$ et que, par {\bf 4.7}, $\vert\{M_{A'}\vert A'\in
[A]_{P,\sS}\}\vert={\vert\sP_{\ssS}(M_A)\vert\ \vert W _{\Delta _{\sso },M_A}\vert\over\vert W(M_A,\so
)\vert}$. Par ailleurs, on a, pour $w\in W(M,\o)$, $\E_2(wM_A, w\sigma
_A)=\{w\pi\vert\pi\in\E_2(M_A,\sigma _A)\}$ et $\gamma (G/M_A)=\gamma (G/wM_A)$. On en d\'eduit avec
[W] V.3.1 que les termes dans la somme (8.4.1) ne d\'ependent pas du choix de $A'\in [A]_{P,\sS}$. Par
suite, (8.4.1) est \'egal \`a
$$\displaylinesno {
\qquad\gamma (G/M_A)\vert W(M_A,\o )\vert^{-1}\sum _{A'\in [A]_{P,\ssS},
M_{A'}=M_A}\sum_{\pi\in\sE_2(M_{A'},\sigma _{A'})}\hfill&(8.4.2)\cr
\hfill \int _{\so_{\pi,0}}\deg(\pi ')E_{P_{A'}}^G(\varphi _{\xi } (P_{A'},\pi'))(g)\mu (\pi
')d\pi'.\qquad   }$$ 

Soit $\pi\in \E_2(M_A,\sigma _A)$. On peut identifier $W(M_A,\o _{\pi, 0})$ \`a un sous-groupe de
$W(M_A, \o)$: soit $w\in W(M_A,\o _{\pi, 0})$. Alors $w\sigma _A$ est dans le support cuspidal de $w\pi
$ et, comme $w\pi\in\o_{\pi ,0}$, il existe $w'\in W^{M_A}$ tel que $ww'\sigma _A\in\o $, d'o\`u
$ww'\o=\o $. La classe \`a droite modulo $W^{M_A}(M,\o )$ de $w'$ est d\'etermin\'ee de fa\c con
unique. 

L'ensemble des $\o _{\pi ',0}$ avec $\pi '$ dans la r\'eunion des $\E _2 (M_{A'},\sigma _{A'})$ avec
$A'\in [A]_{P,\sS }$, $M_{A'}=M_A$ et $\o _{\pi',0}$ conjugu\'es \`a $\o _{\pi ,0}$ par un \'el\'ement
de $W$ est sans multiplicit\'e, puisque $A'$ ne peut \^etre conjugu\'e \`a $A$ par un \'el\'ement de
$W^{M_A}(M,\o )$ (les \'el\'ements de $[A]_{P,\sS }$ \'etant des repr\'esentants de ces classe de
conjugaison). On v\'erifie par ailleurs que cet ensemble est \'egal \`a l'orbite de $\o _{\pi ,0}$ pour
l'action par conjugaison par $W(M_A, \o)$. Inversement, toute orbite $\o _{\pi ',0}$ avec $\pi '$ dans
la r\'eunion des $\E _2(M_{A'},\sigma _{A'})$, $A'\in [A]_{P,\sS }$ et $M_{A'}=M_A$, est conjugu\'ee
par un \'el\'ement de $W$ \`a une orbite de la forme $\o _{\pi ,0}$ avec $\pi\in\E _2 (M_A,\sigma _A)$. 

En rassemblant les repr\'esentations $\pi '$ dans la r\'eunion des $\E _2 (M_{A'},\sigma _{A'})$ avec
$A'\in [A]_{P,\sS }$ et $M_{A'}=M_A$, dont les orbites $\o _{\pi ',0}$ sont conjugu\'ees et qui ont par
[W] V.3.1 la m\^eme contribution dans (8.4.2), l'identit\'e (8.4.2) devient donc
$$\gamma (G/M_A)\sum _{\pi\in\sE_2(M_A,\sigma _A)}\vert W(M_A,\o_{\pi,0})\vert^{-1} \int_{\so_{\pi,
0}}\deg(\pi')E_{P_A}^G(\varphi _{\xi }(P_A,\pi'))(g)d\pi'.\eqno (8.4.3)$$ 

Notons $L^2(G)_{[A]_{P,\ssS}}$ le sous-espace de l'espace hilbertien $L^2(G)$ engendr\'e par les
combinaisons lin\'eaires des fonctions de la forme
$$g\mapsto\int _{\so _{\pi,0}}E_{P_A}^G(\varphi _{\xi }(P_A,\pi '))(g)d\pi ', \eqno (8.4.4)$$  avec
$\pi\in\E_2(M_A,\sigma _A)$. En comparant (8.4.3) avec (2.5.1), on voit en utilisant la proposition
VII.2.2 de [W] que (8.4.3) est la projection orthogonale de $f_{\xi }$, consid\'er\'e comme \'el\'ement
de $L^2(G)$, sur $L^2(G)_{[A]_{P,\ssS}}$. 

Si on change de classe $[A]_{P,\sS}$, (8.4.3) d\'efinit une fonction qui est orthogonale \`a l'espace
$L^2(G)_{[A]_{P,\ssS}}$, puisque des fonctions de la forme (8.4.4) correspondant \`a des orbites
unitaires non conjugu\'ees sont orthogonales (cf. [W] proposition VII.2.2).

Pour conclure que la somme des termes dans (8.2.3) index\'es par les $\Omega \subseteq\Delta _{\so }$,
$M_{\Omega }\ne G$, est \'egale \`a la partie continue de la formule de Plancherel, il reste \`a
v\'erifier que, pour tout
$(P'=M'U',\o _0')\in\Theta _2(\o)$ avec $M'\ne G$, il existe $A\in[\A_{\mu }(\S )]$, $\sigma _A\in A$
et $\pi\in\E_2(M_A,\sigma _A)$, tel que $M_A$ soit $(P,\S )$-standard et que $\o _0'$ soit conjugu\'e
\`a $\o _{\pi ,0}$. 

Par hypoth\`ese de r\'ecurrence, pour que $\E _2(M',\sigma )\ne\emptyset $ pour un $\sigma\in \o $, il
faut qu'il existe $A\in\A_{\mu ^{M'}}(\S _{\mu ^{M'}})\subseteq \A(\S )$ tel que $M'=M_A$. Par le lemme
{\bf 4.6}, $M'$ est conjugu\'e \`a un sous-groupe de Levi $(P,\S )$ -standard de $G$ par un \'el\'ement
$w$ de $W(M,\o )$. Soit $\pi $ dans $w\o _0'$. Par hypoth\`ese de r\'ecurrence, il existe alors
$A'\in[\A _{\mu^{M'}}(\S _{\mu ^{M'}})]\subseteq [\A _{\mu}(\S )]$, $\sigma _{A'}\in A'$ tels que
$\pi\in\E_2(M_{A'},\sigma _{A'})$. 

Tous les termes de la partie continue de la formule de Plancherel apparaissent donc bien dans la partie
continue de (8.2.3).

\null
\null  {\bf 8.5} Supposons maintenant $M_{\Delta _{\sso }}=G$ et consid\'erons les termes dans (8.2.3)
qui sont index\'es par $\Delta _{\so }$. La fonction 
$f_{\xi }:G\rightarrow\CC $, $g\mapsto \int_{\Re (\sigma )=r\gg _P0} E_P^G(J_{\ol{P}\vert
P}^{-1}(\sigma )\xi (\sigma ))(g^{-1})d_{\so }\Im (\sigma)$ est lisse \`a support compact, donc de
carr\'e int\'egrable. D'apr\`es [W] proposition VI.3.1, pour $M'$ un sous-groupe de Levi $\supseteq M$,
$\pi'\in\E_2(M')$, la fonction
$$g\mapsto \int _{\so _{\pi,0}}E_{P'}^G(\varphi _{\xi }(\pi ',P'))(g)\mu (\pi ')d\pi'$$  est dans
l'espace de Schwartz-Harish-Chandra et donc \'egalement de carr\'e int\'egrable. On d\'eduit alors de
(8.2.3) et des identit\'es (8.2.2) que l'on vient d'\'etablir pour $M_{\Omega }\ne G$ que la fonction
sur $G$ d\'efinie pour $g\in G$ par 
$$\displaylinesno{
\qquad\sum _{A\in [\sA (\sS)], M_A=G}\gamma (G/M)\int _{\Re(\sigma ')=r(A)+ \epsilon _A}\deg(\sigma
')\hfill&(8.5.1)\cr
\hfill\vert\Stab (A)\vert^{-1} \sum _{w'\in W^(M,\so )} \Res_{w'A}^P(E_P^G(J_{\ol{P}\vert
P}^{-1}(w'\sigma ')\xi (w'\sigma '))(g))d_A\Im (\sigma ')\qquad\cr }$$  est de carr\'e int\'egrable
pour toute application polynomiale $\xi $ sur $\o$ et, comme la somme des termes index\'es par les
$\Omega\subseteq\Delta_{\so }$, $M_{\Omega }\ne G$, dans (8.2.3) s'identifie \`a la partie continue de
la formule de Plancherel, l'expression (8.5.1) est \'egale \`a 
$$\sum _{(G,\so_0')\in\Theta _2(\so )}\int _{\sso_0'}\deg(\pi') E_G^G(\varphi_{\xi }(G,\pi
'))(g)d(\pi').\eqno{(8.5.2)}$$

Remarquons d'abord que, si $M_A=G$, alors $\P_{\sS}(M_A)=\{a_G^*\}$. Par suite, on peut choisir
$\epsilon _A=0$ dans (8.5.1).

Soit $A\in [\A (\S )]$ et $\sigma\in A$ avec $\Re (\sigma )=r(A)$. Il faut \'etablir les identit\'es
(8.2.2) et montrer que la valeur de chaque c\^ot\'e est nulle, si $A\in [\A(\S )]-[\A _{\mu }(\S )]$.
Notons $\chi _{\sigma }$ la restriction \`a $T_G$ du caract\`ere central de $\sigma $ qui est par choix
de $\sigma $ un caract\`ere unitaire. Rappelons l'application
$p_{\chi _{\sigma }}$ d\'efinie dans [W] VIII.4: elle associe \`a une fonction $f$ de l'espace de
Schwartz-Harish-Chandra la fonction $p_{\chi _{\sigma }}(f):G\rightarrow \CC $ bi-invariante par un
sous-groupe ouvert compact d\'efinie par $p_{\chi _{\sigma }}(f)(g)=\int _{T_G}\chi _{\sigma }^{-1}(t)
f(tg)dt$. Elle v\'erifie $p_{\chi _{\sigma }}(f)(tg)=\chi _{\sigma }(t)f(g)$ pour tout $t\in T_G$ et
$g\in G$. Comme les fonctions d\'efinies sur $G$ par les identit\'es (8.5.1) et (8.5.2) sont dans
l'espace de Schwartz-Harish-Chandra, on peut leur appliquer la transformation $p_{\chi _{\sigma }}$.

Notons  $\E _2(G,\o)$ l'ensemble des repr\'esentations de carr\'e int\'egrable de $G$ dont le support
cuspidal contient un \'el\'ement de $\o $, et $\E _2 (G,\o)_{\chi _{\sigma }}$ (resp. $A_{\chi _{\sigma
}}$) le sous-ensemble de $\E _2(G,\o )$ (resp. d'un \'el\'ement $A$ de $\A(\S )$), form\'e des
repr\'esentations $\pi '$ (resp. $\sigma '$) dont la restriction \`a $T_G$ agit par le caract\`ere
$\chi _{\sigma }$. En refaisant les arguments de la preuve de [W] th\'eor\`eme VIII.4.2, on voit que la
fonction d\'efinie par 8.5.2 devient, si on lui applique $p_{\chi _{\sigma }}$, \'egale \`a celle
d\'efinie pour $g\in G$ par 
$$\sum _{\pi'\in\sE_2(G,\so )_{\chi _{\sigma }}}\deg(\pi ')E_G^G(\varphi _{\xi }(G,\pi '))(g).\eqno
(8.5.3)$$  Remarquons que l'espace $A$ est form\'e des $\sigma\otimes\chi' $ avec $\chi'\in\X (G)$.
Comme la valeur de $J_{P\vert\ol{P}}$ en une repr\'esentation ne change pas si on tensorise celle-ci
par un \'el\'ement de $\X (G)$, les fonctions sous l'int\'egrale (8.5.1) sont polynomiales. On peut
donc refaire de nouveau les arguments de [W] th\'eor\`eme VIII.4.2 pour montrer que la fonction
d\'efinie par (8.5.1) devient, si on lui applique $p_{\chi _{\sigma }}$, \'egale \`a celle d\'efinie
pour $g\in G$ par 
$$\displaylinesno {
\gamma (G/M)\vert\Stab (A)\vert ^{-1}\sum _{\sigma'\in A_{\chi _{\sigma }}}\deg (\sigma ')\sum _{w\in
W(M,\so )} \Res _{wA}^P(E_P^G(J_{\ol{P}\vert P}^{-1}(w\sigma') \xi (w\sigma
'))(g)).\qquad\qquad&(8.5.4)\cr  }$$  Choisissons un polyn\^ome $p$ sur $\o $ invariant pour l'action
par $W(M,\o)$ et pour celle par $\X (G)$, non nul en $\sigma $ et qui s'annule en un ordre \'elev\'e en
tout $\sigma '\in\o$ qui ne soit pas conjugu\'e \`a $\sigma $ et qui appartienne ou au support cuspidal
d'un \'el\'ement $\pi '$ de $\E_2(G)_{\chi _{\sigma }}$ avec $\varphi _{\xi }(G,\sigma ')\ne 0$ ou \`a
$A_{\chi _{\sigma }}$ pour un $A\in [\A(\S )]$ v\'erifiant $M_A=G$. L'ensemble de ses $\sigma '\in\o $
est fini par d\'efinition de $\A (\S )$ et par [W] VIII.1.1 (1). On peut par ailleurs supposer que les
d\'eriv\'ees de $p$ en $\sigma $ s'annulent \`a un ordre suffisamment \'elev\'e. (L'existence d'une
fonction polynomiale $p_1$ sur $\X (M)$ invariante pour l'action par $\X (G)$ s'annulant \`a un ordre
\'elev\'e en un nombre fini de points donn\'es et ne s'annulant pas sur l'ensemble fini correspondant
aux conjugu\'es de $\sigma $ est \'el\'ementaire, puisque le quotient de $\X (M)$ par $\X (G)$ est un
produit de tores $\CC ^{\times }$ (cf. {\bf 1.2}). En prenant le produit de $p_1$ avec ses conjugu\'es,
on en d\'eduit une fonction polynomiale $p$ sur $\o $ qui est invariante pour l'action $W(M,\o )$ et par
$\X (G)$ et qui a toutes les propri\'et\'es demand\'ees sauf peut-\^etre que ses d\'eriv\'ees en
$\sigma $ ne s'annulent pas \`a un ordre suffisamment \'elev\'e. En remplacant $p$ par son produit avec
$-p+2p(\sigma )$ qui est un polyn\^ome $W(M,\o )$-invariant sur $\o $ et en r\'ep\'etant ce proc\'ed\'e
un nombre suffisant de fois, on aboutit finalement \`a une fonction polynomiale $p$ qui poss\`ede
toutes les propri\'et\'es demand\'ees.)

On d\'eduit de (8.5.3) et (8.5.4), en calculant $\Res _{wA}^P$ \`a l'aide d'op\'erateurs
diff\'erentiels holomorphes \`a coefficients constants (cf. {\bf 3.7}) en utilisant les r\`egles de
Leibniz, l'\'egalit\'e
$$\displaylinesno {
\qquad\sum _{\pi '\in\sE_2(G,\sigma )}\deg (\pi')E_G^G(\varphi _{p\xi } (G,\pi'))(g)\hfill&(8.5.5)\cr
\hfill =\gamma(G/M)\deg(\sigma )\vert\Stab (A)\vert^{-1}\sum _{w\in W(M,\so )}
\Res _{wA}^P(E_P^G(J_{\ol{P}\vert P}^{-1}(w\sigma )((p\xi )(w\sigma )))(g))\cr }$$  pour toute
application polynomiale $\xi:\o\rightarrow i_{\ol{P}\cap K}^K E_{\so }\otimes i_{P\cap K}^KE_{\so
}^{\vee }$ \`a image dans un espace de dimension finie, o\`u on a utilis\'e que $\varphi _{p\xi }(\pi
',G)=p(\sigma ')\varphi _{\xi }(\pi ',G)$ pour tout $\pi '$ dans $\E_2(G,\sigma ')$ et tout $\sigma '$
dans la r\'eunion de $A_{\chi _{\sigma }}$ et de $\E _2(G, \o)_{\chi _{\sigma }}$, ainsi que le fait
que les d\'eriv\'ees de $p$ en $\sigma '$ s'annulent \`a un ordre \'elev\'e. Les m\^emes arguments
montrent que l'on peut "diviser par $p$" pour obtenir l'identit\'e (8.2.2) en $\sigma $ avec $M'=G$ (et
$\chi'=1$). On en d\'eduit que l'identit\'e (8.2.2) est v\'erifi\'ee relative \`a $\sigma $ pour tout
$\chi '\in\Xim (G)$ et donc, puisque les deux applications dans (8.2.2) sont rationnelles, en tout
$\chi'\in\X (G)$.

Choisissons $\xi =\xi _{P,\so }^H$ (cf. {\bf 7.1}). Comme l'application $\o\rightarrow \CC,\
\sigma'\mapsto E_P^G(J_{P\vert\ol{P}}(\sigma ')\xi _{P,\so }^H$ $(\sigma '))(1)$ est constante (cf.
{\bf 7.2}) et que $E_P^G(J_{\ol {P}\vert P}^{-1}(\sigma ')\xi _{P,\so }^H(\sigma '))=
E_P^G(J_{P\vert\ol{P}}(\sigma ')\xi _{P,\so }^H (\sigma '))\mu (\sigma ')$, on voit par d\'efinition de
$\A_{\mu }(\S )$ que 
$$\Res_A^P((E_P^G(J_{\ol {P}\vert P}^{-1}(\sigma )\xi _{P,\so }^H(\sigma )) (1))=0,$$   si $A\in\A
(\S)-\A_{\mu }(\S )$. On en d\'eduit que chacun des termes dans (8.2.2) \'evalu\'es en $g=1$ est nul si
$A\in [\A (\S)]-[\A_{\mu }(\S )]$, i.e. que l'on a 
$$\sum _{\pi\in\sE_2(G,\sigma )} \deg(\pi) \tr(\pi(f_{\varphi _{\sso }}))=0,$$ avec $f_{\varphi_{\sso
}}=f_{\xi _{P, \sso}^H}$ (cf. {\bf 2.2}, {\bf 2.3} et {\bf 7.1} pour les notations). Comme
$\tr(\pi(f_{\varphi _{\sso }}))$ est par {\bf 7.1} un entier strictement positif pour tout
$\pi\in\E_2(G,\sigma )$, il s'en suit que $\E_2(G,\sigma )= \emptyset $. Ceci prouve les derni\`eres
assertions du th\'eor\`eme.\hfill{\fin 2}

\null
\null {\bf 8.6 Corollaire:} \it Pour qu'un \'el\'ement $\sigma $ de $\o $ appartienne au support
cuspidal d'une repr\'esentation de carr\'e int\'egrable de $G$, il faut que $\sigma $ soit un p\^ole de
la  fonction $\mu $ de Harish-Chandra d'ordre \'egal au rang parabolique de $M$. Ces p\^oles sont
d'ordre maximal. R\'eciproquement, un p\^ole $\sigma $ de la fonction $\mu $ de Harish-Chandra d'ordre
\'egal au rang parabolique de $M$ correspond ainsi \`a une repr\'esentation de carr\'e int\'egrable de
$G$, si et seulement si l'espace affine  $A=\o_{\sigma, G}$ v\'erifie $\sum _{w\in W(M,\so )}(\Res
_{wA}^P\mu_{sp,\sigma })(w\sigma )\ne 0$ et si $\Re (\sigma )=r(A)$.

\null Preuve: \rm La premi\`ere partie est une cons\'equence directe de la derni\`ere assertion du
th\'eor\`eme {\bf 8.2} avec la caract\'erisation de $\A _{\mu}(\S)$ donn\'ee par la proposition {\bf
8.1}. Concernant la deuxi\`eme assertion, l'identit\'e (8.2.2) du th\'eor\`eme {\bf 8.2} et le lemme
{\bf 7.2} prouvent qu'avec les notations du paragraphe {\bf 7} et $A=\o _{\sigma ,G}$
$$\eqalignno { &\gamma (G/M)\deg(\sigma )\vert\Stab (A)\vert^{-1}E_P^G((J_{P\vert\overline {P}}(\sigma )\otimes 1)\xi _{P,\so }^H(\sigma ))(1)
\sum _{w\in W(M,\so )} (\Res _{wA}^P\mu )(w\sigma )\cr =&\sum _{\pi\in\sE_2(G,\sigma )}\deg (\pi
)\tr(\pi(f_{\varphi _{\sso }})). &(8.6.1)\cr }$$  si $\Re (\sigma )=r(A)$. Comme
$E_P^G((J_{P\vert\overline {P}}(\sigma )\otimes 1)\xi _{P,\so }^H(\sigma ))(1)$ et $\tr(\pi(f_{\varphi _{\sso
}}))$ sont des nombres $>0$ par les lemmes {\bf 7.1} et {\bf 7.2} respectivement, $\E_2(G,\sigma
)\ne\emptyset $  \'equivaut \`a $\sum _{w\in W(M,\so )}$ $(\Res _{wA}^P\mu )(\sigma )\ne 0$. On peut
finalement remplacer $\mu $ par $\mu _{sp, \sigma }$ suite \`a l'identit\'e (8.1.2), puisque $\mu
_{nsp, \sigma }$ est r\'eguli\`ere et non nulle en $A_{\sreg }$ et invariante par $W(M,\o )$ (cf. {\bf
4.1}).

Il reste \`a remarquer qu'une condition n\'ecessaire, pour que $\sigma $ soit dans le support cuspidal
d'une repr\'esentation de carr\'e int\'egrable de $G$, est que $\sigma $ agit sur $T_G$ par un
caract\`ere unitaire, ce qui \'equivaut \`a $\Re (\sigma )=r(A)$. \hfill{\fin 2}

\null
\it Remarque: \rm La somme $\sum _{w\in W(M,\so )}(\Res_{wA}^P\mu_{sp,\sigma })(w.)$ comporte
beaucoup de termes qui sont en effet nuls. D'abord, comme on l'a d\'ej\`a remarqu\'e en {\bf 4.12}, on
devrait avoir $\Res _{wA}^P\equiv 0$, si $r(wA)$ n'est pas une somme d'\'el\'ements de $\Delta _{\so
}$ \`a coefficients $\geq 0$. Si $\Re (\sigma )$ est dans la composante $(P,\S )$-standard $P_A$ de
$\P_{\sS}(M_A)$ (d\'efinie dans {\bf 4.6}) (Silberger appelle $\sigma $ dans [S3] alors un point de
Casselman), il est facile de voir (cf. par exemple les arguments utilis\'es pour la remarque 3.16 dans
[HO1]) que $\Res _{wA}^P\mu =0$ si $w\not\in\Stab (A)$. Par ailleurs, $J_{P\vert\ol{P}}$ est dans ce
cas r\'egulier en $\sigma $. Notons $\iota $ l'inclusion de l'image de $J_{P\vert\ol{P}}(\sigma )$ dans
$i_{P\cap K}^KE$ et $p$ un projecteur de $i_{P\cap K}^KE$ sur l'image de $J_{P\vert\ol{P}}(\sigma )$.
Par la relation (\#) de l'introduction, on trouve avec $\Stab (\Delta _{\so })=W_{\Delta _{\sso }}$ et les identifications de {\bf 2.1}
$$\varphi _{\xi }(\sigma )\circ\iota =
\sum _{w\in\sStab(\Delta _{\sso })}(J_{P\vert\ol{wP}}(\sigma )\lambda (w)\xi (w^{-1}\sigma )\lambda
(w^{-1})J_{wP\vert P}(\sigma ))\circ\iota .$$ Remarquons que, pour $w\in\Stab (\Delta _{\so })$,
l'op\'erateur $J_{P\vert\ol{wP}}(\sigma )\lambda (w)$ a la m\^eme image que $J_{P\vert\ol{P}}(\sigma
)$. On en d\'eduit, puisque l'image de $J_{P\vert\ol{P}}(\sigma )$ est stable par $\varphi _{\xi
}(\sigma )$ (cet endomorphisme venant d'un \'el\'ement de l'espace de Paley-Wiener), que
$$\eqalign { E_P^G(\varphi _{\xi }(\sigma )\circ\iota )&=\sum _{w\in\sStab (\Delta _{\sso })}
E_P^G((pJ_{P\vert w\ol{P}}(\sigma )\lambda (w)\otimes\iota ^{\vee }J_{P\vert wP}(\sigma ^{\vee
})\lambda (w))\xi (w^{-1}\sigma ))\cr &=\sum _{w\in\sStab (\Delta _{\sso
})}E_P^G(J_{P\vert\ol{P}}(w\sigma )\xi (w\sigma ))\cr }$$ Il r\'esulte alors de l'identit\'e (8.2.2)
avec les lemmes {\bf 4.8} et {\bf 4.10}, compte tenu de (8.1.2) et de l'invariance de $\mu $ par
$W(M,\o )$, que
$$\sum _{\pi\in\sE_2(G,\sigma )}\deg (\pi )E_P^G(\varphi _{\xi }(G,\pi )) =\gamma (G/M)\deg (\sigma
)\vert\Stab (A)\vert^{-1}E_P^G(\varphi _{\xi }(\sigma )\circ\iota )(\Res _A^P\mu )(\sigma ).$$ Par
suite, l'image de $\iota $ dans $i_{P\cap K}^KE$ correspond \`a une somme directe de repr\'esentations
de carr\'e int\'egrable (ce qui a d\'ej\`a \'et\'e montr\'e par Silberger (cf. [S3] Theorem 4.1.1)). Si
on note $m(\pi, \sigma )$ la multiplicit\'e d'une repr\'esentation irr\'eductible de carr\'e
int\'egrable $\pi $ dans l'image de $\iota $, on voit de plus que l'on a, par l'ind\'ependance
lin\'eaire des coefficients matriciels de repr\'esentations irr\'eductibles non \'equivalentes,
$$\deg (\pi ) =\gamma (G/M)\deg(\sigma )\vert\Stab (A)\vert ^{-1} m(\pi ,\sigma) (\Res _A^P\mu )(\sigma
).$$ En particulier, si $\Re (\sigma )>_P0$ (ce que Silberger appelle dans [S3] un point propre de
Casselman), il est facile de voir que l'image de $J_{P\vert\ol{P}}(\sigma )$ correspond \`a une seule
sous-repr\'esentation irr\'eductible $\pi $ de $i_P^G\sigma $ qui est de carr\'e int\'egrable par ce
qui pr\'ec\'edait. On trouve donc dans ce cas $\deg (\pi )=\gamma (G/M) \deg(\sigma )(\Res _A^P\mu
)(\sigma ).$

Si $\Re (\sigma )$ n'est pas dans $P_A$, on peut obtenir des expressions assez compliqu\'ees (voir par
exemple le cas de l'identit\'e (2)+(7)+(10) dans {\bf A.3}), dont il semble difficile de d\'eduire des
informations pr\' ecises sur le degr\'e formel des sous-quotients de carr\'e int\'egrable ou sur leur
position dans la repr\'esentation induite. Il semble toutefois \'eventuellement possible de d\'eduire
de  ces expressions une certaine propri\'et\'e d'invariance des degr\'es formels des repr\'esentations
de carr\'e int\'egrable ayant le m\^eme support cuspidal comme cela a \'et\'e fait dans [O2] pour des
identit\'es semblables dans le cas de la s\'erie principale non ramifi\'ee. Mentionnons \'egalement
[HO2] o\`u les degr\'es formels de certaines repr\'esentations de carr\'e int\'egrable apparaissant
dans la s\'erie principale non ramifi\'ee de $G$ sont calcul\'es.

\null {\bf 8.7} Le r\'esultat suivant est d\^u \`a E. Opdam: 

\null {\bf Th\'eor\`eme:} (cf. [O2] Theorem 3.29) \it Soit $\sigma\in \o$ un p\^ole de la fonction $\mu
$ de Harish-Chandra d'ordre \'egal au rang parabolique de $M$. Alors 
$\sum _{w\in W(M,\so )}\Res _{wA}^P\mu_{sp,\sigma }(w\sigma )\ne 0$.

\null Remarque: \rm  La convention 2.1 dans l'article [O2] qui n'a \'et\'e faite que pour des raisons de notations
et qui ne joue aucun r\^ole dans [O2] est d'ailleurs dans notre situation \'equiva-lente \`a la rationnalit\'e 
des points de r\'eductibilit\'e pour les induites paraboliques de repr\'esen-tations cuspidales. Elle a \'et\'e 
d\'emontr\'ee pour les groupes classiques d\'eploy\'es par C. Moeglin (cf. [M3]). Pour les repr\'esentations 
g\'en\'eriques, elle r\'esulte des travaux de Shahidi [Sh].

\null {\bf Corollaire:} \it Soit $\sigma\in \o$ tel que $\Re (\sigma )_G=0$. Pour que $\sigma $ appartienne au  
support cuspidal d'une repr\'esentation de carr\'e int\'egrable, il faut et il suffit que $\sigma $ soit un 
p\^ole de la fonction $\mu $ de Harish-Chandra d'ordre \'egal au rang parabolique de $M$. Ces p\^oles sont 
d'ordre maximal. \rm

\null
\null
\null
\vfill
\eject
{\bf Annexe A: L'exemple d'un groupe d\'eploy\'e semi-simple de type $G_2$}

\null
\null  Dans cet annexe, on va traiter l'exemple de la s\'erie principale non ramifi\'ee d'un groupe
semi-simple de type $G_2$ en d\'etail et "\`a la main". C'est le groupe de rang semi-simple minimal
pour lequel l'analyse des combinaisons lin\'eaires de coefficients matriciels qui apparaissent dans le
th\'eor\`eme {\bf 8.2} n'est pas sans difficult\'e. Ce probl\`eme appara\^\i t \'egalement pour des
groupes classiques (d\'eploy\'e ou non) de rang plus
\'elev\'e et d'autres groupes. Comme on l'a d\'ej\`a remarqu\'e, c'est l'analogue local des analyses
effectu\'ees dans l'annexe 3 de  [MW1]. Nous n'allons toutefois ici pas aussi loin dans le cas
difficile des identit\'es (2)+(7)+(10) (cf. {\bf A.3}) pour des raisons expliqu\'ees dans
l'introduction.

Remarquons finalement que le cas d'un groupe de type $A_n$ peut se traiter enti\`erement \`a la main,
comme cela a
\'et\'e fait pour la s\'erie principale non ramifi\'ee dans le manuscrit [H3], o\`u on utilise un
chemin d'int\'egration analogue
\`a celui effectu\'e dans [MW2]. D\'ej\`a pour la s\'erie principale non ramifi\'ee de Sp$_4$, nous
avons rencontr\'e des p\^oles qui s'annulent, ce qui paraissait compliqu\'e \`a montrer par une
m\'ethode directe. Disons tout de suite que ce dernier cas de figure n'appara\^\i t pas pour un groupe
de type $G_2$.

\null {\bf A.1.} Soit $G$ un tel groupe. Fixons un sous-groupe parabolique minimal $B=MU$ de $G$. Le
sous-groupe de Levi $M$ est un tore d\'eploy\'e maximal de $G$. Remarquons que tout sous-groupe de Levi
$M_{\gamma }$ de $G$ obtenu en associant \`a $M$ une racine $\gamma $ est isomorphe \`a $\GL _2$.  On
notera $\St _{2,\gamma }$ la repr\'esentation de Steinberg de $M_{\gamma }$. L'ensemble des racines
positives simples associ\'e \`a $B$ sera not\'e $\Delta =\{\alpha ,\beta \}$, $\alpha $ d\'esignant la
racine courte. L'ensemble des racines pour $G$ relatives \`a
$M$ qui sont  positives pour $B$, not\'e $\Sigma ^+$, est alors
$$\{\alpha ,\beta , \alpha +\beta ,2\alpha +\beta , 3\alpha +\beta ,3\alpha +2\beta \}.$$  Remarquons
que le syst\`eme des coracines $\Sigma ^{\vee }$ est \'egalement de type $G_2$. L'ensemble $\Sigma
^{\vee +}$ des coracines positives pour $B$ s'\'ecrit 
$$\{\alpha ^{\vee },\beta ^{\vee },\alpha ^{\vee }+\beta ^{\vee },\alpha ^{\vee }+2\beta ^{\vee
},\alpha ^{\vee }+3\beta ^{\vee }, 2\alpha ^{\vee }+3\beta ^{\vee }\},$$  en remarquant que $\alpha
^{\vee }$ est une racine longue.

Plus pr\'ecis\'ement, $(\alpha +\beta )^{\vee }=\alpha ^{\vee }+3\beta ^{\vee }$, $(2\alpha +\beta
)^{\vee }=2\alpha ^{\vee }+3\beta ^{\vee }$,
$(3\alpha +\beta )^{\vee }=\alpha ^{\vee }+\beta ^{\vee }$ et $(3\alpha +2\beta )^{\vee }=\alpha ^{\vee
}+2\beta ^{\vee }$, o\`u on a utilis\'e que $\langle \beta ,\alpha ^{\vee }\rangle =-3$ et $\langle
\alpha ,\beta ^{\vee }\rangle =-1$.

Le groupe de Weyl $W$ correspondant \`a $A_M$ est engendr\'e par les sym\'etries
$s_{\alpha }$ et $s_{\beta }$. C'est en fait le groupe di\'edral d'ordre 12 engendr\'e par $s_{\beta }$
et
$s_{\alpha }s_{\beta }$ avec $s_{\beta }^2=(s_{\alpha }s_{\beta })^6=1$. On a
$s_{\alpha }(\beta )=3\alpha +\beta $, $s_{\beta }(\alpha )=\alpha +\beta $.

Notons $\{\omega _{\alpha },\omega _{\beta }\}$ la base de $a_M^*$ qui est duale
\`a $\{\alpha ^{\vee },\beta ^{\vee }\}$. On a $\omega _{\alpha }=2\alpha +\beta
$, $\omega _{\beta }=3\alpha +2\beta $. Par suite,
$s_{\alpha }(\omega _{\alpha })=\omega _{\beta }-\omega _{\alpha }$,
$s_{\alpha }(\omega _{\beta })=\omega _{\beta }$, $s_{\beta }(\omega _{\alpha })=\omega _{\alpha }$ et
$s_{\beta }(\omega _{\beta })=3\omega _{\alpha }-\omega _{\beta }$.

On en d\'eduit que 
$$\displaylines {
\qquad\mu (\chi _{z_{\alpha }\omega _{\alpha }+z_{\beta }\omega _{\beta }}) ={(1-q^{z_{\alpha
}})(1-q^{-z_{\alpha }})
  (1-q^{z_{\beta }})(1-q^{-z_{\beta }})
  \over (1-q^{-1+z_{\alpha }})(-1-q^{-1-z_{\alpha }})
        (1-q^{-1+z_{\beta }})(-1-q^{-1-z_{\beta }})}
\qquad\cr
\hfill \times \ {(1-q^{z_{\alpha }+z_{\beta }})(1-q^{-z_{\alpha }-z_{\beta }})
        (1-q^{z_{\alpha }+2z_{\beta }})(1-q^{-z_{\alpha }-2z_{\beta }})
       \over (1-q^{-1+z_{\alpha }+z_{\beta }})(1-q^{-1-z_{\alpha }-z_{\beta }})   (1-q^{-1+z_{\alpha
}+2z_{\beta }})(1-q^{-1-z_{\alpha }-2z_{\beta }})} \qquad\cr
\hfill \times \ {(1-q^{z_{\alpha }+3z_{\beta }})
 (1-q^{-z_{\alpha }-3z_{\beta }})(1-q^{2z_{\alpha }+3z_{\beta }})
        (1-q^{-2z_{\alpha }-3z_{\beta }})
\over  (1-q^{-1+z_{\alpha }+3z_{\beta }})
        (1-q^{-1-z_{\alpha }-3z_{\beta }})
        (1-q^{-1+2z_{\alpha }+3z_{\beta }})
        (1-q^{-1-2z_{\alpha }-3z_{\beta }})}\qquad\cr  }$$  Les hyperplans affines radiciels dans $\o
=\X(M)$ qui sont singuliers pour $\mu $ sont donn\'es par $S_{\gamma }:=\o _{\chi _{\gamma/2},\gamma }$
avec $\gamma \in\Sigma ^+$. (On munit ces hyperplans donc de l'orientation donn\'ee par le choix de
$\Sigma ^+$.) L'ensemble $\S =\S _{\mu }$ est donc form\'e de ces hyperplans $S_{\gamma }$ ainsi que de
leurs parties imaginaires (\'egales \`a $\o _{1,\gamma }$).

Les hyperplans singuliers pour l'op\'erateur d'entrelacement $J_{B\vert\ol{B}}$ se d\'eterminent en le
d\'ecomposant en op\'erateurs d'entrelacement d\'efinis par rapport \`a des sous-groupes parabo-liques
adjacents. Soient $B'$ et
$B''$ deux sous-groupes de Borel adjacents dont l'intersection est un sous-groupe de Borel de
$M_{\gamma }$. Rappelons que $M_{\gamma }$ est isomorphe \`a $\GL _2$. Pour ce groupe, les op\'erateurs
d'entrelacement d\'efinis relatifs aux caract\`eres non ramifi\'es de $M$ sont bien connus. On en
d\'eduit que $J_{B'\vert B''}$ est r\'egulier en $\chi _{\lambda }$, si et seulement si
$q^{-\langle\lambda ,\gamma ^{\vee }\rangle}\ne 1$. L'op\'erateur d\'efini en $\chi _{\lambda }$ est
alors alors bijectif, si et seulement si $q^{-\langle\lambda ,\gamma ^{\vee }\rangle}\not\in
\{q^{-1},q\}$. Remarquons qu'un op\'erateur d'entrelacement qui est le compos\'e de plusieurs
op\'erateurs d'entrelacement dont certains ont des p\^oles en $\chi _{\lambda }$ et d'autres ne sont
pas injectifs en
$\chi _{\lambda }$ peut bien \^etre d\'efini en $\chi _{\lambda }$.

\null {\bf A.2} Notons $E$ l'espace de la repr\'esentation unit\'e de $M$. Pour calculer l'int\'egrale
$$\int _{\Re(\chi )=r\gg _B0} E_B^G(J_{B\vert\ol{B}},\xi ,\chi ) \mu (\chi )d\Im(\chi )\eqno {\rm
(A.1.1)}$$  pour une application polynomiale
$\xi:\o\rightarrow i_{\ol{B}\cap K}^KE\otimes i_{B\cap K}^KE^{\vee }$ donn\'ee, on va effectuer le
chemin $\C$ dans
$a_M^*$ ci-dessous d\'ej\`a d\'ecrit dans [MW2], appendice 3. On identifiera dans la suite parfois la
valeur de
$\xi $ en un caract\`ere non ramifi\'e $\chi $ \`a un
\'el\'ement de $\Hom (i_{B\cap K}^KE, i_{\ol{B}\cap K}^KE)$.

\vfill
\eject

\vskip 8.5cm A cette occasion, on trouve des p\^oles, lorsque l'on coupe les hyperplans $H_{\gamma ,1}$
de $a_M^*$ (qui sont donn\'es par $\langle
\lambda ,\gamma ^{\vee }\rangle =1$, $\gamma\in\Sigma ^+$). L'origine de ces hyperplans est le point
${\gamma\over 2}$. Soit $\delta $ une racine dans $\Sigma $, tel que $\{\gamma ,\delta \}$ soit un
ensemble de racines simples pour un certain ordre sur $a_M^*$. Soit $\{\omega _{\gamma },\omega
_{\delta }\}$ la base duale \`a la base $\{\gamma ^{\vee },\delta ^{\vee }\}$ de $a_M$. Tout point de
$H_{\gamma ,1}$ s'\'ecrit alors sous la forme ${\gamma\over 2}+t_{\delta }\omega _{\delta }$. \'Ecrits
sous cette forme, les points d'intersection du chemin $\C $ avec les \'el\'ements de $\H(\S )$ sont:
$$\eqalign {
\gamma =\alpha ,\ \delta =\beta :\qquad &{\alpha\over 2}+t_{\delta }\omega _{\beta },\qquad (t_{\delta
}\gg 0);\cr
\gamma =\beta ,\ \delta =\alpha :\qquad &{\beta\over 2}+t _{\delta }\omega _{\alpha },\qquad ({3\over
2}< t_{\delta }<{5\over 2});\cr
\gamma = 3\alpha +\beta ,\ \delta =-(2\alpha +\beta ):\qquad &{3\alpha +\beta\over 2}+(-t _{\delta
})(-\alpha -\beta ),\qquad (-{3\over 2}< t _{\delta }< -{1\over 2});\cr
\gamma =3\alpha +2\beta ,\ \delta =-(\alpha +\beta ):\qquad &{3\alpha +2\beta\over 2}+ t_{\delta
}\alpha ,\qquad (0<t_{\delta }<{1\over 2});\cr
\gamma =\alpha +\beta ,\ \delta =-(3\alpha +2\beta ):\qquad &{\alpha +\beta\over 2}+(-t_{\delta
})(-3\alpha -\beta ),\qquad (-{1\over 6}<t_{\delta }<0)\cr
\gamma =2\alpha +\beta ,\ \delta =-(3\alpha +\beta ):\qquad &{2\alpha +\beta\over 2}+t_{\delta }\beta
,\qquad (-{1\over 2}<t_{\delta }<{1\over 2}).\cr }$$

\null  Posons $r_{\delta }=t_{\delta }\omega _{\delta }$. Comme $H_{\gamma ,1}$ est la partie r\'eelle
d'un unique hyperplan singulier $S_{\gamma }$, le calcul de l'int\'egrale (A.1.1) le long du chemin $\C
$ donne
$$\displaylines {
\ \ \ \ \int _{S_{\alpha ,r_{\beta }}} E_B^G(J_{B\vert\ol{B}},\xi ,\chi ) (\Res_{S _{\alpha }}\mu
)(\chi )\ d\Im(\chi )\hfill\cr 
\cr +\int _{S_{\beta,r_{\alpha }}} E_B^G(J_{B\vert\ol{B}},\xi ,\chi )  (\Res_{S _{\beta }}\mu )(\chi )\
d\Im(\chi )\hfill\cr
\cr +\int _{S_{3\alpha +\beta , r_{-2\alpha -\beta }}} E_B^G(J_{B\vert\ol{B}},\xi ,\chi )
(\Res_{S_{3\alpha +\beta }}\mu )(\chi )\ d\Im(\chi )\hfill\cr
\cr +\int _{S _{3\alpha +2\beta ,r_{-\alpha -\beta }}} E_B^G(J_{B\vert\ol{B}},\xi ,\chi
)(\Res_{S_{3\alpha +2\beta }}\mu )(\chi )\ d\Im(\chi )\hfill\cr 
\cr +\int _{S _{\alpha +\beta ,r_{-3\alpha -2\beta }}} E_B^G(J_{B\vert\ol{B}},\xi ,\chi )(\Res_{S
_{\alpha +\beta }}\mu )(\chi )\ d\Im(\chi )\hfill\cr
\cr +\int _{S _{2\alpha +\beta ,r_{-(3\alpha +\beta )}}} E_B^G(J_{B\vert\ol{B}},\xi ,\chi )(\Res_{S
_{2\alpha +\beta }}\mu )(\chi )\ d\Im(\chi )\hfill\cr
\cr+\int _{\Re(\chi )=0} E_B^G(J_{B\vert\ol{B}},\xi ,\chi ) \mu (\chi )\ d\Im(\chi )\hfill\cr }$$

Il reste maintenant \`a rejoindre l'axe unitaire, i.e. \`a effectuer \`a chaque int\'egrale un
d\'ecalage vers
$r_{\delta }=0$. Pour d\'eterminer les p\^oles qui apparaissent, posons $\ti {B}_{\alpha }=(M_{\alpha
}\cap B)\ol {U}_{\alpha }$ et $\ti {B}_{\beta }=(M_{\beta }\cap B)\ol{U}_{\beta }$. On a
$$\displaylines {
\qquad\mu _{B\vert\ti{B}_{\alpha }}(\chi _{{\alpha \over 2}+z_{\beta }\omega _{\beta
}})={(1-q^{{1\over 2}+z_{\beta }})(1-q^{{1\over 2}-z _{\beta }}) (1-q^{2z_{\beta }})(1-q^{-2z_{\beta
}}) 
\over (1-q^{-{3\over 2}+z_{\beta }})(-1-q^{-{3\over 2}-z_{\beta }}) (1-q^{-1+2z_{\beta
}})(-1-q^{1-2z_{\beta }})} \hfill\cr
\hfill \times  {(1-q^{{1\over 2}+3z_{\beta }})(1-q^{{1\over 2}-3z_{\beta }})
\over (1-q^{-{3\over 2}+3z_{\beta }})(1-q^{-{3\over 2}-3z_{\beta }})}\qquad\cr }$$  et
$$\mu _{B\vert\ti {B}_{\beta }}(\chi _{z_{\alpha }\omega _{\alpha }+{\beta\over 2}}) ={(1-q^{{3\over
2}-z_{\alpha }})(1-q^{{3\over 2}+z_{\alpha }}) (1-q^{2z_{\alpha }})(1-q^{-2z_{\alpha }})
\over (1-q^{-{5\over 2}-z_{\alpha }})(1-q^{-{5\over 2}+z_{\alpha }}) (1-q^{-1+2z_{\alpha
}})(1-q^{-1-2z_{\alpha }}). }$$

Ces formules restent valables, si on remplace $\alpha $ par une racine courte ou
$\beta $ par une racine longue, tout en remplacant $\omega _{\beta }$ ou $\omega _{\alpha }$ de fa\c
con appropri\'ee, les constantes restant les m\^emes.

On remarque qu'avec $B_{\gamma }$ un sous-groupe de Borel pour lequel $\gamma $ est simple positif et
$\ti {B}_{\gamma }=(M_{\gamma }\cap B)\ol {U}_{\gamma }$, on a         $$(\Res _{S_{\gamma }}\mu )=\mu
_{B_{\gamma }\vert\ti {B}_{\gamma }}(\Res _{S_{\gamma }}\mu ^{M_{\gamma }}),$$  puisque les p\^oles de
$\mu ^{M_{\gamma }}$ sont simples. Les p\^oles de $\Res _{S_{\gamma }}\mu $ sont donc donn\'es par ceux
de $\mu _{B\vert
\ti {B}_{\gamma }}$.

Effectuons maintenant \`a chaque int\'egrale ci-dessus le changement de contours de $r_{\delta }$ \`a
$0$. Le calcul de $\mu _{B_{\gamma }\vert\ti{B}_{\gamma }}$ ci-dessus montre qu'aucun p\^ole de
$J_{B\vert\ol {B}}$ n'appara\^\i t lors de ce proc\'ed\'e qui n'est pas \'egalement un p\^ole de la
fonction $\mu _{B_{\gamma }\vert\ti{B}_{\gamma }}$ (compte tenu de ce qui a \'et\'e dit sur les
hyperplans singuliers de
$J_{B\vert\ol{B}}$ \`a la fin de {\bf A.1}).  Notons $\vert\omega _{\delta }\vert$ l'unique \'el\'ement
de $\Sigma ^+\cap \{\omega _{\delta },-\omega _{\delta }\}$, fixons $z_0\in \CC$ et notons $z_0'$
l'unique nombre complexe qui v\'erifie $z_0\omega _{\delta }=z_0'\vert \omega _{\delta }\vert$. Alors
on d\'efinit $\Res ^+_{{\gamma\over 2}+z_0\omega _{\delta }}$ comme \'etant l'op\'erateur $\R
(S_{\gamma },\S (S_{\gamma }))\rightarrow\CC $ qui associe
\`a $\psi $ la valeur $\Res _{z=z_0'}\psi (\chi _{{\gamma\over 2}+z \vert\omega _{\delta }\vert})$.
(Remarquons qu'en notant $\delta _1$ la racine dans $\Sigma (A_{M_{\gamma }})$ qui est la restriction
\`a $A_{M_{\gamma }}$ de l'unique racine positive dans
$\{\delta ,-\delta \}$, l'\'el\'ement de $a_{M_{\gamma }}^*$ que l'on a not\'e $\ti{\delta _1}$ en {\bf
3.2} est \'egal \`a
$\vert\omega _{\delta }\vert
$.) 

Pour l'int\'egrale (A.1.1) on obtient alors 
$$\displaylinesno {
\ \ \ \Res ^+_{{\alpha\over 2}+{3\over 2}\omega _{\beta }} (E_B^G(J_{B\vert\ol{B}},\xi ,.)\ (\Res_{S
_{\alpha }}\mu ))\hfill &(1) \cr 
\cr + \Res ^+_{{\alpha\over 2}+{1\over 2}\omega _{\beta }} (E_B^G(J_{B\vert\ol{B}},\xi ,.) (\Res_{S
_{\alpha }}\mu ))\hfill &(2)\cr 
\cr + \Res ^+_{{\alpha\over 2}+({1\over 2}+{\pi i\over\log q})\omega _{\beta }}
(E_B^G(J_{B\vert\ol{B}},\xi ,.)(\Res_{S_{\alpha }}\mu ))\hfill &(3)\cr
\cr + \Res ^+_{{\alpha\over 2}+({1\over 2}+{2\pi i\over3\log q})\omega _{\beta
}}(E_B^G(J_{B\vert\ol{B}},\xi ,.)(\Res_{S_{\alpha }}\mu ))\hfill &(4)\cr
\cr + \Res ^+_{{\alpha\over 2}+({1\over 2}+{4\pi i\over3\log q})\omega _{\beta
}}(E_B^G(J_{B\vert\ol{B}},\xi ,.)(\Res_{S_{\alpha }}\mu ))\hfill &(5)\cr
\cr + \int _{S_{\alpha ,0}} E_B^G(J_{B\vert\ol{B}},\xi ,\chi ) (\Res_{S _{\alpha }}\mu )(\chi )
d\Im(\chi )\hfill &(6)\cr 
\cr + \Res ^+_{{\beta\over 2}+{1\over 2}\omega _{\alpha }} (E_B^G(J_{B\vert\ol{B}},\xi
,.)(\Res_{S_{\beta }}\mu ))\hfill &(7)\cr
\cr + \Res ^+_{{\beta\over 2}+({1\over 2}+{\pi i\over \log q})\omega _{\alpha
}}(E_B^G(J_{B\vert\ol{B}},\xi ,.)(\Res_{S_{\beta }}\mu ))\hfill &(8)\cr
\cr +\int _{S_{\beta ,0}} E_B^G(J_{B\vert\ol{B}},\xi ,\chi ) (\Res_{S _{\beta }}\mu )(\chi ) d\Im(\chi
)\hfill &(9)\cr
\cr +\Res ^+_{{3\alpha +\beta\over 2}+{\alpha +\beta\over 2}} (E_B^G(J_{B\vert\ol{B}},\xi
,.)(\Res_{S_{3\alpha +\beta }}\mu ))\hfill &(10)\cr
\cr +\Res ^+_{{3\alpha +\beta\over 2}+({1\over 2}+{\pi i\over\log q})(\alpha +\beta
)}(E_B^G(J_{B\vert\ol{B}},\xi ,.)(\Res_{S_{3\alpha +\beta }}\mu )) \hfill &(11)\cr 
\cr +\int _{S_{3\alpha +\beta, 0}} E_B^G(J_{B\vert\ol{B}},\xi, \chi ) (\Res_{S_{3\alpha +\beta }}\mu
)(\chi )d\Im(\chi )\hfill &(12)\cr 
\cr +\int _{S_{3\alpha +2\beta ,0}}E_B^G(J_{B\vert\ol{B}},\xi ,\chi ) (\Res_{S_{3\alpha +2\beta }}\mu
)(\chi )d\Im(\chi )\hfill &(13)\cr +\int _{S_{\alpha +\beta ,0}}E_B^G(J_{B\vert\ol{B}},\xi ,\chi ) 
(\Res _{S_{\alpha +\beta }}\mu )(\chi )d\Im(\chi )\hfill &(14)\cr 
\cr +\int _{S_{2\alpha +\beta ,0}}E_B^G(J_{B\vert\ol{B}},\xi ,\chi ) (\Res_{S_{2\alpha +\beta }}\mu
)(\chi )d\Im(\chi )\hfill &(15)\cr  +\int _{\Re(\chi )=0} E_B^G(J_{B\vert\ol{B}},\xi ,\chi ) \mu (\chi
) d\Im(\chi )\hfill &(16)\cr  }$$

\null {\bf A.3} Il reste \`a rassembler les diff\'erents termes obtenus. Remarquons que, si $\iota $
est l'injection d'une sous-repr\'esentation de $i_B^G\chi $ dans $i_{B\cap K}^KE$, alors l'image de
$\varphi _{\xi }(\chi )\circ\iota
$ est incluse dans celle de $\iota $. Soit $p$ un projecteur de $i_{P\cap K}^KE$ sur l'image de $\iota
$. On va utiliser \`a plusieurs reprises le fait suivant:

(A.3.1) si $J$ est un \'el\'ement de $i_{B\cap K}^KE\otimes i_{B\cap K}^KE^{\vee }$ tel que
l'endomorphisme de
$i_{B\cap K}^KE$ correspondant ait l'image incluse dans celle de $p$, alors $E_{B,\chi }^G
((p\circ\iota ^{\vee })J)=E_{B,\chi } ^G(J)$.

\null
$\underline {(16):}$ Comme par un changement de variable,
$$(16) = \int _{\Re (\chi )=0} E_B^G(J_{B\vert\ol{B}},\xi ,w^{-1}\chi )\mu (\chi ) d\Im(\chi ),$$ 
pour tout $w\in W$, il r\'esulte de la relation (\#) de l'introduction que 
$$(16) = {1\over \vert W\vert}\int _{\Re(\chi )=0} E_B^G(\varphi _{\xi }(\chi )) \mu (\chi ) d\Im(\chi
).$$

\null
$\underline {(6)+(14)+(15):}$ Notons $\iota $ l'injection de $i_{P_{\alpha }}^G(\St _{2,\alpha
}\otimes\chi _{z_{\beta }\omega _{\beta }})$ dans
$i_B^G (\chi _{{\alpha\over 2}+z_{\beta }\omega _{\beta }})$ et $p$ la projection de $i_{B\cap K}^KE $
sur l'image de
$\iota $ de noyau \'egal \`a
$Im(J_{B\vert s_{\alpha }B}(\chi _{-{\alpha\over 2}+z_{\beta }\omega _{\beta }}))$. Si $w\in W$
v\'erifie
$l(s_{\alpha }w)<l(w)$, alors $J_{wB\vert B}(\chi _{{\alpha\over 2}+z_{\beta }\omega _{\beta }})\circ
\iota =0$. On en d\'eduit pour $z_{\beta }$ dans un ouvert dense de $\CC $
$$\displaylines { 
\ \ \varphi (\chi _{{\alpha\over 2}+z_{\beta }\omega _{\beta }})\circ\iota  \hfill\cr  =\sum
_{\scriptstyle w\in W\atop\scriptstyle l(s_{\alpha }w)>l(w)}   p J_{B\vert\ol {wB}}(\chi
_{{\alpha\over 2}+z_{\beta }\omega _{\beta }})\lambda (w)\xi (w^{-1}\chi _{{\alpha\over 2}+z_{\beta
}\omega _{\beta }}))\lambda (w^{-1})J_{wB\vert B}(\chi _{{\alpha\over 2}+z_{\beta }\omega _{\beta
}}^{\vee })\ \iota\hfill\cr }$$ Comme $Im(J_{B\vert
\ol{wB}}(\chi _{{\alpha\over 2}+z_{\beta }\omega _{\beta }}))\subseteq Im(p)$, si $l(s_{\alpha
}w)>l(w)$, on en d\'eduit par (A.3.1)
$$E_B^G(\varphi (\chi _{{\alpha\over 2}+z_{\beta }\omega _{\beta }})\circ\iota )=\sum _{\scriptstyle
w\in W\atop\scriptstyle l(s_{\alpha }w)>l(w)} E_B^G(J_{B\vert\ol{B}}(w^{-1}\chi _{{\alpha\over
2}+z_{\beta }\omega _{\beta }})\xi (w^{-1}\chi _{{\alpha\over 2}+z_{\beta }\omega _{\beta }}))$$  En
calculant les $w^{-1}\chi _{{\alpha\over 2}+z_{\beta }\omega _{\beta }}$, on trouve alors
$$(6)+(14)+(15)={1\over 2}\int _{S_{\alpha ,0}} E_B^G(\varphi _{\xi }(\chi )\circ \iota )
(Res_{S_{\alpha }}\mu )(\chi ) dIm(\chi ).$$

\null
$\underline {(9)+(12)+(13):}$ C'est analogue au cas pr\'ec\'edent, en
\'echangeant $\alpha $ et $\beta $.

\null
$\underline {(1):}$ L'op\'erateur d'entrelacement $J_{B\vert\ol {B}}(\chi _{{\alpha\over 2}+z_{\beta
}\omega _{\beta }})$ est r\'egulier en
$z_{\beta }={3\over 2}$. On a ${\alpha\over 2}+{3\over 2}\omega _{\beta }=\omega _{\alpha }+\omega
_{\beta }$. Soit $\iota $ l'injection de l'image de
$J_{B\vert\ol {B}}(\chi _{\omega _{\alpha }+\omega _{\beta }})$. Il est imm\'ediat par la relation (\#
) de l'introduction que
$$\varphi _{\xi }(\chi _{\omega _{\alpha }+\omega _{\beta }})\circ\iota =pJ_{B\vert\ol{B}}(\chi
_{\omega _{\alpha }+\omega _{\beta }})\xi (\chi _{\omega _{\alpha }+\omega _{\beta }})\iota,$$   d'o\`u
l'on d\'eduit par (A.3.1) que
$$(1)=E_B^G(\varphi _{\xi }(\chi _{\omega _{\alpha }+\omega _{\beta }}\circ\iota ))\Res
^+_{{\alpha\over 2}+{3\over 2}\omega _{\beta }}(\Res _{S_{\alpha }}(\mu )).$$ La repr\'esentation
$i_B^G\chi _{\omega _{\alpha }+\omega _{\beta }}$ poss\`ede (par exemple par la remarque 1) de {\bf
8.7.1}) un unique sous-quotient irr\'eductible de carr\'e int\'egrable $\pi (\chi_{\omega _{\alpha
}+\omega _{\beta }})$ qui correspond \`a l'image de l'op\'e-rateur d'entrelacement $J_{B\vert\ol
{B}}(\chi _{\omega _{\alpha }+\omega _{\beta }})$. Par l'identit\'e (8.6.1),  on trouve $\deg(\pi
(\chi_{\omega _{\alpha }+\omega _{\beta }}))=\gamma(G/M)$ $\Res ^+_{{\alpha\over 2}+{3\over 2}\omega
_{\beta }}(\Res _{S_{\alpha }}(\mu )).$

\null
$\underline {(4)+(5):}$ L'op\'erateur d'entrelacement $J_{B\vert\ol {B}}(\chi _{{\alpha\over
2}+z_{\beta }\omega _{\beta }})$ est r\'egulier en
$z_{\beta }={1\over 2}+{2\pi i\over 3\log q}$. On a 
${\alpha\over 2}+({1\over 2}+{2\pi i\over 3\log q})\omega _{\beta }=\omega _{\alpha }+{2\pi i\over
3\log q}\omega _{\beta }$. Soit $\iota $ l'injection de $Im(J_{B\vert\ol {B}}(\chi _{\omega _{\alpha
}+{2\pi i\over 3\log q}\omega _{\beta }}))$ dans $i_{B\cap K}^KE$ et $p$ la projection de $i_{B\cap
K}^KE$ sur l'image de $\iota $ de noyau \'egal
\`a celui de $J_{B\vert\ol {B}}(\chi _{\omega _{\alpha }+{2\pi i\over 3\log q}\omega _{\beta }})$. On
v\'erifie que
$J_{wB\vert B}(\chi _{\omega _{\alpha }+{2\pi i\over 3\log q}\omega _{\beta }})\circ \iota =0$ sauf si
$w\in \{1, s_{\beta }\}$. Par suite,
$$\displaylines { 
\ \ \varphi_{\xi } (\chi _{\omega _{\alpha }+{2\pi i\over 3\log q}\omega _{\beta }})\circ
\iota\hfill\cr  =pJ_{B\vert\ol {B}}(\chi _{\omega _{\alpha }+{2\pi i\over 3\log q}\omega _{\beta }})\xi
(\chi _{\omega _{\alpha }+{2\pi i\over 3\log q}\omega _{\beta }})\iota\hfill\cr  +
pJ_{B\vert\ol{s_{\beta }B}}(\chi _{\omega _{\alpha }+{2\pi i\over 3\log q}\omega _{\beta }})\lambda
(s_{\beta })\xi (s_{\beta }\chi _{\omega _{\alpha }+{2\pi i\over 3\log q}\omega _{\beta }})\lambda
(s_{\beta })J_{s_{\beta }B\vert B} (\chi _{\omega _{\alpha }+{2\pi i\over 3\log q}\omega _{\beta
}})\iota.\hfill\cr  }$$

\null  Remarquons que l'op\'erateur $J_{\ol{s_{\beta }B}\vert\ol{B}}(\chi _{\omega _{\alpha }+{2\pi
i\over 3\log q}\omega _{\beta }})$ est bijectif. L'image de $J_{B\vert\ol{s_{\beta }B}}$ $(\chi
_{\omega _{\alpha }+{2\pi i\over 3\log q}\omega _{\beta }})$ est donc incluse dans
$Im(J_{B\vert\ol {B}}(\chi _{\omega _{\alpha }+{2\pi i\over 3\log q}\omega _{\beta }}))$. On en
d\'eduit par (A.3.1) que
$$\displaylines { 
\ \ E_B^G(\varphi (\chi _{\omega _{\alpha }+{2\pi i\over 3\log q}\omega _{\beta }})\circ
\iota)\hfill\cr  =E_B^G(J_{B\vert\ol {B}}(\chi _{\omega _{\alpha }+{2\pi i\over 3\log q}\omega _{\beta
}})\xi (\chi _{\omega _{\alpha }+{2\pi i\over 3\log q}\omega _{\beta }}))+E_B^G((J_{B\vert\ol
{B}}(s_{\beta }\chi _{\omega _{\alpha }+{2\pi i\over 3\log q}\omega _{\beta }})\xi (s_{\beta }\chi
_{\omega _{\alpha }+{2\pi i\over 3\log q}\omega _{\beta }})),\hfill\cr  }$$  d'o\`u
$$(4)+(5)=E_B^G(\varphi _{\xi }(\chi _{\omega _{\alpha }+{2\pi i\over 3\log q}\omega _{\beta
}})\circ\iota) (\Res ^+_{{\alpha\over 2}+({1\over 2}+{2\pi i\over\log q})\omega _{\beta }}(\Res _{S
_{\alpha }}\mu )),$$  les deux r\'esidus qui apparaissent dans (4) et (5) \'etant \'egaux.

On peut montrer facilement (comparer [Mu] proposition 4.2) que la repr\'esentation $i_B^G\chi _{\omega
_{\alpha }+{2\pi i\over 3\log q}\omega_{\beta }}$ poss\`ede un unique sous-quotient irr\'eductible de
carr\'e int\'egrable que l'on notera $\pi (\chi _{\omega _{\alpha }+{2\pi i\over 3\log q}\omega_{\beta
}})  $ et que celui-ci correspond \`a l'image de $J_{B\vert\ol {B}}(\chi _{\omega _{\alpha }+{2\pi
i\over 3\log q}\omega _{\beta }}))$. Par l'identit\'e (8.6.1), on en d\'eduit que $\deg(\pi (\chi
_{\omega _{\alpha }+{2\pi i\over 3\log q}\omega_{\beta }}))=\gamma(G/M)(\Res ^+_{{\alpha\over
2}+({1\over 2}+{2\pi i\over\log q})\omega _{\beta }}\ $ $(\Res _{S _{\alpha }}\mu ))$.

\null
$\underline {(3)+(8)+(11):}$ L'op\'erateur d'entrelacement $J_{B\vert\ol {B}}(\chi _{{\alpha\over
2}+z_{\beta }\omega _{\beta }})$ est r\'egulier en $z_{\beta }={1\over 2}+{\pi i\over\log q}$. On a
${\alpha\over 2}+({1\over 2}+{\pi i\over\log q})\omega _{\beta }=\omega _{\alpha }+{\pi i\over\log
q}\omega _{\beta }$. Soit $\iota $ l'injection de
$Im(J_{B\vert\ol {B}}(\chi _{\omega _{\alpha }+{\pi i\over\log q}\omega _{\beta }}))$ dans $i_{B\cap
K}^KE$ et $p$ la projection de $i_{B\cap K}^KE$ sur l'image de $\iota $ de noyau \'egal \`a celui de
$J_{B\vert\ol {B}}(\chi _{\omega _{\alpha }+{\pi i\over\log q}\omega _{\beta }})$. On v\'erifie que
$J_{wB\vert B}(\chi _{\omega _{\alpha }+ {\pi i\over\log q}\omega _{\beta }})\circ\iota =0,$ sauf si
$w\in \{1, s_{\beta }, s_{\beta }s_{\alpha }\}$. Par suite,
$$\displaylines { 
\varphi (\chi _{\omega _{\alpha }+{\pi i\over\log q}\omega _{\beta }})\circ \iota\hfill\cr 
=pJ_{B\vert\ol {B}}(\chi _{\omega _{\alpha }+{\pi i\over\log q}\omega _{\beta }})\xi (\chi _{\omega
_{\alpha }+{\pi i\over\log q}\omega _{\beta }})\iota\hfill\cr  +pJ_{B\vert \ol{s_{\beta }B}}(\chi
_{\omega _{\alpha }+{\pi i\over\log q}\omega _{\beta }})\lambda (s_{\beta })\xi (s_{\beta }\chi
_{\omega _{\alpha }+{\pi i\over\log q}\omega _{\beta }})\lambda (s_{\beta })J_{s_{\beta }B\vert B}(\chi
_{\omega _{\alpha }+{\pi i\over\log q}\omega _{\beta }}))\iota\hfill\cr  +pJ_{B\vert\ol{s_{\beta
}s_{\alpha }B}}(\chi _{\omega _{\alpha }+{\pi i\over\log q}\omega _{\beta }})\lambda (s_{\beta
}s_{\alpha })\xi (s_{\alpha }s_{\beta }\chi _{\omega _{\alpha }+{\pi i\over\log q}\omega _{\beta
}})\lambda(s_{\alpha }s_{\beta })J_{s_{\beta }s_{\alpha }B\vert B}(\chi _{\omega _{\alpha }+{\pi
i\over\log q}\omega _{\beta }}))\iota\hfill\cr  }$$ Remarquons que les op\'erateurs $J_{\ol{s_{\beta
}B}\vert\ol{B}}(\chi _{\omega _{\alpha }+{\pi i\over\log q}\omega _{\beta }})$ et $J_{\ol{s_{\beta
}s_{\alpha }B}\vert\ol{B}}(\chi _{\omega _{\alpha }+{\pi i\over\log q}\omega _{\beta }})$ sont
inversi-bles. Par suite, les images de $J_{B\vert \ol{s_{\beta }B}}(\chi _{\omega _{\alpha }+{\pi
i\over\log q}\omega _{\beta }})$ et de $J_{B\vert\ol{s_{\beta }s_{\alpha }B}}(\chi _{\omega _{\alpha
}+{\pi i\over\log q}\omega _{\beta }})$ sont \'egales \`a celle de $J_{B\vert\ol {B}}(\chi _{\omega
_{\alpha }+{\pi i\over\log q}\omega _{\beta }})$. On en d\'eduit avec (A.3.1) que 
$$\displaylines { 
\ \ E_B^G(\varphi (\chi _{\omega _{\alpha }+{\pi i\over\log q}\omega _{\beta }})\circ\iota)\hfill\cr 
=E_B^G(J_{B\vert\ol {B}}(\chi_{\omega _{\alpha }+{\pi i\over\log q}\omega _{\beta }} )\xi(\chi _{\omega
_{\alpha }+{\pi i\over\log q}\omega _{\beta }}))\hfill\cr  +E_B^G(J_{B\vert\ol {B}}(s_{\beta }\chi
_{\omega _{\alpha }+{\pi i\over\log q}\omega _{\beta }})\xi (s_{\beta }\chi _{\omega _{\alpha }+{\pi
i\over\log q}\omega _{\beta }}))\hfill\cr +E_B^G(J_{B\vert\ol {B}}(s_{\alpha }s_{\beta }\chi _{\omega
_{\alpha }+{\pi i\over\log q}\omega _{\beta }})\xi (s_{\alpha }s_{\beta }\chi _{\omega _{\alpha }+{\pi
i\over\log q}\omega _{\beta }})),\hfill\cr  }$$  d'o\`u l'on d\'eduit
$$(3)+(8)+(11)=E_B^G(\varphi _{\xi }(\chi _{\omega _{\alpha }+{\pi i\over\log q}\omega _{\beta }})\circ
\iota ) (\Res ^+_{{\alpha\over 2}+({1\over 2}+{\pi i\over\log q})\omega _{\beta }} (\Res_{S_{\alpha
}}\mu )),$$  les r\'esidus de
$\mu $ qui apparaissent dans (3), (8) et (11) \'etant
\'egaux.

On peut montrer facilement (comparer [Mu] proposition 4.1) que la repr\'esentation $i_B^G\chi _{\omega
_{\alpha }+{\pi i\over\log q}\omega_{\beta }}$ poss\`ede un unique sous-quotient irr\'eductible de
carr\'e int\'egrable que l'on notera $\pi (\chi _{\omega _{\alpha }+{\pi i\over \log q}\omega_{\beta
}})$ et que celui-ci correspond \`a l'image de $J_{B\vert\ol {B}}(\chi _{\omega _{\alpha }+{\pi i\over
\log q}\omega _{\beta }}))$. Par l'identit\'e (8.6.1) on en d\'eduit alors que $\deg (\pi (\chi
_{\omega _{\alpha }+{\pi i\over \log q}\omega_{\beta }}))=\gamma (G/M)(\Res ^+_{{\alpha\over
2}+({1\over 2}+{\pi i\over\log q})\omega _{\beta }}$ $(\Res_{S_{\alpha }}\mu ))$. 

\null
$\underline {(2)+(7)+(10):}$ Ce cas est bien plus compliqu\'e, puisque l'op\'erateur d'entrelacement
$J_{B\vert\ol{B}}$ $(\chi _{{\alpha\over 2}+z _{\beta } \omega _{\beta }})$ n'est pas r\'egulier en
$z_{\beta }={1\over 2}$. Remarquons que ${\alpha\over 2}+{1\over 2}\omega _{\beta }=\omega _{\alpha }$

Il est connu (cf. [Mu]) que la repr\'esentation $i_B^G\chi  _{\omega _{\alpha }}$ poss\`ede exactement
deux sous-quotients qui sont de carr\'e int\'egrable. Ils sont non isomorphes. Notons $St _{2,\alpha
}$  et $St _{2,\beta }$ respectivement les repr\'esentations de Steinberg de $M_{\alpha }$ et de
$M_{\beta }$. Alors la sous-repr\'esenta-tion semi-simple  maximale de $i_{P_{\alpha }}^G(St _{2,\alpha
}\chi _{{1\over 2}\omega _{\beta }})$ est de longueur $1$ et de carr\'e int\'egrable. Notons-la $\pi
_1(\chi _{\omega_{\alpha }})$. La  sous-repr\'esentation semi-simple maximale de $i_{P_{\beta }}^G(St
_{2,\beta }\chi _{{1\over 2}\omega _{\alpha }})$ est somme directe de deux repr\'esentations
irr\'eductibles de carr\'e int\'egrables dont l'une est isomorphe \`a $\pi _1(\chi _{\omega_{\alpha
}})$. D\'esignons l'autre par $\pi _2(\chi_{\omega_{\alpha }})$. Le degr\'e formel de ces deux
repr\'esentations est calcul\'e dans [R] (avec une normalisation des mesures diff\'erente de la
n\^otre). On a $\deg(\pi _1(\chi _{\omega_{\alpha }}))={1\over 2}\deg (\pi _2(\chi _{\omega_{\alpha
}}))$. 

On peut retrouver ces r\'esultats par calcul direct en prenant pour $\xi $ l'application polynomiale
constante de valeur $\lambda (w_0)1_{\sI w_0\sI} \otimes 1_{\sI}$, $\I $ d\'esignant le sous-groupe
d'Iwahori qui correspond \`a $B$. Pour d\'eterminer les termes (2), (7) et (10), on peut alors utiliser
la description de Bernstein et Rogawski des op\'erateurs d'entrelacement d\'efinis relatifs \`a
l'alg\`ebre d'Iwahori-Hecke, comme cela est expliqu\'e dans l'annexe 3 de [MW1]. Tout est alors
explicite et on retrouve les r\'esultats ci-dessus par calcul (\'eventuellement \`a l'aide d'un
logiciel).

Il serait probablement possible de simplifier les calculs et de rendre le tout un peu plus
intelligible, en employant une m\'ethode semblable \`a celle utilis\'ee dans l'annexe 3 [MW1].

Nous renoncons ici \`a aller plus loin dans l'analyse de ce cas, d'une part puisque les r\'esultats
sont d\'ej\`a connus, et d'autre part puisque l'on n'a pas trouv\'e de m\'ethode vraiment intelligente
pour traiter la somme des trois termes qui serait par ailleurs susceptible de se g\'en\'eraliser.

\null
\null
\null
{\bf Annexe B: Preuve du th\'eor\`eme de Opdam} 
 
\null 
\null 
{\bf B.1} Pour la commodit\'e du lecteur, nous donnons ci-apr\`es avec l'autorisation
de l'auteur une preuve du th\'eor\`eme {\bf 8.7} qui est - comme d\'ej\`a signal\'e dans 
le texte, d\^u \`a E. Opdam (cf. [O2] theorem 3.29). La preuve ci-dessous est essentiellement 
la sienne. La motivation pour l'ins\'erer ici vient du fait que la preuve de ce r\'esultat dans 
l'article [O2] est assez concise et r\'edig\'ee dans un langage diff\'erent du n\^otre.  
  
Remarquons que, comme nous travaillons dans le cadre des groupes $p$-adiques, 
notre r\'esultat concerne la fonction $\mu $ de Harish-Chandra, alors 
que Opdam s'int\'eresse aux fonctions $\mu $ venant des alg\`ebres d'Iwahori-Hecke. 
 
\null 
{\bf B.2} Par le th\'eor\`eme des r\'esidus {\bf 3.10}, la proposition 
{\bf 4.11} et les notations de {\bf 8.1} et {\bf 8.2}, il existe une unique 
famille de donn\'ees de r\'esidus $\Res _A^P$, $A\in \A(\S )$, v\'erifiant 
$\Res _{wA}(\psi\circ w^{-1})=(\Res _A\psi )\circ w^{-1}$ pour tout $w\in 
W(M,\o)$, $w\Sigma _{\so }^{M_A}(P\cap M_A)\subseteq\Sigma (P)$, (cf. 
lemmes {\bf 4.8}   et {\bf 4.11}), telles que, pour toute fonction 
rationnelle $\psi $ dans $\R (\o ,\S )$, on ait  
$$\eqalign { 
&\int _{\Re (\sigma )=r\gg _P0} \psi (\sigma ) d\Im (\sigma ) \cr 
=&\sum _{\Omega\subseteq \Delta _{\sso}} \sum _{A\in [\sA(\sS)], M_A=M_{\Omega}} 
\vert W_{\Delta _{\sso },M_A}\vert ^{-1} \vert \P_{\sS}(M_A)\vert ^{-1} 
\vert \Stab (A)\vert ^{-1}\cr 
&\int _{\Re (\sigma )=r(A)+\epsilon _A}\sum _{w'\in W^{M_A}(M,\so )} 
\sum _{w\in W_{M_A}^+(M,\so )}(\Res _{w'A}^P\psi (w.))(w'\sigma )d_A\Im 
(\sigma ).\cr  }$$ 
(cf. {\bf 8.2.3}). 
 
La premi\`ere partie de la preuve du th\'eor\`eme {\bf 3.10} donne par 
ailleurs un proc\'ed\'e r\'ecursif de calcul de ces donn\'ees de r\'esidu: 
appelons point singulier d'un sous-espace affine $L\in\L(\S )$ tout point 
de $L-L_{\sreg}$. Fixons $r\gg _P0$. Choisissons pour tout 
$Q\in\P_{\sS}(M)$ un chemin de $r$ \`a $r(\o )+\epsilon _Q$ qui ne contient 
qu'un nombre fini de points singuliers de $a_M^*$, chacun se trouvant sur un 
unique hyperplan affine $H\in\H(\S )$. (Rappelons que $r(\o )=0$.)  Chaque  
point singulier travers\'e donne lieu \`a des op\'erateurs r\'esidus  
$\pm \Res _S$, correspondant aux $S\in \S$ dont la partie r\'eelle 
contient ce point.  
 
On continue alors ce proc\'ed\'e avec chaque point 
singulier rencontr\'e sur les chemins ci-dessus relatif \`a l'hyperplan 
affine $H\in\H $ qui le contient, en rempla\c cant, pour tout $S\in\S $ de 
partie r\'eelle $H$, $r(\o )$, $M$ et $\P_{\sS}(M)$ respectivement par 
$r(H)$, $M_H$ et  $\P_{\sS}(M_H)$. On obtient alors de nouveaux op\'erateurs 
de r\'esidu relatifs \`a des hyperplans affines $A$ de $S$. On compose 
ceux-ci avec $\Res _S$ pour obtenir des op\'erateurs r\'esidus $Res _A$ 
relatifs \`a $\o $. On continue ce proc\'ed\'e jusqu'\`a ce que l'on 
aboutisse \`a des sous-espaces affines v\'erifiant $A=A_{\sreg}$. (C'est 
certainement v\'erifi\'e si $A$ est un singleton.) 
 
La donn\'ee de r\'esidu $\Res _A^P$ en un \'el\'ement $A\in\A(\S )$ se d\'eduit 
alors, apr\`es normalisation par la constante $\vert\P_{\sS }(M_A)\vert ^{-1}$ 
(cf. {\bf 3.10}), de la somme des op\'erateurs r\'esidus en $A$ qui apparaissent 
lors du proc\'ed\'e ci-dessus. 
 
Remarquons que, pour tout $A\in\A (\S )$, toute donn\'ee de r\'esidus $Res _A$ 
en $A$,  tout $\psi\in\R (\o,\S )$ et tout $\sigma\in A_{\sreg }$, la fonction 
$\lambda\mapsto (Res _A\psi )  (\sigma\otimes\chi _{\lambda })$ est 
r\'eguli\`ere sur $a_G^*$. Pour le calcul des donn\'ees de r\'esidu $\Res 
_A^P$, seule importe donc la projection du chemin sur $a_G^*$.  
 
\null  
{\bf B.3} L'objet de cette annexe est la preuve du th\'eor\`eme suivant:  
  
\null  
{\bf Th\'eor\`eme:} \it Soit $A\in [\A _{\mu }(\S )]$ tel que $M_A=G$. Alors   
$$\sum _{w\in W(M,\so )} (\Res _{wA}^P\mu )(w\sigma )\ne 0$$  
pour tout $\sigma \in A$. \rm  
 
\null 
\it Remarque: \rm Le r\'esultat figurant dans le th\'eor\`eme {\bf 8.7} 
concerne la fonction $\mu _{sp, \sigma }$, $\sigma \in A$, et non pas la 
fonction $\mu $. Les deux assertions sont \'equivalentes gr\^ace \`a (8.1.2),  
puisque $\mu _{nsp ,\sigma }$ est r\'eguli\`ere sur $A$ et invariant par  
$W(M,\o )$. La formulation du th\'eor\`eme {\bf 8.7} avait \'et\'e motiv\'e 
par le fait que le th\'eor\`eme 3.29 de l'article [O2], o\`u on travaille 
avec des fonctions $\mu $ venant des alg\`ebres d'Iwahori-Hecke, s'y applique 
directement.    
 
\null  {\bf B.4} On r\'eduira d'abord la preuve du th\'eor\`eme \`a celle de 
la  proposition suivante:     
 
\null  {\bf Proposition:} \it Il existe un sous-groupe parabolique $P'$ de 
sous-groupe  de Levi $M$ et $\tau \in A$ tel que $(\Res _A^{P'}\mu )(\tau  
)\ne 0$.  \rm     
 
\null   
{\bf B.5} Il faut donc prouver la proposition suivante: 
 
\null 
{\bf Proposition:} \it Supposons la proposition {\bf B.4} v\'erifi\'ee. Alors
le  th\'eor\`eme {\bf B.3} est vraie. \rm 
 
\null 
Comme les chambres dans $a_M^*$ relatives \`a $\H (\S )$ sont conjugu\'ees 
par les \'el\'ements de $W(M,\o )$, il existe $w\in W(M,\o )$ tel que 
$\Res _A^{w^{-1}P}=\Res _A^{P'}$. Par {\bf 4.8}, on a donc $0\ne (\Res _A^{P'} 
\mu )(\sigma )=(\Res _A^{w^{-1}P}\mu )(\tau )=(\Res _{wA}^P\mu )(w\tau 
)$.   Le lemme suivant est \'el\'ementaire: 
 
\null 
{\bf Lemme:} \it Il existe une fonction polynomiale $\psi $ sur $\o $  
qui ne s'annulle pas en $w\tau $ et qui v\'erifie $\psi (w'\tau )=0$ pour 
tout  $w\in W(M,\o )$ tel que $w'\tau \ne w\tau $. \rm 
 
\null 
Soit $\psi $ une fonction sur $\o $ v\'erifiant le lemme et appliquons 
le th\'eor\`eme {\bf 8.2} \`a $\xi = \psi \xi _{P,\so }^H$ et $g=1$. Le terme 
de gauche de l'\'egalit\'e (8.2.2) de ce th\'eor\`eme vaut alors $\sum  
_{w'\in W(M,\so )}$ $(\Res _{w'A}^P(\psi\mu ))(w'\tau )$. 
Comme   $(\Res _{w'A}^P(\psi\mu ))(w'\tau )=\psi (w'\tau )(\Res _{w'A}^P\mu 
)$ $(w'\tau )$ par {\bf 8.1.2}, cette expression est non nulle par choix 
de $\psi $. L'ensemble $\E _2(G,\tau )$ est donc non vide, d'o\`u $\sum 
_{w'\in W(M,\so )}(\Res _{w'A}^P\mu )  (w'\tau )\ne 0$ par {\bf 8.6}. \hfill 
{\fin 2}     
 
\null 
{\bf B.6} La preuve de la proposition {\bf B.4} se fait en plusieurs \'etapes 
par r\'ecurrence sur le rang parabolique de $M$. Si $\dim (a_M^{G*})=1$, on se 
trouve dans le cadre des fonctions complexes d'une variable o\`u la proposition  
est \'el\'ementaire.  
 
Supposons donc $\dim (a_M^{G*})>1$. Comme par la remarque \`a la fin de {\bf B.2}, 
seule la projection des chemins sur $a_M^{G*}$ importe, il suffit 
de consid\'erer le cas o\`u $G$ est semi-simple (et donc $a_M^{G*}=a_M^*$). 
Notons $\S _{\tau }$ l'ensemble des hyperplans affines de $\S $ passant par 
$\tau $, $\H _{\tau }$ l'ensemble des parties r\'eelles des \'el\'ements de 
$\S _{\tau }$ et $S$ la sph\`ere dans $a_M^{*}$ de centre $r_{\tau }:=\Re 
(\tau )$ passant par l'origine de $a_M^*$. (Le symb\^ole $S$ ne
d\'esignera donc plus un sous-espace affine de $\o $.) Toute droite $L$ qui
appartient \`a $\L(\H _{\tau})$ est partag\'ee par $r_{\tau }$ en deux
demi-droites. Le lemme suivant (qui est analogue \`a la remarque {\bf 3.14}
de [Ho2]) montre qu'une et une seule des deux demi-droites contient $r(L)$.

\null
{\bf Lemme:} \it On a $r(L)\ne r_{\tau }$.

\null
Preuve: \rm Soit $A\in \A(\S _{\tau })$ avec $L=\Re (A)$. Notons $\mu
_{M_{A}}$ la fonction rationnelle \'egale \`a $\mu (\mu^{M_{A}})^{-1}$.
Alors, par les r\'esultats du paragraphe {\bf 8}, l'hypoth\`ese de
r\'ecurrence, et l'holomorphie de la fonction $\mu $ de Harish-Chandra en une
repr\'esentation de carr\'e int\'egrable, la fonction $\sum _{w\in
W^{M_{A}}(M,\so )} (\Res _{wA}\mu )(w\sigma )=\mu _{M_{A}}(\sigma )\sum _{w\in
W^{M_{A}}(M,\so )} (\Res _{wA}\mu ^{M_{A}})(w\sigma )$ est r\'egulier pour $\Re
(\sigma )=r(L)$, d'o\`u $r(L)\ne r_{\tau }$. \hfill {\fin 2}

\null
Dessinons en blanc la demi-droite qui contient $r(L)$ et l'autre en
noir.
    
Choisissons un point $r_1$ sur la sph\`ere $S$ qui se trouve sur une demi-droite  
noire et qui est de distance \`a $0$ minimale avec cette propri\'et\'e. Notons  
$B$ la boule dans $a_M^{*}$ de centre $0$ passant par $r_1$.  
   
La droite $l$ passant par $r_1$ et $r_{\tau }$ peut s'\'ecrire sous la forme  
$r_1+a_{M_1}^{*}$ avec $M_1$ un sous-groupe de Levi maximal de $G$. Notons $A_1$ 
l'\'el\'ement de $\A (\S _{\tau })$ de partie r\'eelle $r_1+a_{M_1}^{*}$. (Il 
est d\'etermin\'e de fa\c con unique, puisque ses \'el\'ements sont n\'ecessairement 
de la forme $\tau\otimes\chi _{\lambda }$, $\lambda\in a_{M_1}^{*}$.) Par hypoth\`ese 
de r\'ecurrence, il existe un sous-groupe parabolique $P'$ de $G$ tel que, 
l'on ait $(\Res _{A_1}^{M_1\cap P'}\mu ^{M_1})(\sigma )\ne 0$ pour tout 
$\sigma \in A_1$. On peut choisir $P'$ de telle sorte que $P'M_1$ soit un
sous-groupe parabolique de Levi $M_1$ et que $r_1-r_{\tau }>_{P'M_1}0$.
    
Remarquons que l'op\'erateur $\Res _A^{P'}$ est d\'etermin\'e, gr\^ace \`a 
{\bf 3.9}, par sa restriction \`a l'espace $\R (\o ,\S _{\tau  })$. En 
particulier, les hyperplans de la forme $\o _{\Im (\sigma ),\alpha }$, 
$\o  _{\sigma ,\alpha  }\in \S_{\tau }$ (cf. {\bf 3.10}) peuvent \^etre 
n\'eglig\'es lors du calcul des r\'esidus. (Ils apparaissent toutefois ci-dessous
dans la d\'efinition des ensembles $\P _{\sS }(M_L)$ et des points $\epsilon
_Q$ (cf. {\bf 3.10}).)
           
\null 
{\bf B.7} Avant de d\'ecrire un proc\'ed\'e qui calcule les donn\'ees de
r\'esidus relatives \`a $\R (\o ,\S _{\tau })$ et \`a $P'$ et qui montre que
$(\Res _A^{P'}\mu )$ $(\tau  )\ne 0$, on va expliciter dans cette section
certains points de $a_M^*$. La proposition simple suivante nous servira \`a
plusieures reprises.

\null
{\bf Proposition:} \it Soit $L\in \L(\S _{\tau })$ et $p\in L_{\sreg }$. La
demi-droite partant de $r_{\tau }$ et passant par $p$ (et priv\'ee du point
$r_{\tau }$) est contenue dans une et une seule composante connexe de $L_{\sreg }$.

\null
Preuve: \rm Sinon, cette demi-droite contiendrait un point singulier $p'$ qui
serait donc contenu dans un hyperplan $H\in \H(\S _{\tau })$ qui ne contient
pas $L$. Mais alors, comme $r_{\tau }\in H$, la droite passant par $p'$ et
$r_{\tau }$ serait incluse dans $H$. Ainsi $p$ serait singulier, d'o\`u une
contradiction. \hfill {\fin 2}  

\null
Pour tout $L\in \L(\S _{\tau })$, notons $\ti {r(L)}$ le point d'intersection 
de $S$ avec la demi-droite partant de $r_{\tau }$ et passant par $r(L)$. On  
v\'erifie facilement que $\ti {r(L)}$ est le point de $L\cap S$ \`a distance 
minimale de $0$. Si $\dim (L)>1$, $L\cap S$ est connexe, et on appelle 
$\ti{r(L)}$ l'origine sph\'erique de $L$ ou de $L\cap S$. Si $L$ est une
droite, alors $L\cap S$ est form\'e de deux points que l'on appelle des
origines sph\'eriques. Le point $\ti {r(L)}$ est blanc et l'autre noir
d'apr\`es la convention faite au num\'ero pr\'ec\'edent.  

Pour tout $Q\in\P_{\sS}(M_L)$, soit $\ti {r(L)}+\ti{\epsilon _Q}$ le point  
d'intersection de $S$ avec la demi-droite partant de $r_{\tau }$ et passant 
par $r(L)+\epsilon _Q$. La proposition ci-dessus montre que, $\ti{r(L)} +\ti
{\epsilon _Q}$ et  $r(L)+\epsilon _Q$ sont dans une m\^eme composante connexe
de $L_{\sreg }$.  
 
Pour tout $L\in\L(\S _{\tau })$ d'origine sph\'erique dans l'int\'erieur de 
$B$ et tout $Q\in\P_{\sS}(M_L)$, on peut supposer $\epsilon _Q$ choisi 
suffisamment pr\`es de $0$, pour que  $\ti{r(L)}+\ti{\epsilon _Q}$ soit \`a 
l'int\'erieur de $B$.   
 
On observe alors que, si $p$ est un point \`a l'int\'erieur du segment 
de $r_{\tau }$ \`a $r_1$, l'intersection de $S$ avec l'hyperplan affine 
$p+a_M^{M_1*}$ est une sph\`ere $S_p^{M_1}$ de centre $p$. On peut choisir $p$ 
de telle sorte que les hyperplans affines $H\in\H(\S _{\tau })$ qui ont une
intersection $\ne\{r_{\tau }\}$ avec le c\^one ferm\'e de centre $r_{\tau
}$ contiennent tous $l$ (o\`u, rappelons-le, $l$ d\'esigne la droite passant
par $r$ et $r_{\tau }$).  

Supposons $p$ choisi avec cette propri\'et\'e. Remarquons que $r_1+a_M^{M_1*}$
 est le plan tangent \`a $S$ en $r_1$. Notons $S_1^{M_1}$ la projection  
orthogonale de $S_p^{M_1}$ sur $r_1+a_M^{M_1*}$ et $0_1^{M_1}$ celle de $0$.
Pour  $L\in \L(\S_{\tau })$, $L\supseteq l$, soit $r(L)_1^{M_1}$ la projection
orthogonale de $r(L)$ et, pour $Q\in \P_{\sS}(M_L)$, $r(L)_1^{M_1}+
\epsilon_{Q,1}^{M_1}$ celle de $r(L)+\epsilon _Q$. Il en suit que
$r(L)_1^{M_1}$ est le point de  $L\cap (r_1+a_M^{M_1*})$ de distance minimale
\`a $0_1^{M_1}$. Par ailleurs

\null
{\bf Lemme:} \it Le point $r(L)_1^{M_1}+\epsilon _{Q,1}^{M_1}$ n'est
contenu dans aucun hyperplan $H\in\H (\S _{\tau })$ qui contient $l$ et non
pas $L$.
 
\null
Preuve: \rm Supposons que $H$ soit un tel hyperplan. Alors $H$ contiendrait
\'egalement la droite par $r(L)_1^{M_1}+\epsilon  _{Q,1}^{M_1}$ 
parall\`ele \`a $l$, et, par suite, $r(L)+\epsilon _Q\in H\cap L_{\sreg }$, ce
qui est impossible. \hfill {\fin 2}

\null
{\bf B.8} Fixons un point $0'$ de la demi-droite allant de $r_1$ \`a
$0_1^{M_1}$ et qui se   trouve \`a l'int\'erieur de la sph\`ere
$S_1^{M_1}$. Consid\'erons la projection centrale de centre $r_1$ et de
rapport  $\lambda =d(r_1,0')$ $d(r_1,0_1^{M_1})^{-1}$. Notons, pour 
$L\in \L(\S _{\tau })$, $L\supseteq l$, et $Q\in \P _{\sS}(M_L)$, $r(L)'$ et
$r(L)'+\epsilon  _Q'$ respectivement les images de $r(L)_1^{M_1}$ et de
$r(L)_1^{M_1}+ \epsilon_{Q,1}^{M_1}$ par cette projection. Remarquons que
$r(a_M^*)'=0'$ et que $r(l)'=r_1$. 

\null
{\bf Lemme:} \it Le point $r(L)'$ est le point de $L\cap (r_1+a_M^{M_1})$ \`a
distance minimale de $0'$. Il appartient \`a l'int\'erieur de $S_1^{M_1}$ et \`a
l'int\'erieur de $B$. Quitte \`a remplacer \'eventuellement $\epsilon _Q$ par 
un point de norme absolue plus petite, $r(L)'+\epsilon _Q'$ se trouve \`a 
l'int\'erieur de $S_1^{M_1}$ ainsi qu'\`a l'int\'erieur de $B$. Il appartient
alors \'egalement \`a $L_{\sreg }$.

\null
Preuve: \rm La premi\`ere assertion est une cons\'equence imm\'ediate du fait
qu'une projection centrale laisse invariant les arcs et multiplie les distances 
par son rapport. Comme $L\cap (r_1+a_M^{M_1*})$ contient $r_1$ et que $r(L)'$ est
donc plus pr\`es de $0'$ que $r_1$, $r(L)'$ se trouve 
par le th\'eor\`eme de Pythagore \`a l'int\'erieur de $S_1^{M_1}$. En appliquant 
le th\'eor\`eme de Pythagore au triangle d'extr\'emit\'es $r_1$, $r(L)_1^{M_1}$ et 
$0_1^{M_1}$, on en d\'eduit \'egalement que $r(L)'$ est plus proche de $0_1^{M_1}$ 
(et donc de $0$) que $r_1$, si $L\ne l$. Par suite $r(L)'$ se trouve \`a l'int\'erieur
de $B$. Il est alors clair que l'on peut \'egalement supposer $r(L)'+\epsilon _Q'$ 
\`a l'int\'erieur de $S_1^{M_1}$ ainsi qu'\`a l'int\'erieur de $B$, quitte \`a 
remplacer $\epsilon _Q$ par un point de norme absolue plus petite. Par choix de 
$S_1^{M_1}$, ce point n'appartient donc \`a aucun hyperplan affine $H\in\H (\S  )$ 
qui ne contient pas $l$. Suite \`a la proposition et le lemme {\bf B.7}, il ne se 
trouve pas non plus sur un hyperplan contenant $l$ et non pas $L$. Il est donc 
r\'egulier.\hfill {\fin 2}

\null 
{\bf B.9} {\bf Lemme:} \it Soit $r$ un point r\'egulier de $S_1^{M_1}$
tel que la projection orthogonale de $r-r_1$ sur $a_M^{M_1*}$ soit dans la
chambre de Weyl positive de $P'\cap M_1$. Alors, $r$ se trouve dans la m\^eme
composante connexe de $a_{M,reg}^*$ que tout \'el\'ement de $a_M^*$ qui est
suffisamment positive dans la chambre de Weyl de $P'$.

\null
Preuve: \rm Consid\'erons la demi-droite partant de $r_{\tau }$ et passant par
$r$. Les points $\ne r_{\tau }$ sont tous r\'egulier gr\^ace \`a la proposition 
{\bf B.7}, puisque $r$ l'est, et ils correspondent aux valeurs de l'application 
$r(t)=r+t \overrightarrow{r_{\tau }r}$, $t>0$. Il suffit donc de montrer que
$\langle r(t),\alpha ^{\vee }\rangle $ est  strictement croissante en $t$ pour
toute racine $\alpha $ qui est ou dans   $\Sigma (P'\cap M_1)$ ou dans $\Sigma
(P'M_1)$.  

Si $\alpha\in\Sigma (P'\cap M_1)$, $\langle r-r_{\tau }, \alpha ^{\vee
}\rangle =\langle r-r_1, \alpha ^{\vee }\rangle $ et donc $\langle r(t),
\alpha ^{\vee }\rangle =\langle r,\alpha ^{\vee}\rangle +t\langle r-r_1,\alpha
^{\vee }\rangle$. Comme $\langle r-r_1,\alpha ^{\vee }\rangle >0$ par choix de 
$P'$, la fonction $r(t)$ est donc strictement croissante dans ce cas. Si $\alpha\in\Sigma
(P'M_1)$, $\langle r-r_{\tau }, \alpha ^{\vee }\rangle =\langle r_1-r_{\tau },
\alpha ^{\vee }\rangle $ et donc $\langle r(t), \alpha ^{\vee }\rangle =\langle
r,\alpha ^{\vee}\rangle +t\langle r_1-r_{\tau } ,\alpha ^{\vee }\rangle$. Comme
$\langle r_1-r_{\tau } ,\alpha ^{\vee }\rangle >0$ par choix de $P'$, on conclut 
de m\^eme que la fonction $r(t)$ est strictement croissante. 

\hfill{\fin 2}

\null
{\bf B.10} Choisissons $r\in a_M^*$ v\'erifiant les hypoth\`eses du lemme {\bf
B.9}. On d\'eduit de ce lemme que l'on peut choisir $r$ comme point de
d\'epart pour calculer les donn\'ees de r\'esidus relatives \`a $P'$ et $\S 
_{\tau }$.

Rappelons (cf. {\bf B.6}) que l'on a not\'e $A_1$ l'unique \'el\'ement de $\A
(\S _{\tau  })$ de partie r\'eelle $l$ et que $\L (\S _{A_1})$ est alors
\'egal \`a l'ensemble des $L\in \L(\S )$ qui contiennent $l$. Observons que, si
$L\in \L (\S _{A_1})$, tout \'el\'ement $Q\in \P _{\sS _{\tau }}(M_L)$ est un
sous-ensemble d'un \'el\'ement $Q'\in \P_{\sS _{A_1}}(M_L)$. Choisissons pour
tout $Q'\in \P_{\sS _{A_1}}(M_L)$ un $Q\in \P _{\sS _{\tau }}(M_L)$
v\'erifiant $Q\subseteq Q'$, et posons $\epsilon _{Q'}:=\epsilon _Q$.  (On
\'ecrira alors en particulier $\epsilon _{Q'}':=\epsilon _Q'$.)

Commencons par effectuer le proc\'ed\'e d\'ecrit en {\bf B.2}, partant de
$r$, mais relatif au choix de $0'$ comme origine et relatif aux points
$r(L)'+\epsilon  _{Q'}'$, $L\in\L (\S _{A_1})$, $Q'\in\P  _{\sS _{A_1}}(M_L)$,
en se rappelant que $r(L)'$ est le point de $L$ \`a distance minimale de $0'$.
(On commence donc par relier $r$ aux points $0'+\epsilon  _{Q'}'$, $Q'\in\P
_{\sS _{A_1}}(M)$, et ainsi de suite. Tous les chemins peuvent donc \^etre
pris \`a l'int\'erieur de $S_1^{M_1}$, en sorte que l'on ne rencontre en effet
que des sous-espaces affines $L\in\L(\S _{A_1})$, et ainsi de suite.)    
Le lemme suivant montre que suite \`a ce proc\'ed\'e la donn\'ee de
r\'esidus  $\Res _{A_1}^{P'}$ a d\'ej\`a \'et\'e calcul\'ee:
      
\null 
{\bf Lemme:} \it Il existe un chemin de $r(L)'+\epsilon _{Q'}'$ \`a
$r(L)+\epsilon _{Q'}$ qui ne rencontre  aucun point singulier de $L$ relatif
\`a $\H (\S _{A_1})$.

\null 
Preuve: \rm Supposons d'abord $L\ne l$. Notons $\ti{r(L)'}+\ti{\epsilon _Q'}$
la projection  de $r(L)'+\epsilon _Q'$ sur $S$ parall\`ele \`a la droite
$l$. Ce point est donc plus proche de $0$ que $r(L)'+\epsilon _Q'$, et il
appartient donc \`a $B$ gr\^ace au lemme {\bf B.8}. Relions $r(L)'+\epsilon _Q'$ 
\`a $\ti {r(L)'}+\ti{\epsilon _Q'}$ par la  droite passant par ces deux points. 
Par choix de $S_1^{M_1}$, ce chemin ne contient aucun point singulier de $L$ 
relatif \`a $\S_{A_1}$ (et d'ailleurs pas non plus relatif \`a $\S _{\tau }$). 
Par ailleurs, il reste  \`a l'int\'erieur de $B$, puisque tout point sur ce
chemin est plus pr\`es de $0$ que $r(L)'+\epsilon _Q'$.      

On remarque alors que le point $\ti {r(L)}+\ti {\epsilon _Q}$ est \`a 
l'int\'erieur de $B$, puisque $\ti{r(L)'}+\ti{\epsilon _Q'}$ y est et que 
$\ti {r(L)}+\ti {\epsilon _Q}$ est (quitte \`a changer \'eventuellement 
$\epsilon _Q$ par un point de norme absolue plus petite) n\'ecessairement plus 
proche de $0$.

L'hyperplan affine qui contient $r_1$ et qui est orthogonale \`a la droite passant par $0$
et $r_{\tau }$ partage l'espace $a_M^*$ en deux r\'egions. On observe que
l'intersection de $S$ avec la r\'egion qui contient $0$ est pr\'ecis\'ement
$B\cap S$. Une des deux arcs sur $S$ reliant $\ti{r(L)'}+ \ti{\epsilon _Q'}$
\`a $\ti{r(L)}+ \ti{\epsilon _Q}$ est donc contenue dans $B$.

Par ailleurs, cet arc ne rencontre aucun hyperplan affine $H\in\H(\S_{A_1})$
ne contenant pas $L$: sinon, il en serait de m\^eme de sa projection
orthogonale  sur  $r_1+a_M^{M_1*}$. Rappelons que $\ti {r(L)}+\ti {\epsilon
_{Q'}}$  est sur  la demi-droite partant de $r_{\tau }$ et passant par
$r(L)+\epsilon  _{Q'}$.  La projection orthogonale de $r(L)+\epsilon _{Q'}$
sur $r_1+a_M^{M_1*}$  \'etant $r(L)_1^{M_1}+\epsilon _{Q',1}^{M_1}$ et celle de
$r_{\tau }$ \'etant  $r_1$, l'image de $\ti {r(L)}+\ti {\epsilon _{Q'}}$ par
cette projection  orthogonale se trouve sur la droite passant par $r_1$ et
$r(L)_1^{M_1}+  \epsilon _{Q',1}^{M_1}$. Cette droite contient \'egalement le
point $r(L)'+\epsilon _{Q'}'$ qui est la projection orthogonale de $\ti
{r(L)'}+\ti{\epsilon  '_{Q'}}$. Il en suit que la projection orthogonal de
l'arc est un segment de la demi-droite partant de $r_1$  et passant par
$r(L)_1^{M_1}+\epsilon _{Q,1}^{M_1}$. Ce segment ne contient pas $r_1$. Comme
$r(L)_1^{M_1}+\epsilon _{Q,1}^{M_1}$ est r\'egulier relatif \`a $\S _{A_1}$ 
par le lemme {\bf B.7}, suite \`a la  proposition {\bf B.7}, aucun point de 
ce segment ne peut appartenir \`a un hyperplan affine  $H\in \H(\S _{A_1})$ 
ne contenant pas $L$.        

On relie alors le point $\ti {r(L)}+\ti {\epsilon _{Q'}}$ \`a $r(L)+ \epsilon
_{Q'}$ par la droite passant par ces points. La proposition {\bf B.7} montre
que celle-ci ne contient aucun point singulier relatif \`a $\S _{A_1}$. 
On fait de m\^eme pour $L=l$. Ceci prouve le lemme. \hfill{\fin 2} 

\null
{\bf B.11} {\bf Lemme:} \it Soient $L\in\L (\S _{A_1})$, $Q'\in\P _{\sS
_{A_1}}(M_L)$ et $Q\in\P _{\sS }(M_L)$ tel que $Q\subseteq Q'$. Alors
$r(L)'+\epsilon _Q'$ et $r(L)'+\epsilon _{Q'}'$ sont dans une m\^eme
composante connexe de $L_{\sreg}$.

\null
Preuve: \rm Par d\'efinition de $\P _{\sS _{A_1}}(M_L)$ et de $\epsilon _{Q'}$,
$r(L)+\epsilon _Q$ et $r(L)+\epsilon _{Q'}$ sont dans une m\^eme composante
connexe de $L$ relatif \`a $\H (\S _{A_1})$. Cette composante connexe contient
par un raisonnement analogue au lemme {\bf B.7} \'egalement 
$r(L)^{M_1}_1+\epsilon _{Q, 1}^{M_1}$, et par suite elle contient les projections
centrales $r(L)'+\epsilon _Q'$ et $r(L)'+\epsilon _{Q'}'$ de ces points. Ils
peuvent donc \^etre reli\'es par un chemin \`a l'int\'erieur de $S_1^{M_1}$ qui
ne peut donc, par choix de $S_1^{M_1}$, recontrer aucun hyperplan 
$H\in\H (\S _{\tau })$ qui contient $L$ et qui ne contient pas $l$.\hfill {\fin 2}

\null   
{\bf B.12} Effectuons maintenant un proc\'ed\'e suivant {\bf B.2}
pour calculer les donn\'ees de r\'esidus relatives \`a $\S _{\tau }$. Il faut
donc prendre en compte tous les hyperplans affines de $\S _{\tau  }$ et
consid\'erer, pour $L\in \L(\S _{A_1})$ les \'el\'ements de $\P _{\sS _{\tau
}}(M_L)$ et pas seulement ceux de $\P _{\sS _{A_1 }}(M_L)$. Soit d'abord $L\in
\L(\S _{A_1})$, $L\ne l$, et $Q\in \P _{\sS _{\tau }}(M_L)$. Il existe alors $Q'\in \P
_{\sS _{A_1 }}(M_L)$ avec $Q\subseteq Q'$. Comme $r(L)'+\epsilon _{Q}'$ et
$r(L)'+\epsilon  _{Q'}'$ se trouvent par {\bf B.11} dans la m\^eme composante 
connexe de $L_{\sreg }$, on peut les relier sans rencontrer de point singulier. 
Ensuite, on relie $r(L)'+\epsilon  _{Q}'$ \`a $r(L)+\epsilon _{Q}$ par le 
m\^eme proc\'ed\'e que $r(L)'+\epsilon  _{Q'}$ \`a $r(L)+\epsilon _{Q}$. (En
effet, $\epsilon _{Q'}$ et $\epsilon _Q$ sont interchangeable.) Ce chemin reste
\`a l'int\'erieur de $B$.

Ce chemin ne rencontre seulement un nombre fini de points singuliers: s'il y
en a, ils se trouvent, par ce qui a \'et\'e dit dans la preuve du lemme 
{\bf B.10}, sur l'arc de
$\ti{r(L)'}+ \ti{\epsilon _Q'}$ \`a $\ti {r(L)}+\ti {\epsilon _Q}$. S'il
existait une infinit\'e de points, au moins deux d'entre eux se trouveraient
sur un m\^eme hyperplan affine. La droite passant par ces deux points serait
donc contenue dans cet hyperplan et donc, comme cet  hyperplan contient
$r_{\tau }$, \'egalement la droite passant par $\ti{r(L)'}+\ti{\epsilon _Q'}$
et $\ti {r(L)}+\ti {\epsilon _Q}$, ce  qui contredit le fait que $\ti
{r(L)}+\ti {\epsilon _Q}$ est r\'egulier. Quitte \`a d\'eformer le chemin dans
un tr\`es petit voisinage, on peut \'egalement supposer que tout point
singulier de ce chemin se trouve sur un et un seul hyperplan affine singulier
de $L$, sans que les autres propri\'et\'es de ce chemin n'aient chang\'e.

On continue alors le proc\'ed\'e avec chacun des points singuliers relatif \`a 
l'hyperplan affine singulier de $L$ qui le contient. Comme ces points se
trouvent \`a l'int\'erieur de $B$, l'origine sph\'erique de cet hyperplan
affine de $L$ se trouve \'egalement \`a l'int\'erieur de $B$. On peut donc
prendre les chemins \`a l'int\'erieur de $B$. Et ainsi de suite.

On aboutit alors \`a des points $\ti {r(L)}+\ti {\epsilon _Q}$, $L\in \L(\S
)$, $L\ne l$, $Q\in \P_{\sS }(M_L)$, \`a l'int\'erieur de $B$ (sans rencontrer 
des origines sph\'eriques noires d'une droite $L$), auxquels il faut ajouter le 
point $r_1$.

Par nos choix et par la proposition {\bf B.7}, on peut relier les points $\ti
{r(L)}+\ti {\epsilon _Q}$, $L\ne l$, \`a $r(L)+\epsilon _Q$ par la droite
passant par ces points, sans rencontrer de point singulier. Si on va de
$r_1$ \`a $r(l)$ par la droite $l$, on rencontre un unique point singulier qui
est $r_{\tau }$. (On peut ensuite continuer vers les points $r(l)+\epsilon _Q$,
$Q\in\P _{\sS }(M_l)$, sans rencontrer de point singulier, puisque $r(L)$ est
n\'ecessairement r\'egulier pour $S_{\tau }$).)

Pour terminer la preuve de la proposition {\bf B.4}, il reste donc \`a
v\'erifier qu'en traversant $r_{\tau }$, on trouve une donn\'ee de r\'esidu
$\Res _A^{P'}$ qui v\'erifie $(\Res _A^{P'}\mu )(\tau )\ne 0$. Suite \`a {\bf
B.7} et ce qui pr\'ec\'edait, on a $\Res _A^{P'}=\pm \Res _{\tau }\Res
_{A_1}^{P'}$, o\`u $\Res _{\tau }$ d\'esigne l'op\'erateur qui associe \`a une
fonction m\'eromorphe sur la droite affine $A_1$ dans $\o $ son r\'esidu en
$\tau $. En notant $\mu _{M_{A_1}}$ la fonction m\'eromorphe $\mu (\mu
^{M_{A_1}})^{-1}$ d\'efinie sur $A_1$, on trouve par (8.1.2) $(\Res
_{A_1}^{P'}\mu )(\sigma )=\mu _{M_{A_1}}(\sigma ) (\Res _{A_1}^{P'}\mu
^{M_{A_1}})(\sigma )$ pour $\sigma\in A_1$. Par {\bf 4.9}, l'hypoth\`ese de
r\'ecurrence et la remarque \`a la fin de {\bf B.2}, $\Res _{A_1}^{P'}\mu
^{M_{A_1}}=\Res _{A_1}^{P'\cap M_1}\mu ^{M_{A_1}}$ est une fonction constante
non nulle sur $A_1$. Comme, par le lemme {\bf 8.1}, $(\Res _{\tau }\mu
_{M_{A_1}})\ne 0$, la proposition {\bf B.4} est prouv\'ee. \hfill {\fin 2}

\null
\null
\null
\centerline {\bf Index de notations}

\null
\null
$A_{\sreg }$, {\bf 1.4}

$[A]$, {\bf 8.2}

$\A (\S )$, {\bf 1.4}

$\A _{\mu }(\S )$, {\bf 8.1}

$[\A (\S )]$, $[\A _{\mu }(\S )]$, {\bf 8.2}

$d_A\Im (\sigma ')$, {\bf 1.6}

$d_S{\psi }$, $d_{\sigma ,\alpha }(\psi )$, {\bf 3.4}

$\deg (\sigma )$, {\bf 2.2}

$E_{\so }$, {\bf 1.3}

$E_P^G$, {\bf 1.5}

$\E_2(M',\sigma )$, {\bf 6.1}

$f_{\varphi }$, {\bf 2.3}

$f_{\xi }$, {\bf 2.2}

$H_M$, {\bf 1.2}

$H_{\alpha ,c}$, {\bf 1.4}

$\ti {h}_{\alpha }$, {\bf 3.2}

$h_{\alpha }$, {\bf 1.2}

$G$, {\bf 1.1}

$\H (S)$, $\H (L)$, {\bf 1.4}

$\H (\S )$, {\bf 1.4}

$\Im (\sigma )$, $\Im _{M'}(\sigma )$, {\bf 1.3}

$J_{P\vert P'}$, {\bf 1.5}

$L_{\sreg }$, {\bf 1.4}

$\L (\S )$, $\L (\H )$, {\bf 1.4}

$M_{\Omega }$, {\bf 4.6}

$M_A$, $M_L$, {\bf 1,4}

$M_{\alpha }$, $M^1$, {\bf 1.2}

$\o _{\sigma ,M}$, $\o _{\sigma ,\alpha }$, {\bf 1.4}

$\o _{\sigma ,M',r'}$, {\bf 1.6}

$p_{\sigma ,\alpha }$, $p_H(\lambda )$, {\bf 3.4}

$P_A$, {\bf 4.6}

$P=MU$, {\bf 1.1}

$\P _{\sS} (M_L)$, {\bf 3.10}

$\P(\o ,V)$, {\bf 1.3}

$\ol{P}$, {\bf 1.1}

$PW(G)$, {\bf 2.3}

$q$ ,{\bf 1.2}

$r(\o )$, $r(A)$, $r(L)$, {\bf 1.4}

$Res _A$, $Res _{\sD }$, {\bf 3.7}

$\Res _A^C$, {\bf 3.10}

$\R (\o , \S )$, {\bf 3.1}

$\R(\o,\S,\underline {d})$, {\bf 3.7}

$\Re (\chi )$, {\bf 1.2}

$\Re (\sigma )$, $\Re _{M'}(\sigma )$, {\bf 1.3}

$\S _A$, $\S (A)$, {\bf 1.4}

$\S _{\mu }=\S $, {\bf 4.3}

$\Stab (A)$, {\bf 6.1}

$\Stab _{M'}(\o )$, $\Stab _{M'}(\sigma )$, {\bf 1.3}

$T_M$, {\bf 1.1}

$W$, $\W (M)$, $W(M)$, {\bf 1.1}

$W_{\so }$, {\bf 4.2}

$W(M,\o )$, {\bf 1.3}, {\bf 4.7}

$W_{\so ,M'}^+$, $W _{\Delta _{\sso }}$, $W _{\Delta _{\sso },M'}$, $W_{M'}^+(M,\o )$, {\bf 4.5}

$\X (M)$, $\Xim (M)$, {\bf 1.2}

$\alpha ^*$, {\bf 1.2}

$\ti {\alpha }$, {\bf 3.2}

$\Delta _{\so }$, {\bf 4.5}

$\lambda (w)$, {\bf 1.5}

$\mu _{P\vert P'}$, $\mu $, {\bf 1.5}

$\mu _{sp,\sigma }$, $\mu _{nsp, \sigma }$, {\bf 4.2}

$\psi _{\sigma }$, {\bf 3.5}

$\varphi (P',\pi ')$, {\bf 2.4}

$\varphi _{\ti {\xi }^{M'}}$, {\bf 6.1}

$\varphi _{\so }$, {\bf 7.1}

$\xi _{\varphi _{P,\sso }^H}$, {\bf 7.1}

$\Sigma ^G(T_M)$, $\Sigma ^G(P)$, $\Sigma _{red}^G(P)$, {\bf 1.1}

$\Sigma _{sp}(\sigma )$, {\bf 4.1}

$\Sigma _{\so }$, {\bf 4.2}

$\Theta $, {\bf 1.2}

$\Theta _2$, $\Theta _2(\o )$, {\bf 2.5}

$>_P$, {\bf 1.2}

\null
\null
\vfill
\eject\centerline {\bf Index terminologique}

\null
\null
\null
\null composante connexe adjacente de $L_{\sreg }$\hfill {\bf 1.4}

\null composante connexe $(P,\S )$-standard\hfill {\bf 4.6}

\null donn\'ee de r\'esidu \hfill {\bf 3.7}

\null hyperplan affine (radiciel)\hfill {\bf 1.4}

\null orbite inertielle \hfill {\bf 1.3}

\null op\'erateur r\'esidu \hfill {\bf 3.7}

\null orbite inertielle de $\sigma $ sous $M$\hfill {\bf 1.4}

\null ordre du p\^ole (d'une fonction rationnelle sur $\o $)\hfill {\bf 3.4} 

\null rang parabolique\hfill {\bf 8.1}

\null sous-espace affine (radiciel) de $\o $ (ou de $a_M^*$)\hfill {\bf 1.4}

\null sous-espace affine r\'esiduel pour $\mu $\hfill {\bf 8.1}
 
\null sous-groupe de Levi associ\'e \`a un sous-espace affine de $\o $ (ou de $a_M^*$)\hfill{\bf 1.4}

\null sous-groupe de Levi $(P,\S )$-standard\hfill {\bf 4.5}

\null sous-groupe parabolique $(P,\S )$-standard\hfill {\bf 8.2}

\vfill
\eject
{\bf R\'ef\'erences bibliographiques:}

\null [A] \it A Paley-Wiener theorem for real reductive groups, \rm J. Arthur, Acta Math. {\bf 150},
1-89, 1983.

\null [B] \it Groupes et alg\`ebres de Lie, chapitres 4,5 et 6, \rm N. Bourbaki, Masson, 1981.

\null [BS] \it A Residue Calculus for Root Systems, \rm E.P. van den Ban and H. Schlichtkrull,
Composito Mathematica {\bf 123}, 27-72, 2000.

\null [Ca] \it Finite Groups of Lie Type, \rm R.W. Carter, Wiley Classics Library, 1993.

\null [Cs] \it Introduction to the theory of admissible representations of $p$-adic reductive groups
\rm - W. CASSELMAN, non publi\'e.

\null  [H1] \it Une formule de Plancherel pour l'alg\`ebre de Hecke d'un groupe r\'eductif p-adique -
\rm V. Heiermann, Comm. Math. Helv. {\bf 76}, 388-415, 2001.

\null  [H2] \it Une caract\'erisation semi-simple des $L^2$-paires de Lusztig - \rm V. Heiermann,
manus-crit, I.A.S. Princeton, NJ, 4p., 1999.

\null  [H3] \it Vers le spectre r\'esiduel d'un groupe $p$-adique - \rm V. Heiermann, manuscrit, 24p.,
2000.

\null  [HO1] \it Yang's system of particles and Hecke algrebras, \rm G.J. Heckman et E.M. Opdam, Ann.
Math. {\bf 145}, pp. 139-173, 1997.

\null [HO2] \it Harmonic analysis for affine Hecke algebras, \rm G. J. Heckman et E.M. Opdam, dans \it
Current  developments in mathematics, 1996 (Cambridge, MA), \rm 37--60, Int. Press, Boston, MA, 1997

\null [KL] \it  Proof of the Deligne-Langlands conjecture for Hecke algebras, \rm D. Kazhdan et G.
Lusztig, Invent. math., {\bf 87}, p.153-215, 1987.

\null  [L] \it On the functional equations satisfied by Eisenstein series,
\rm R.P. Langlands, Lecture Notes in Mathematics {\bf 544}, Springer, 1976.

\null  [M1] \it Normalisation des op\'erateurs d'entrelacement et r\'eductibilit\'e des induites de
cuspidales; le cas des groupes classiques
$p$-adiques, \rm C. Moeglin, Ann. Math. {\bf 151}. pp. 817-847, 2000.

\null  [M2] \it Sur la classification des s\'eries discr\`etes des groupes classique $p$-adiques;
param\`e-tres de Langlands et exhaustivit\'e, \rm C. Moeglin, pr\'epublication, 1999.

\null
[M3] \it Points de r\'eductibilit\'e pour les induites de cuspidales, \rm C. Moeglin,
pr\'epublication, 2001.

\null [Mo] \it Eisenstein series for reductive groupes over global function fields I:the cusp form
case; II: the general case, \rm L.E. Morris, Can. J. Math. {\bf 34}, pp. 91--168, pp. 1112--1182, 1982.

\null  [MT] \it Construction of discrete series for classical $p$-adic groups, \rm C. Moeglin et M.
Tadic, pr\'epublication, 1999.

\null  [Mu] \it The unitary dual of $p$-adic $G_2$, \rm G. Muic, Duke Math. J., {\bf 90}, pp. 465--493,
1997.

\null   [MW1] \it D\'ecomposition spectrale et s\'eries d'Eisenstein. Une paraphrase de l'\'ecriture,
\rm C. Moeglin et J.-L. Waldspurger, Progress in Mathematics {\bf 113}, Birkh\"auser, Basel, 341p.,
1994.

\null  [MW2] \it Le spectre r\'esiduel de $\GL _n$, \rm C. Moeglin et J.-L. Waldspurger, Ann. ENS t.
{\bf 22}, pp. 605-674, 1989.

\null [O1] \it A generating function for the trace of the Iwahori-Hecke algebra, \rm E. M. Opdam, dans
\it Studies in Memory of Issai Schur (Joseph, Melnikov, Rentschler eds.) \rm, Progress in Mathematics,
Vol. {\bf 210}, Birkh\"auser, 2003. (arXiv:math.RT/0101006 v1 31 Dec 2000)

\null [O2] \it On the spectral decomposition of affine Hecke algebras, \rm E. M. Opdam,
pr\'epublica-tion, arXiv:math.RT/0101007 v3 11 Oct 2002.

\null  [R] \it On the Iwahori-spherical discrete series for $p$-adic Chevalley groups; formal degrees
and $L$-packets, \rm M. Reeder, Ann. ENS t. {\bf 27}, pp. 463 - 491, 1994.

\null [S1] \it Introduction to harmonic analysis on reductive $p$-adic groups, \rm A. Silberger,
Mathematical Notes of Princeton University No. {\bf 23}, Princeton, NJ, 1979.

\null [S2] \it Special representations of reductive $p$-adic groups are not integrable, \rm A.
Silberger, Ann. of Mathematics, {\bf 111}, 571 - 587, 1980.

\null [S3] \it Discrete Series and classification of $p$-adic groups I, \rm A. Silberger, Amer. J.
Mathematics, {\bf 103}, 1241 - 1321, 1981.

 \null [Sh] \it  A proof of Langlands' conjecture on Plancherel measures; complementary series for
$p$-adic groups, \rm F. Shahidi,  Ann. of Mathematics, {\bf 132}, 273-330, 1990. 

\null  [W] \it La formule de Plancherel pour les groupes p-adiques (d'apr\`es Harish-Chandra), \rm
J.-L. Waldspurger, \`a para\^\i tre dans Journal de l'Institut de Math. de Jussieu.

\null
[Z] \it Induced Representations of reductive $\wp $-adic groups II. On irreducible
representations of $\GL (n)$, \rm A.V. Zelevinsky, Ann. ENS t. 13., p. 165-210, 1980.

\bye